\documentclass{amsart}
\usepackage{amssymb,latexsym}
\usepackage{amsfonts}

\newtheorem{thm}{Theorem}[section]
\newtheorem{lem}[thm]{Lemma}

\newtheorem{rem}[thm]{Remark}
\newtheorem{pro}[thm]{Proposition}

\newcommand{\ba}{\begin{array}}
\newcommand{\ea}{\end{array}}

\def \qed{\cqfd}

\newcommand*{\QEDB}{\hfill\ensuremath{\square}}
\def\qed{\vbox{\hrule
\hbox{\vrule\hbox to 5pt{\vbox to 8pt{\vfil}\hfil}\vrule}\hrule}}

\newcommand{\beg}{\begin{eqnarray*}}
\newcommand{\begn}{\begin{eqnarray}}
\newcommand{\en}{\end{eqnarray*}}
\newcommand{\enn}{\end{eqnarray}}

\begin{document}
\title{The conical K\"ahler-Ricci flow on Fano manifolds}
\keywords{twisted K\"ahler-Ricci flow, conical K\"ahler-Ricci flow,
conical K\"ahler-Einstein metric.
}
\author{JiaWei Liu}
\address{Jiawei Liu\\School of Mathematical Sciences\\
University of Science and Technology of China\\
Hefei, 230026, P.R. China\\} \email{liujw24@mail.ustc.edu.cn}
\author{Xi Zhang}
\address{Xi Zhang\\Key Laboratory of Wu Wen-Tsun Mathematics\\ Chinese Academy of Sciences\\School of Mathematical Sciences\\
University of Science and Technology of China\\
Hefei, 230026, P.R. China\\ } \email{mathzx@ustc.edu.cn}
\thanks{AMS Mathematics Subject Classification. 53C55,\ 32W20.}
\thanks{}

\begin{abstract} In this paper, we study the long-term behavior of the conical K\"ahler-Ricci flow on Fano manifold $M$. First, based on our work of locally uniform regularity for the twisted K\"ahler-Ricci flows, we  obtain a long-time solution to the conical K\"ahler-Ricci flow by limiting a sequence of these twisted flows. Second, we study the uniform Perelman's estimates of the twisted K\"ahler-Ricci flows. After that, we  prove that the conical K\"ahler-Ricci flow must converge to a conical K\"ahler-Einstein metric if there exists one.

\end{abstract}

\maketitle

\section{Introduction}
\setcounter{equation}{0}

Let M be a compact complex manifold with K\"ahler metric $\omega_{0}$. Finding a K\"ahler-
Einstein metric in a given K\"ahler class $[\omega_{0}]$ is an important problem in K\"ahler geometry,
that is, when $2\pi c_{1}(M) = \lambda [\omega_{0}]$, establishing whether there exists a unique K\"ahler metric $\omega\in[\omega_{0}]$, such that $Ric(\omega) = \lambda \omega$. One approach to this problem is the continuity method, see the works of T. Aubin
and S.T. Yau ( \cite{TA}, \cite{STY}). The other approach is the K\"ahler-Ricci flow, which was first used by H.D. Cao in \cite{HDC} to give a parabolic proof of the Calabi-Yau theorem. There are some interesting results on  the convergence of K\"ahler-Ricci flow,  see references:  \cite{CCZ, CT1, CT2, SD1, GP, PSSW, PS, SOT},  etc. In particular, on Fano manifold, G. Tian and X.H. Zhu (\cite{TZ}, \cite{TZ2}) proved that if there exists a K\"ahler-Einstein metric, then the K\"ahler-Ricci flow with any initial metric in the  first Chern class must converge to a K\"ahler-Einstein metric in the $C^{\infty}$-topology. The main result in this paper extends theirs to the conical K\"ahler-Ricci flow.

Let $M$ be a Fano manifold of complex dimension $n$ and $D\in |-\lambda K_{M}|$ be a smooth divisor. By saying a closed positive $(1,1)$-current $\omega\in2\pi c_{1}(M)$ is a conical K\"ahler metric with angle $2\pi \beta $ ($0<\beta \leq 1$ ) along $D$, we mean that $\omega$ is a smooth K\"ahler metric on $M\setminus D$, and near each point $p\in D$, there exists local holomorphic coordinate $(z^{1}, \cdots, z^{n})$ in a neighborhood $U$ of $p$ such that locally $D=\{z^{n}=0\}$, and $\omega$ is asymptotically equivalent to the model conic metric
$$\sqrt{-1} |z^{n}|^{2\beta -2} dz^{n}\wedge d\overline{z}^{n} +\sqrt{-1} \sum_{j=1}^{n-1} dz^{j}\wedge d\overline{z}^{j}$$
on $U$.

We call $\omega $ a conic K\"ahler-Einstein metric with conic angle $2\pi \beta $ along $D$ if it is a conic K\"ahler metric and satisfies
\begin{eqnarray}\label{2015.1.21.1}
Ric (\omega )=\mu \omega +2\pi (1-\beta ) [D]
\end{eqnarray}
on $M$, where $[D]$ is the current of integration along $D$ and $\mu=1-(1-\beta)\lambda$. Here the equation (\ref{2015.1.21.1}) is classical outside $D$ and it holds in the sense of currents globally on $M$.
There are other definitions of metrics with conical singularity (see \cite{SD2} \cite{JMR}, etc.). But for conical K\"ahler-Einstein metrics, these definitions turn out to be equivalent (see Theorem $2$ in \cite{JMR}). The conical K\"ahler-Einstein metric was studied on the Riemann surfaces by R. McOwen \cite{Mco} and M. Troyanov
\cite{Tro},  and was first considered in higher dimensions by G. Tian in \cite{T0}. The renewed interest has been sparked by S. Donaldson's project which aims to solve smooth K\"ahler-Einstein problem on Fano manifold by using conical K\"aher-Einstein metrics
as a continuity method in \cite{SD2}. Recently,  the Yau-Tian-Donaldson's conjecture has been proved by G. Tian in \cite{T1}, X.X. Chen, S. Donaldson and S. Sun in \cite{CDS1, CDS2, CDS3} respectively. The existence of conical K\"ahler-Einstein metric still has its own interest, there is by now a large body of works, see references \cite{RB, SB, CGP, EGZ, GP1, JMR, LS, SW} etc.

 In this paper, we study the following conical K\"ahler-Ricci flow
\begin{eqnarray}\label{CKRF1}\frac{\partial \omega }{\partial t}=-Ric(\omega)+\beta\omega+(1-\beta)[D]\end{eqnarray}
which starts with a conical K\"ahler metric with cone angle $2\pi \beta $ along the divisor $D$. Here we assume that smooth divisor $D\in|-K_M|$, i.e. $\lambda =1$. In fact, our argument in the following are also valid for $\lambda \geq 1$, only if the coefficient $\beta$ before $\omega $ in (\ref{CKRF1}) is replaced by $1-(1-\beta )\lambda $. By saying $\omega (t)$ ($t\in [0, +\infty)$) is a long-time solution of the above conical K\"ahler-Ricci flow, we mean that for any $t$, $\omega (t)$  is a conical K\"ahler metric with conic angle $2\pi \beta $ along $D$, it satisfies (\ref{CKRF1}) in the sense of currents globally on $M$ and can be simplified to the classical K\"ahler-Ricci flow outside $D$, i.e.
$$\frac{\partial \omega }{\partial t}=-Ric(\omega)+\beta\omega$$
on $M\setminus D$. In \cite{CW}, X.X. Chen and Y.Q. Wang introduced the strong conical K\"ahler-Ricci flow and established the short-time existence. When $n=1$, R. Mazzeo, Y. Rubinstein and N. Sesum in \cite{MRS}, H. Yin in \cite{Y1}\ \cite{Y2} did it with different function spaces.

For our research, we will combine the conical K\"ahler-Ricci flow with the twisted K\"ahler-Ricci flow. By Assuming that the K\"ahler class and the first Chern class satisfy $2\pi c_{1}(M)-\beta [\omega_{0}]=[\alpha ]\neq 0$ and then fixing a closed $(1, 1)$-form $\theta \in [\alpha ]$, J. Song and G. Tian firstly introduced the twisted K\"ahler-Einstein metric
\begin{eqnarray}
Ric (\omega )=\beta \omega +\theta
\end{eqnarray}
in \cite{SOT}, there are many subsequent work, see \cite{RB, St, XZ}.

The twisted K\"ahler-Ricci flow
 \begin{eqnarray}
 \frac{\partial \omega }{\partial t}=-Ric(\omega)+\beta\omega +\theta
 \end{eqnarray}
 was studied respectively by  the first author in \cite{JWL} \cite{LW}, T. Collins and G. Sz\'ekelyhidi in \cite{TC}. In paper \cite{T1},
G. Tian approximated the conical K\"ahler-Einstein metric by a sequence of smooth twisted K\"ahler-Einstein metrics. It is natural to use G. Tian's idea to approximate the conical K\"ahler-Ricci flow by a sequence of smooth twisted K\"ahler-Ricci flows.

Let $\omega_{0}$ be a smooth K\"ahler metric in $2\pi c_{1}(M)$, $h$ be a smooth Hermitian metric on the line bundle $-K_{M}$ with curvature $\omega_{0}$ and $s$ be the defining section of $D$. It is well known that, for small $k$, \begin{eqnarray}\label{3.27.1}\omega^{\ast}=\omega_{0}+k\sqrt{-1}\partial \overline{\partial }|s|_{h}^{2\beta}\end{eqnarray}
is a conical K\"ahler metric with cone angle $2\pi \beta $ along $D$. As in \cite{CGP}, we also denote
 \begin{eqnarray}\label{3.27.2}\omega_{\varepsilon}=\omega_{0}+\sqrt{-1} k \partial\overline{\partial}\chi(\varepsilon^{2}+|s|^{2}_{h}),\end{eqnarray}
where  \begin{eqnarray}\label{3.27.3}\chi(\varepsilon^{2}+t)=\frac{1}{\beta}\int_{0}^{t}\frac{(\varepsilon^{2}+r)^{\beta}-\varepsilon^{2\beta}}{r}dr,\end{eqnarray}
$k$ is a sufficiently small number such that $\omega_{\varepsilon}$ is a K\"ahler form for each $\varepsilon>0$. It is easy to see that $\omega_{\varepsilon }$ converges to $\omega^{\ast}$ in the sense of currents globally on $M$ and in $C_{loc}^{\infty}$ topology outside $D$. From \cite{CGP}, we know that the function $\chi(\varepsilon^{2}+t)$ is smooth for each $\varepsilon>0$, and there exist uniform constants (independent of $\varepsilon$) $C>0$ and $\gamma>0$ such that
\begin{eqnarray}\label{3.19.1}0\leq\chi(\varepsilon^{2}+t)\leq C,\end{eqnarray}
provided that $t$ belongs to a bounded interval and
\begin{eqnarray}\label{3.19.2}\omega_{\varepsilon}\geq\gamma \omega_{0}.\end{eqnarray}

Now, we consider the following  twisted K\"ahler-Ricci flow:
\begin{eqnarray}\label{GKRF1}
\begin{cases}
  \frac{\partial \omega_{\varphi_{\varepsilon}}}{\partial t}=-Ric(\omega_{\varphi_{\varepsilon}})+\beta\omega_{\varphi_{\varepsilon}}+
  (1-\beta)(\omega_{0}+\sqrt{-1}\partial\overline{\partial}\log(\varepsilon^{2}+|s|_{h}^{2})),\\
  \\
  \omega_{\varphi_{\varepsilon}}|_{t=0}=\omega_{\varepsilon}\\
  \end{cases}
  \end{eqnarray}
where $\omega_{\varphi_{\varepsilon}}=\omega_{\varepsilon}+\sqrt{-1}\partial\overline{\partial}\varphi_{\varepsilon}$. We can see that $(1-\beta)(\omega_{0}+\sqrt{-1}\partial\overline{\partial}\log(\varepsilon^{2}+|s|_{h}^{2}))$ is a smooth closed semi-positive (1,1)-form.
Since the twisted K\"ahler-Ricci flow preserves the K\"ahler class, we can write this flow as the parabolic Monge-Amp\'ere equation on potentials:
\begin{eqnarray}\label{CMAF1}
\begin{cases}
  \frac{\partial \varphi_{\varepsilon}}{\partial t}=\log\frac{\omega_{\varphi_{\varepsilon}}^{n}}{\omega_{0}^{n}}+F_{0}+\beta(k\chi+\varphi_{\varepsilon})
  +\log(\varepsilon^{2}+|s|_{h}^{2})^{1-\beta},\\
  \\
  \varphi_{\varepsilon}|_{t=0}=c_{\varepsilon0}\\
  \end{cases}
\end{eqnarray}
where the constant $c_{\varepsilon 0}$ (its representation will be given in section $5$) is uniformly bounded for $\varepsilon$, $F_{0}$ satisfies $-Ric(\omega_{0})+\omega_{0}=\sqrt{-1}\partial\overline{\partial}F_{0}$ and $\frac{1}{V}\int_{M}e^{-F_{0}}dV_{0}=1$, and $\chi$ denotes the function $\chi(\varepsilon^{2}+|s|^{2}_{h})$. Sometimes, we will rewrite the flow (\ref{CMAF1}) as follows:
\begin{eqnarray}\label{CMAF2}
\begin{cases}
  \frac{\partial \varphi_{\varepsilon}}{\partial t}=\log\frac{\omega_{\varphi_{\varepsilon}}^{n}}{\omega_{\varepsilon}^{n}}+F_{\varepsilon}+\beta(k\chi+\varphi_{\varepsilon}),\\
  \varphi_{\varepsilon}|_{t=0}=c_{\varepsilon0}\\
  \end{cases}
\end{eqnarray}
where $F_{\varepsilon}=F_{0}+\log(\frac{\omega_{\varepsilon}^{n}}{\omega_{0}^{n}}\cdot(\varepsilon^{2}+|s|_{h}^{2})^{1-\beta})$.

In our paper, on the basis of proving locally uniform estimates for equation (\ref{CMAF2}), we obtain a long-time solution to the conical K\"ahler-Ricci flow (\ref{CKRF1}) on Fano manifolds by limiting a sequence of the twisted K\"ahler-Ricci flows (\ref{GKRF1}) as $\varepsilon \rightarrow 0$. For any $\beta\in(0,1)$, we prove uniform Perelman's estimates (when $t\geq1$) and uniform Sobolev inequalities (when $t\geq0$) along the twisted K\"ahler-Ricci flows (\ref{GKRF1}). Here the uniformity means that the constants in the estimates and inequalities are independent of $\varepsilon$ and $t$. Using these estimates, we prove that the conical K\"ahler-Ricci flow (\ref{CKRF1}) must converge to a conical K\"ahler-Einstein metric if there exists one and the convergence is in $C_{loc}^{\infty}$ topology outside the divisor $D$ and globally in the sense of currents on $M$. In fact, we prove the following theorem:

\begin{thm}
\label{1.1}
Let $\omega_{\varphi_{\varepsilon }} $ be a long-time solution of the twisted K\"ahler-Ricci flow (\ref{GKRF1}), then there must exist a sequence $\varepsilon_{i}\rightarrow0$ such that $\omega_{\varphi_{\varepsilon_{i} }} $ converges to a solution of the conical K\"ahler-Ricci flow
\begin{eqnarray}\label{CKRF2}
\begin{cases}
  \frac{\partial \omega_{\varphi}}{\partial t}=-Ric(\omega_{\varphi})+\beta\omega_{\varphi}+(1-\beta)[D],\\
  \omega_{\varphi}|_{t=0}=\omega^*\\
  \end{cases}
  \end{eqnarray}
where $\omega^*=\omega_{0}+k\sqrt{-1}\partial\bar{\partial}|s|_{h}^{2\beta}$. The convergence is in $C_{loc}^{\infty}$ topology outside the divisor $D$ and in the sense of currents on $M\times[0,+\infty)$. Furthermore, the potential $\varphi(t)$ is H$\ddot{o}$lder continuous with respect to the smooth metric $\omega_{0}$ on $M$.

Moreover, if there exists a conical K\"ahler-Einstein metric with cone angle $2\pi \beta $ ($0<\beta<1$) along $D$, then the long-time solution $\omega_{\varphi }(\cdot, t)$ must converge to a conical K\"ahler-Einstein metric in $C_{loc}^{\infty}$ topology outside the divisor $D$ and globally in the sense of currents on $M$.
\end{thm}

\medskip

In \cite{YQW}, Y.Q. Wang also considered the long-time existence of a weak conical K\"ahler-Ricci flow $(\ref{CKRF1})$ by using the limiting method. The difference is that we further study the local uniform  higher order estimates of the twisted K\"ahler-Ricci flows $(\ref{CMAF1})$. In fact, we can get local uniform $C^{\infty}$ estimates outside the divisor $D$ on any finite time interval $[0, T]$(see Proposition \ref{2.2}). Our argument is based on elliptic estimate which is superior to the parabolic Schauder estimates, because the latter can only provide us with a local uniform $C^{\infty}$ estimates on $\overline{B_{r}}\times [\delta,T]$ for some $\delta>0$, and the fact that these estimates depend on $\delta$.

\medskip

 In \cite{CW1}, X.X. Chen and Y.Q. Wang  proved the existence of long-time solution for the strong conical K\"ahler-Ricci flow, and obtained the convergence result when $\mu =1-(1-\beta )\lambda \leq 0$, i.e. the twisted first Chern class is negative or zero. In this paper, we consider the convergence with positive twisted first Chern class. Recently,  R. Berman, S. Boucksom, P. Eyssidieux, V. Guedj and A. Zeriahi \cite{BBEGZ} studied the convergence of the K\"ahler-Ricci flow on $\mathbb{Q}$-Fano variety with log terminal singularities. They proved that if the Mabuchi functional is proper, then the solution of K\"ahler-Ricci flow converges to the unique K\"ahler-Einstein metric in some weak sense.   Although their work is independent of Perelman's estimates, their weakly convergence can't yield the convergence in $C^\infty$ topology even if $M$ is non-singular. Here, our main goal is to prove the local $C^\infty$ convergence by obtaining the uniform Perelman's estimates along flows (\ref{GKRF1}). By the arguments in \cite{JWL} or \cite{NSGT}, we know that these estimates mainly depend on the bound of the initial twisted scalar curvature $R(g_{\varepsilon}(0))-tr_{g_{\varepsilon}(0)}\theta_{\varepsilon}$. But it may not be bounded uniformly when $\beta\in(\frac{1}{2}, 1)$. In order to overcome this difficulty, we need the following key observation (Proposition \ref{1.8.1}) that  the twisted scalar curvature $R(g_{\varepsilon}(t))-tr_{g_{\varepsilon}(t)}\theta_{\varepsilon}$  is bounded uniformly from below along the flows $(\ref{GKRF1})$ when $t\geq 1$. Using this observation, we can get the uniform Perelman's estimates on $[1, +\infty )$, which is enough for us to study the convergence of the conical K\"ahler-Ricci flow. For details, one can see section $4$.

 \medskip

The paper is organized as follows. In section $2$, we prove the uniform Laplacian estimate and local $C^{\infty}$ estimates for the twisted K\"ahler-Ricci flow (\ref{GKRF1}). Then, in section $3$,  we get a long-time solution to the conical K\"ahler-Ricci flow (\ref{CKRF2}) by limiting a sequence of twisted K\"ahler-Ricci flows. In section $4$, for any $\beta\in(0,1)$, we obtain uniform Perelman's estimates along twisted K\"ahler-Ricci flows $(\ref{GKRF1})$ when $t\geq1$. Making use of these estimates, we choose a suitable initial value $\varphi_\varepsilon(0)$, and then obtain uniform $C^{0}$ estimates for the metric potentials with the uniform properness of the twisted Mabuchi $\mathcal{K}$-energy functional in section $5$. At the last section, we first give a remark to C.J. Yao's paper \cite{CJY} which provides an alternative proof of S. Donaldson's openness theorem. Next, we show that the properness of  Log Mabuchi $\mathcal{K}$-energy functional $\mathcal{M}_{ \omega_{0},\ (1-\beta)D }$ implies the uniform properness of the twisted Mabuchi $\mathcal{K}$-energy functional $\mathcal{M}_{ \omega_{0},\ \theta_{\varepsilon} }$. Then we prove that the conical K\"ahler-Ricci flow (\ref{CKRF2}) must converge to a conical K\"ahler-Einstein metric in $C_{loc}^{\infty}$ topology outside the divisor $D$ and in the sense of currents on $M$ if there exists one.\\
\medskip

{\bf  Acknowledgement:} We would like to thank Professor J.Y. Li and Professor X.H. Zhu for their useful conversations and suggestions. We are also grateful to the referees for their careful reading and valuable suggestions. The authors are supported in part by NSF in China No.11131007 and the Hundred Talents Program of CAS.

\section{The local estimates for the twisted K\"ahler-Ricci flows}
\setcounter{equation}{0}

In this section, we will give the uniform Laplacian estimate and local higher order estimates for the parabolic Monge-Amp\'ere equation (\ref{CMAF2}). In the following sections, by saying a uniform constant, we mean that it is independent of $\varepsilon$ and $t$. We shall use the letter $C$ for a uniform constant which may differ from line to line. We first follow Guenancia-Paun's argument ( in \cite{GP1}) to obtain the Laplacian estimate, we have:

\begin{pro}\label{2.1} Let $\varphi_{\varepsilon}$ be a solution of equation $(\ref{CMAF2})$.
Assume that there exists a uniform constant $C>0$ such that

(1) $\sup\limits_{M\times[0,T]}|\varphi_{\varepsilon}|\leq C$;

(2) $\sup\limits_{M\times[0,T]}|\dot{\varphi}_{\varepsilon}|\leq C$.\\

Then there exists a uniform constant $A$ only depending on $\omega_{0}$, $n$, $\beta$ and $C$, such that
$$A^{-1}\omega_{\varepsilon}\leq \omega_{\varepsilon}+\sqrt{-1}\partial\overline{\partial}\varphi_{\varepsilon}\leq A \omega_{\varepsilon}$$
on $M\times[0,T]$.
\end{pro}

We notice that the estimates are independent of time $T$, so the above result holds also for time interval $[0,+\infty)$. In local coordinates,
$$\omega=\sqrt{-1}g_{i\bar{j}}dz^{i}\wedge d\bar{z}^{j},$$
with
\begin{eqnarray}\label{2015.1}R_{i\bar{j}k\bar{l}}=-\frac{\partial^{2}g_{i\bar{j}}}{\partial z^{k}\partial \bar{z}^{l}}+g^{r\bar{s}}\frac{\partial g_{i\bar{s}}}{\partial z^{k}}\frac{\partial g_{r\bar{j}}}{\partial \bar{z}^{l}},\end{eqnarray}
as its corresponding components of the curvature tensor, and the Ricci curvature
\begin{eqnarray}\label{2015.2}R_{i\bar{j}}=g^{k\bar{l}}R_{i\bar{j}k\bar{l}}.\end{eqnarray}

Let's first recall the appropriate coordinate system ( see Lemma $4.1$ in \cite{CGP}).

\begin{lem}\label{1.8.11} Let $(L,h)$ be the hermitian line bundle associated to a smooth divisor $D$, and $s$ be a section of $L$ such that
$$D:=\{s=0\}.$$
Let $p_0\in D$, then there exists a constant $C>0$ and an open set $\Omega\subset M$ centered at $p_0$, such that for any point $p\in\Omega$ there exists a coordinate system $z=(z^1,\cdots,z^n)$ and a trivialization $\eta$ for $L$ such that:

$(1)$ $D\bigcap\Omega=\{z^n=0\}$;

$(2)$ With respect to the trivialization $\eta$, the metric $h$ has the weight $\varphi$, such that
\begin{eqnarray}\label{3.22.1}\varphi(p)=0,\ \  d\varphi(p)=0,\ \  |\frac{\partial^{|\alpha|+|\beta|}\varphi}{\partial z^\alpha\partial \bar{z}^{\beta}}(p)|\leq C_{\alpha,\beta}
\end{eqnarray}
for some constant $C_{\alpha,\beta}$ depending only on the multi indexes $\alpha, \beta$.
\end{lem}

{\bf Proof of Proposition \ref{2.1}:}\ \ We let $\varphi_{\varepsilon}$ evolve along the parabolic Monge-Amp\'ere equation $(\ref{CMAF2})$. By direct computation, we have
\begin{eqnarray}\label{3.19.3}&\ &(\frac{d}{dt}-\triangle_{\omega_{\varphi_{\varepsilon}}})\log tr_{\omega_{\varepsilon}}\omega_{\varphi_{\varepsilon}}=\frac{1}{tr_{\omega_{\varepsilon}}\omega_{\varphi_{\varepsilon}}}
(\triangle_{\omega_{\varepsilon}}(\dot{\varphi_{\varepsilon}}-\log\frac{\omega_{\varphi_{\varepsilon}}^{n}}{\omega_{\varepsilon}^{n}})+R_{\omega_{\varepsilon}})
\\ \nonumber
&-&\frac{1}{tr_{\omega_{\varepsilon}}\omega_{\varphi_{\varepsilon}}}(g_{\varphi_{\varepsilon}}^{p\bar{q}}g_{\varphi_{\varepsilon} m\bar{j}}R_{\omega_{\varepsilon}p\bar{q}}^{\bar{m}j})+
\{\frac{g_{\varphi_{\varepsilon}}^{\delta\bar{k}}\partial_{\delta}tr_{\omega_{\varepsilon}}\omega_{\varphi_{\varepsilon}}\partial_{\bar{k}}tr_{\omega_{\varepsilon}}\omega_{\varphi_{\varepsilon}} }{(tr_{\omega_{\varepsilon}}\omega_{\varphi_{\varepsilon}})^{2}}-\frac{g_{\varepsilon}^{\gamma\bar{s}}\varphi^{\ \ \ t}_{\varepsilon\gamma\ p}\varphi_{\varepsilon\bar{s}t}^{\ \ \ p}}{tr_{\omega_{\varepsilon}}\omega_{\varphi_{\varepsilon}}}\}.
\end{eqnarray}

Then we choose a local coordinate system $w=(w^1,\ldots, w^n)$, to make $(g_{\varepsilon i\bar{i}})$ be identity and $(g_{\varphi_{\varepsilon} i\bar{i}})$ be a diagonal matrix. Since $(g_{\varphi_{\varepsilon} i\bar{i}})$ is positive definite, we have $g_{\varphi_{\varepsilon} i\bar{i}}=1+\varphi_{\varepsilon i\bar{i}}>0$. It was shown by T. Aubin \cite{TA} and S.T. Yau \cite{STY} that

\begin{eqnarray}\label{3.19.4}\frac{g_{\varphi_{\varepsilon}}^{\delta\bar{k}}\partial_{\delta}tr_{\omega_{\varepsilon}}\omega_{\varphi_{\varepsilon}} \partial_{\bar{k}}tr_{\omega_{\varepsilon}}\omega_{\varphi_{\varepsilon}} }{(tr_{\omega_{\varepsilon}}\omega_{\varphi_{\varepsilon}})^{2}}-\frac{g_{\varepsilon}^{\gamma\bar{s}}\varphi^{\ \ \ t}_{\varepsilon\gamma\ p}\varphi_{\varepsilon\bar{s}t}^{\ \ \ p}}{tr_{\omega_{\varepsilon}}\omega_{\varphi_{\varepsilon}}}\leq0.
\end{eqnarray}
On the other hand,
\begin{eqnarray}\label{3.19.5}n=tr_{\omega_{\varepsilon}}\omega_{0}+k\triangle_{\omega_{\varepsilon}}\chi\geq k\triangle_{\omega_{\varepsilon}}\chi.
\end{eqnarray}
By substituting $(\ref{CMAF2})$, $(\ref{3.19.4})$ and $(\ref{3.19.5})$ into $(\ref{3.19.3})$, we have
\begin{eqnarray}\label{15.15.15}\nonumber(\frac{d}{dt}-\triangle_{\omega_{\varphi_{\varepsilon}}})\log tr_{\omega_{\varepsilon}}\omega_{\varphi_{\varepsilon}}&\leq&-\frac{1}{tr_{\omega_{\varepsilon}}\omega_{\varphi_{\varepsilon}}}\sum_{i\leq j}(\frac{1+\varphi_{\varepsilon i\overline{i}}}{1+\varphi_{\varepsilon j\overline{j}}}+\frac{1+\varphi_{\varepsilon j\overline{j}}}{1+\varphi_{\varepsilon i\overline{i}}}-2)R_{\omega_{\varepsilon}i\overline{i}j\overline{j}}(w)\\\nonumber
&\ &+\frac{1}{tr_{\omega_{\varepsilon}}\omega_{\varphi_{\varepsilon}}}(\triangle_{\omega_{\varepsilon}}(F_{\varepsilon}+\beta\varphi+k\beta\chi))\\\nonumber
&=&-\frac{1}{tr_{\omega_{\varepsilon}}\omega_{\varphi_{\varepsilon}}}\sum_{i\leq j}(\frac{1+\varphi_{\varepsilon i\overline{i}}}{1+\varphi_{\varepsilon j\overline{j}}}+\frac{1+\varphi_{\varepsilon j\overline{j}}}{1+\varphi_{\varepsilon i\overline{i}}}-2)R_{\omega_{\varepsilon}i\overline{i}j\overline{j}}(w)\\\nonumber
&\ &+\frac{1}{tr_{\omega_{\varepsilon}}\omega_{\varphi_{\varepsilon}}}(\triangle_{\omega_{\varepsilon}}(F_{\varepsilon}+k\beta\chi))
+\beta\frac{\sum\limits_{i}\varphi_{\varepsilon i\bar{i}}}{\sum\limits_{i}(1+\varphi_{\varepsilon i\bar{i}})}\\\nonumber
&\leq&-\frac{1}{tr_{\omega_{\varepsilon}}\omega_{\varphi_{\varepsilon}}}\sum_{i\leq j}(\frac{1+\varphi_{\varepsilon i\overline{i}}}{1+\varphi_{\varepsilon j\overline{j}}}+\frac{1+\varphi_{\varepsilon j\overline{j}}}{1+\varphi_{\varepsilon i\overline{i}}}-2)R_{\omega_{\varepsilon}i\overline{i}j\overline{j}}(w)\\\nonumber
&\ &+\frac{1}{tr_{\omega_{\varepsilon}}\omega_{\varphi_{\varepsilon}}}(\triangle_{\omega_{\varepsilon}}F_{\varepsilon})+\frac{\beta n}{tr_{\omega_{\varepsilon}}\omega_{\varphi_{\varepsilon}}}
+\beta-\beta\frac{n}{\sum\limits_{i}(1+\varphi_{\varepsilon i\bar{i}})}\\\nonumber
\end{eqnarray}
\begin{eqnarray}
\nonumber&\leq&-\frac{1}{tr_{\omega_{\varepsilon}}\omega_{\varphi_{\varepsilon}}}\sum_{i\leq j}(\frac{1+\varphi_{\varepsilon i\overline{i}}}{1+\varphi_{\varepsilon j\overline{j}}}+\frac{1+\varphi_{\varepsilon j\overline{j}}}{1+\varphi_{\varepsilon i\overline{i}}}-2)R_{\omega_{\varepsilon}i\overline{i}j\overline{j}}(w)\\
&\ &+\frac{1}{tr_{\omega_{\varepsilon}}\omega_{\varphi_{\varepsilon}}}(\triangle_{\omega_{\varepsilon}}F_{\varepsilon})
+\frac{\beta n}{tr_{\omega_{\varepsilon}}\omega_{\varphi_{\varepsilon}}}+\beta.
\end{eqnarray}

First of all, we deal with the term $\triangle_{\omega_{\varepsilon}}F_0$. We know that there exists a uniform constant $C$ such that $\sqrt{-1}\partial\bar{\partial}F_0\geq-C\omega_{0}$. Then by $(\ref{3.19.2})$, we have
$$0\leq tr_{\omega_{\varepsilon}}(\sqrt{-1}\partial\bar{\partial}F_0+C\omega_{0})\leq\gamma^{-1}(Cn+\triangle_{\omega_{0}}F_0)$$
and thus
\begin{eqnarray}\label{3.19.6}-C\gamma^{-1}\leq\triangle_{\omega_{\varepsilon}}F_0\leq\gamma^{-1}(Cn+\triangle_{\omega_{0}}F_0),
\end{eqnarray}
which shows that $\triangle_{\omega_{\varepsilon}}F_0$ is uniformly bounded.

Now we deal with the terms $R_{\omega_{\varepsilon}i\overline{i}j\overline{j}}(w)$ following the argument by H. Guenancia and M. P$\breve{a}$un in \cite{GP1}. To reader's convenience, we give the proof here briefly. At point $p$, we choose the H. Guenancia and M. P$\breve{a}$un's coordinate system $z=(z^1,\ldots,z^n)$ in Lemma \ref{1.8.11}, then the coefficients of holomorphic bisectional curvature change as follows
\begin{eqnarray}\label{15.1.1}R_{\omega_\varepsilon i\bar{i}l\bar{l}}(w)=R_{\omega_\varepsilon p\bar{q}r\bar{s}}(z)\frac{\partial z^p}{\partial w^i}\overline{\frac{\partial z^q}{\partial w^i}}\frac{\partial z^r}{\partial w^l}\overline{\frac{\partial z^s}{\partial w^l}}.
\end{eqnarray}
At the point $p$, we have
\begin{eqnarray}\label{15.1.2}\omega_\varepsilon\geq C\sqrt{-1}\frac{dz^n\wedge d\bar{z}^n}{(\varepsilon^2+|z^n|^2)^{1-\beta}}
\end{eqnarray}
for some uniform constant $C$ independent of $\varepsilon$ and the point $p$. Since $(\frac{\partial}{\partial w^k})$ is unit at the point $p$ with respect to the metric $\omega_\varepsilon$, we have the estimate
\begin{eqnarray}\label{15.1.3}|\frac{\partial z^n}{\partial w^k}|^2\leq C(\varepsilon^2+|z^n|^2)^{1-\beta}.
\end{eqnarray}

From the computation in \cite{GP1}, at the point $p$, we have
\begin{eqnarray}\label{15.1.11.1}R_{\omega_{\varepsilon}i\overline{i}j\overline{j}}(w)\geq-C_1((\uppercase\expandafter{\romannumeral1})+
(\uppercase\expandafter{\romannumeral2})+(\uppercase\expandafter{\romannumeral3}))-C_2,
\end{eqnarray}
where $C_1$ and $C_2$ are uniform constants independent of $\varepsilon$ and the point $p$,
\begin{eqnarray*}
(\uppercase\expandafter{\romannumeral1})&=&\sum_{i,j}\frac{1}{(\varepsilon^2+|z^n|^2)^{\frac{1}{2}}}|\frac{\partial z^n}{\partial w^i}||\frac{\partial z^n}{\partial w^j}|,\\
(\uppercase\expandafter{\romannumeral2})&=&\sum_{i,j}\frac{1}{(\varepsilon^2+|z^n|^2)^{\frac{1}{2}}}|\frac{\partial z^n}{\partial w^i}|^2|\frac{\partial z^n}{\partial w^j}|,\\
(\uppercase\expandafter{\romannumeral3})&=&\sum_{i,j}\frac{1}{\varepsilon^2+|z^n|^2}|\frac{\partial z^n}{\partial w^i}|^2|\frac{\partial z^n}{\partial w^j}|^2.
\end{eqnarray*}

Now we need to deal with $(\uppercase\expandafter{\romannumeral1})$, $(\uppercase\expandafter{\romannumeral2})$ and $(\uppercase\expandafter{\romannumeral3})$.
Take the coefficients of $\frac{1+\varphi_{\varepsilon i\overline{i}}}{1+\varphi_{\varepsilon j\overline{j}}}$ in (\ref{15.15.15}) as an example. By (\ref{15.1.3}), its coefficients can be dominated as follows:
\begin{eqnarray}\label{15.1.11.2}\nonumber(\uppercase\expandafter{\romannumeral1})&\leq&\sum_{j}\frac{1}{(\varepsilon^2+|z^n|^2)^{\frac{\beta}{2}}}|\frac{\partial z^n}{\partial w^j}|\leq\sum_{j}\frac{C}{(\varepsilon^2+|z^n|^2)^{\beta}}|\frac{\partial z^n}{\partial w^j}|^2+C,\\
(\uppercase\expandafter{\romannumeral2})&\leq&\sum_{j}\frac{1}{(\varepsilon^2+|z^n|^2)^{\beta-\frac{1}{2}}}|\frac{\partial z^n}{\partial w^j}|\leq\sum_{j}\frac{C}{(\varepsilon^2+|z^n|^2)^{2\beta-1}}|\frac{\partial z^n}{\partial w^j}|^2+C,\\\nonumber
\end{eqnarray}
\begin{eqnarray}\nonumber
(\uppercase\expandafter{\romannumeral3})&\leq&\sum_{j}\frac{C}{(\varepsilon^2+|z^n|^2)^{\beta}}|\frac{\partial z^n}{\partial w^j}|^2,
\end{eqnarray}
where constant $C$ is independent of $\varepsilon$ and the point $p$. At the same time, $|\triangle_{\omega_{\varepsilon}}\log(\frac{\omega_{\varepsilon}^{n}}{\omega_{0}^{n}}\cdot(\varepsilon^{2}+|s|_{h}^{2})^{1-\beta})|$ can be dominated by
\begin{eqnarray}\label{15.1.11.3}\sum\limits_j\frac{C}{(\varepsilon^2+|z^n|^2)^{\tilde{\beta}}}|\frac{\partial z^n}{\partial w^j}|^2,\end{eqnarray}
where $\tilde{\beta}=\max(\beta, 1-\beta)$ and $C$ is independent of $\varepsilon$ (see section 5.2 in \cite{GP1}).
We denote $\Psi_{\varepsilon,\rho}=\tilde{C}\chi_{\rho}(\varepsilon^{2}+|s|_{h}^{2})$, where
\begin{eqnarray}\chi_{\rho}(\varepsilon^{2}+|s|_{h}^{2})=\frac{1}{\rho}\int_{0}^{|s|_{h}^{2}}\frac{(\varepsilon^{2}+r)^{\rho}-\varepsilon^{2\rho}}{r}dr.\end{eqnarray}
and $\tilde{C}$ will be determined later. The choice of the function $\chi_{\rho}$ above is motivated by the following equality:
\begin{eqnarray}\label{15.1.4}\sqrt{-1}\partial\bar{\partial}\chi_\rho(\varepsilon^2+|s|_{h}^2)=\sqrt{-1}\frac{\langle D's, D's\rangle}{(\varepsilon^2+|s|_{h}^2)^{1-\rho}}-\frac{1}{\beta}((\varepsilon^2+|s|_{h}^2)^{\rho}-\varepsilon^2)\omega_0.
\end{eqnarray}
Corresponding to $\omega_{\varphi_\varepsilon}$, we evaluate the Laplacian of the function $\Psi_{\varepsilon,\rho}$ by using the $(w)$-coordinates, then
\begin{eqnarray}\label{15.1.5}&\ &\\\nonumber
&\ &tr_{\omega_{\varepsilon}}\omega_{\varphi_{\varepsilon}}\Delta_{\omega_{\varphi_\varepsilon}}\Psi_{\varepsilon,\rho}\geq
-\tilde{C}tr_{\omega_{\varepsilon}}\omega_{\varphi_{\varepsilon}}tr_{\omega_{\varphi_\varepsilon}}\omega_\varepsilon+\tilde{C}\sum_{j=1}^{n}(\frac{1}{(\varepsilon^2+|s|_{h}^2)^{1-\rho}}|\frac{\partial z^n}{\partial w^j}|^2\frac{tr_{\omega_{\varepsilon}}\omega_{\varphi_{\varepsilon}}}{1+\varphi_{\varepsilon j\bar{j}}}).
\end{eqnarray}
Hence, after taking sufficiently large uniform constants $\tilde{C}$ and $1-\rho>\tilde{\beta}$, we can cancel the terms in (\ref{15.1.11.2}) and (\ref{15.1.11.3}) by (\ref{15.1.5}). In fact, we have
\begin{eqnarray}\label{3.19.7}&\ &-\sum_{i\leq j}(\frac{1+\varphi_{\varepsilon i\overline{i}}}{1+\varphi_{\varepsilon j\overline{j}}}+\frac{1+\varphi_{\varepsilon j\overline{j}}}{1+\varphi_{\varepsilon i\overline{i}}}-2)R_{\omega_{\varepsilon}i\overline{i}j\overline{j}}(w)-tr_{\omega_{\varepsilon}}\omega_{\varphi_{\varepsilon}}\triangle_{\omega_{\varphi_{\varepsilon}}}\Psi_{\varepsilon,\rho}\\ \nonumber
&\ &+\triangle_{\omega_{\varepsilon}}\log(\frac{\omega_{\varepsilon}^{n}}{\omega_{0}^{n}}\cdot(\varepsilon^{2}+|s|_{h}^{2})^{1-\beta})
\\ \nonumber
&\leq&C\sum_{i\leq j}(\frac{1+\varphi_{\varepsilon i\overline{i}}}{1+\varphi_{\varepsilon j\overline{j}}}+\frac{1+\varphi_{\varepsilon j\overline{j}}}{1+\varphi_{\varepsilon i\overline{i}}})+Ctr_{\omega_{\varphi_{\varepsilon}}}\omega_{\varepsilon}\cdot tr_{\omega_{\varepsilon}}\omega_{\varphi_{\varepsilon}}+C
\end{eqnarray}
for some uniform constant $C$. Combining $(\ref{3.19.6})$ with $(\ref{3.19.7})$, we have
\begin{eqnarray*}(\frac{d}{dt}-\triangle_{\omega_{\varphi_{\varepsilon}}})(\log tr_{\omega_{\varepsilon}}\omega_{\varphi_{\varepsilon}}+\Psi_{\varepsilon,\rho})&\leq&\frac{C}{tr_{\omega_{\varepsilon}}\omega_{\varphi_{\varepsilon}}}\sum_{i\leq j}(\frac{1+\varphi_{\varepsilon i\overline{i}}}{1+\varphi_{\varepsilon j\overline{j}}}+\frac{1+\varphi_{\varepsilon j\overline{j}}}{1+\varphi_{\varepsilon i\overline{i}}})+\frac{C}{tr_{\omega_{\varepsilon}}\omega_{\varphi_{\varepsilon}}}\\
&\ &+Ctr_{\omega_{\varphi_{\varepsilon}}}\omega_{\varepsilon}+C\\
&=&\frac{C}{tr_{\omega_{\varepsilon}}\omega_{\varphi_{\varepsilon}}}[(\sum_{i}\frac{1}{1+\varphi_{\varepsilon i\overline{i}}})(\sum_{j}(1+\varphi_{\varepsilon j\overline{j}}))+n]\\
&\ &+Ctr_{\omega_{\varphi_{\varepsilon}}}\omega_{\varepsilon}+\frac{C}{tr_{\omega_{\varepsilon}}\omega_{\varphi_{\varepsilon}}}+C\\
&\leq&Ctr_{\omega_{\varphi_{\varepsilon}}}\omega_{\varepsilon}+C.
\end{eqnarray*}
Here we use the fact $n\leq tr_{\omega_{\varphi_{\varepsilon}}}\omega_{\varepsilon}\cdot tr_{\omega_{\varepsilon}}\omega_{\varphi_{\varepsilon}}$ in the last inequality. Then we have
\begin{eqnarray*}(\frac{d}{dt}-\triangle_{\omega_{\varphi_{\varepsilon}}})(\log tr_{\omega_{\varepsilon}}\omega_{\varphi_{\varepsilon}}+\Psi_{\varepsilon,\rho}-B\varphi_{\varepsilon})&\leq&C tr_{\omega_{\varphi_{\varepsilon}}}\omega-B\dot{\varphi_{\varepsilon}}+B\triangle_{\omega_{\varphi_{\varepsilon}}}\varphi_{\varepsilon}+C\\
&\leq&-tr_{\omega_{\varphi_{\varepsilon}}}\omega_{\varepsilon}+C,
\end{eqnarray*}
where $B=C+1$.

By the maximum principle, at the maximum point p of $\log tr_{\omega_{\varepsilon}}\omega_{\varphi_{\varepsilon}}+\Psi_{\varepsilon,\rho}-B\varphi_{\varepsilon}$, we have
$$tr_{\omega_{\varphi_{\varepsilon}}}\omega_{\varepsilon}(p)\leq C.$$
Connecting with the fact that $F_\varepsilon$ is uniformly bounded (see $(25)$ in \cite{CGP}), we obtain $$tr_{\omega_{\varepsilon}}\omega_{\varphi_{\varepsilon}}(p)\leq\frac{1}{(n-1)!}(tr_{\omega_{\varphi_{\varepsilon}}}
\omega_{\varepsilon})^{n-1}(p)\frac{\omega_{\varphi_{\varepsilon}}^{n}}{\omega_{\varepsilon}^{n}}(p)\leq C\exp(\dot{\varphi_{\varepsilon}}-F_{\varepsilon}-\beta\varphi_{\varepsilon}-k\beta\chi)(p)\leq C.$$
Hence we have
\begin{eqnarray}\label{2015.1.24.1}tr_{\omega_{\varepsilon}}\omega_{\varphi_{\varepsilon}}\leq \exp(C+B\varphi_{\varepsilon}-B\varphi_{\varepsilon}(p))\leq C.\end{eqnarray}

On the other hand, considering the assumptions on $\varphi_{\varepsilon}$ and $\dot{\varphi_{\varepsilon}}$, we can conclude that
\begin{eqnarray}\label{2015.1.24.2}C^{-1}\leq\frac{(\omega_{\varepsilon}+\sqrt{-1}\partial\bar{\partial}\varphi_{\varepsilon})^{n}}{\omega_{\varepsilon}^{n}}=\exp(\dot{\varphi_{\varepsilon}}
-F_{\varepsilon}-\beta\varphi_{\varepsilon}-k\beta\chi)\leq C.\end{eqnarray}
By (\ref{2015.1.24.1}) and (\ref{2015.1.24.2}), there exists a uniform constant $A$ such that
\begin{eqnarray}\label{EQU}A^{-1}\omega_{\varepsilon}\leq \omega_{\varepsilon}+\sqrt{-1}\partial\overline{\partial}\varphi_{\varepsilon}\leq A \omega_{\varepsilon}\end{eqnarray}
for any $\varepsilon$ and $t$.\QEDB

\medskip

Now we consider the local Calabi's $C^{3}$ estimate and higher order estimates to the twisted K\"ahler-Ricci flow:
\begin{eqnarray}\label{GKRF2}\frac{\partial \omega_{\varphi}}{\partial t}=-Ric(\omega_{\varphi})+\lambda \omega_{\varphi}+\theta,\end{eqnarray}
where $c_{1}(M)=\lambda[\omega_{0}]+[\theta]$, $\omega_{\varphi}=\omega_{0}+\sqrt{-1}\partial\overline{\partial}\varphi$ and $\theta$ is a smooth semi-positive closed $(1,1)$-form. The above flow is equivalent to the following parabolic Monge-Amp$\acute{e}$re equation
\begin{eqnarray}\label{GKRF3}\frac{\partial\varphi}{\partial t}=\log\frac{\omega^{n}_{\varphi}}{\omega^{n}_{0}}+f+\lambda\varphi ,
\end{eqnarray}
where $f$ is the twisted Ricci potential, i.e. $\sqrt{-1}\partial \overline{\partial } f=-Ric (\omega_{0}) +\lambda \omega_{0} +\theta $.
Let
$$S=|\nabla_{0}g_{\varphi}|^{2}_{\omega_{\varphi}}=g_{\varphi}^{i\bar{j}}g_{\varphi}^{k\bar{l}}g_{\varphi}^{p\bar{q}}\nabla_{0i}g_{\varphi k\bar{q}}\overline{\nabla}_{0\bar{j}}g_{\varphi p\bar{l}},$$
where $\nabla_0$ denotes the covariant derivative with respect to the metric $\omega_{0}$. Define $h^{i}_{\ k}=g_{0}^{i\overline{j}}g_{\varphi\overline{j}k}$ and $X^{k}_{il}=(\nabla_{i}h\cdot h^{-1})^{k}_{\ l}$, by direct computation, we have
\begin{eqnarray}\label{2015.1.24.3}X^{k}_{il}&=&\Gamma^{k}_{\varphi il}-\Gamma^{k}_{0 il},\\
\label{2015.1.24.4}S&=&|X|^{2}_{\omega_{\varphi}},\\
\label{2015.1.24.5}\nabla_{\varphi m}V^{k}_{\ l}-\nabla_{0m}V^{k}_{\ l}&=&X^{k}_{ms}V^{s}_{\ l}-X^{s}_{ml}V^{k}_{\ s}.
\end{eqnarray}
Here we let $\nabla_{\varphi}$ and $\Gamma_{\varphi}$ be the covariant derivative and Christoffel symbol respectively under the metric $\omega_\varphi$, and $\Gamma_{0}$ be the Christoffel symbol with respect to the metric $\omega_{0}$. In the following, the norms $\|\cdot\|_{C^{k}}$ and $\|\cdot\|_{C^{k,\alpha}}$ are all related to the fixed metric $\omega_{0}$ unless there is a special statement. We denote the curvature tensor of $\omega_{\varphi}$ by $Rm_{\varphi}$ for convenience.

\begin{pro}\label{2.2} Let $\varphi (\cdot, t)$ be a solution of the equation $(\ref{GKRF3})$ and satisfy
$$N^{-1}\omega_{0}\leq\omega_{\varphi}\leq N\omega_{0}\ \ \ \ \ \ on\ \ \  B_{r}(p)\times[0,T].$$
 Then there exist constant $C'$ and $C''$ such that
\begin{eqnarray*}S&\leq&\frac{C'}{r^{2}},\\
|Rm_{\varphi}|_{\omega_{\varphi}}^{2}&\leq&\frac{C''}{r^{4}}
\end{eqnarray*}
on $B_{\frac{r}{2}}(p)\times[0,T]$. The constant $C'$ depends only on $\omega_{0}$, $N$, $\lambda$, $\|\varphi(\cdot,0)\|_{C^3(B_{r}(p))}$ and $\|\theta\|_{C^1(B_{r}(p))}$; constant $C''$ depends only on $\omega_{0}$, $N$, $\lambda$, $\|\varphi(\cdot,0)\|_{C^4(B_{r}(p))}$ and $\|\theta\|_{C^2(B_{r}(p))}$.

Furthermore,  there exist constants $C_{k}^{1}$, $C_{k}^{2}$ and $C_{k}^{3}$ such that
\begin{eqnarray*}
|D^{k} Rm_{\varphi}|_{\omega_{\varphi}}^{2}&\leq& C_{k}^{1},\\
\|\dot{\varphi}\|_{C^{k+1,\alpha}}&\leq& C_{k}^{2},\\
\|\varphi\|_{C^{k+3,\alpha}}&\leq& C_{k}^{3}
\end{eqnarray*}
for any $k\geq0$ on $B_{\frac{r}{2}}(p)\times[0,T]$. Here constants $C_{k}^{1}$, $C_{k}^{2}$ and $C_{k}^{3}$ depend only on $\omega_{0}$, $N$, $\lambda$, $\|\varphi(\cdot,0)\|_{C^{k+4}(B_{r}(p))}$, $\|\theta\|_{C^{k+2}(B_{r}(p))}$,  $\|\varphi\|_{C^0(B_{r}(p)\times[0,T])}$ and $\|f\|_{C^0(B_{r}(p))}$.
\end{pro}

{\bf Proof:}\ \ By direct calculation, we have
\begin{eqnarray*}(\frac{d}{dt}-\triangle_{\omega_{\varphi}})S&=&g_{\varphi}^{m\overline{\gamma}}g_{\varphi\overline{\mu}\beta}g_{\varphi}^{l\overline{\alpha}}
((g_{\varphi}^{\beta \bar{s}}\nabla_{\varphi m}\theta_{\bar{s}l}-\nabla_{\varphi}^{\overline{q}}R_{0\ l\overline{q}m}^{\ \beta})X^{\overline{\mu}}_{\overline{\gamma}\overline{\alpha}}+X^{\beta}_{ml}(g_{\varphi}^{\overline{\mu}s}\nabla_{\varphi \overline{\gamma}}\theta_{ s \overline{\alpha}}-\nabla_{\varphi}^{q}R_{0\ \overline{\alpha}q\overline{\gamma}}^{\ \overline{\mu}}))\\
&\ &-X^{\beta}_{ml}X^{\overline{\mu}}_{\overline{\gamma}\overline{\alpha}}
(\theta_{p\bar{q}}g_{\varphi}^{p\overline{\gamma}}g_{\varphi}^{m\overline{q}}g_{\varphi\overline{\mu}\beta}g_{\varphi}^{l\overline{\alpha}}
-g_{\varphi}^{m\overline{\gamma}}\theta_{\overline{\mu}\beta}g_{\varphi}^{l\overline{\alpha}}+
g_{\varphi}^{m\overline{\gamma}}g_{\varphi\overline{\mu}\beta}g_{\varphi}^{p\overline{\alpha}}g_{\varphi}^{l\overline{q}}\theta_{p\bar{q}})\\
&\ &-|\nabla_{\varphi}X|_{\omega_{\varphi}}^{2}-|\overline{\nabla}_{\varphi}X|_{\omega_{\varphi}}^{2}-\lambda S.
\end{eqnarray*}
By (\ref{2015.1.24.5}), we know\begin{eqnarray*}
\nabla_{\varphi m}\theta_{l\bar{q}}&=&\nabla_{0 m}\theta_{l\bar{q}}-X^{s}_{ml}\theta_{s\bar{q}},\\
\nabla_{\varphi p}R_{0\ l\overline{q}m}^{\ \beta}&=&\nabla_{0 p}R_{0\ l\overline{q}m}^{\ \beta}+X^{\beta}_{ps}R_{0\ l\overline{q}m}^{\ s}-X^{s}_{pl}R_{0\ s\overline{q}m}^{\ \beta}-X^{s}_{pm}R_{0\ l\overline{q}s}^{\ \beta}.
\end{eqnarray*}
Hence, the evolution equation of $S$ can be written as
\begin{eqnarray}\label{3.20.1}(\frac{d}{dt}-\triangle_{\omega_{\varphi}})S\leq C(S+1)-|\nabla_{\varphi}X|_{\omega_{\varphi}}^{2}-|\overline{\nabla}_{\varphi}X|_{\omega_{\varphi}}^{2},
\end{eqnarray}
where $C$ depends only on $N$, $\lambda$, $\|Rm (\omega_{0})\|_{C^1(B_{r}(p))}$ and $\|\theta\|_{C^1(B_{r}(p))}$. Let $r=r_{0}>r_{1}>\frac{r}{2}$ and $\psi$ be a nonnegative $C^{\infty}$ cut-off function that is identically equal to $1$ on $\overline{B_{r_{1}}}(p)$ and vanishes outside $B_{r}(p)$. We may assume that
$$|\partial\psi|_{\omega_{0}}^{2},\ |\sqrt{-1}\partial\bar{\partial}\psi|_{\omega_{0}}\leq\frac{C}{r^{2}}.$$
Through computation, we have
\begin{eqnarray}\label{3.20.2}(\frac{d}{dt}-\triangle_{\omega_{\varphi}})(\psi^{2}S)\leq \frac{C}{r^{2}}S+C,
\end{eqnarray}
\begin{eqnarray}\label{3.20.3}\nonumber(\frac{d}{dt}-\triangle_{\omega_{\varphi}})tr h&=&\lambda tr h+g_{0}^{i\overline{j}}\theta_{i\overline{j}}-g_{\varphi}^{p\overline{q}}
g_{0}^{\beta\overline{\gamma}}g_{\varphi \overline{\gamma}\alpha}R_{0\ \beta\overline{q}p}^{\ \alpha}-g_{0}^{j\overline{s}}g_{\varphi}^{p\overline{q}}g_{\varphi}^{m\overline{k}}\varphi_{j\overline{k}p}\varphi_{\overline{s}m\overline{q}}\\
&\leq&C-\frac{1}{N}S.
\end{eqnarray}
From $(\ref{3.20.2})$ and $(\ref{3.20.3})$, we obtain
\begin{eqnarray}\label{3.20.4}(\frac{d}{dt}-\triangle_{\omega_{\varphi}})(\psi^{2}S+Btr h)\leq(\frac{C}{r^{2}}-\frac{B}{N})S+(B+1)C.\end{eqnarray}

Let $(x_{0},t_{0})$ be the maximum point of $\psi^{2}S+Btr h$ on $\overline{B_{r}(p)}\times[0,T]$. If $t_{0}=0$, then $S$ is bounded by the initial data $\|\varphi(\cdot,0)\|_{C^3(B_{r}(p))}$. Then we assume that $t_{0}>0$ and that $x_{0}$ doesn't lie in the boundary of $B_{r}(p)$. By maximum principle,
\begin{eqnarray}\label{3.20.5}0\leq(\frac{C}{r^{2}}-\frac{B}{N})S(x_{0},t_{0})+(B+1)C.\end{eqnarray}
Taking $B=\frac{N(C+1)}{r^{2}}$, we conclude that $S(x_{0},t_{0})\leq C$, where $C$ is independent of $T$. Since $0\leq tr h\leq nN$, we have
\begin{eqnarray}\label{3.20.6}S\leq C+BnN\leq\frac{C}{r^{2}}\ \ \ \ \ \ \ \ \ on\ \overline{B_{r_{1}}(p)}\times[0,T],\end{eqnarray}
where the constant $C$ depends only on $N$, $\lambda$, $\|\varphi(\cdot,0)\|_{C^3(B_{r}(p))}$, $\|\theta\|_{C^1(B_{r}(p))}$ and $\omega_0$. By (\ref{2015.1.24.3}) and (\ref{2015.1.24.5}), we know
\begin{eqnarray}\label{3.20.7}\nonumber(\frac{d}{dt}-\triangle_{\omega_{\varphi}})R_{\varphi\bar{j}i\bar{l}k}&=&+R_{\varphi\bar{j}i}^{\ \ \ p\bar{q}}R_{\varphi\bar{l}k\bar{q}p}+R_{\varphi\bar{l}i}^{\ \ \ p\bar{q}}R_{\varphi\bar{j}k\bar{q}p}-R_{\varphi\bar{j}p\bar{l}}^{\ \ \ \ \bar{q}}R_{\varphi\ i\bar{q}k}^{\ \ p}-R_{\varphi p\bar{l}}R_{\varphi\bar{j}i\ k}^{\ \ \ \ p}\\
&\ &-R_{\varphi\bar{j}h}R_{\varphi\ i\bar{l}k}^{\ \ h}-\nabla_{\varphi\bar{l}}\nabla_{\varphi k}\theta_{i\bar{j}}+\lambda R_{\varphi\bar{j}i\bar{l}k}-\theta_{\bar{j}h}R_{\varphi\ ik\bar{l}}^{\ h}.\\
\label{3.20.8}\nabla_{\varphi\bar{l}}\nabla_{\varphi k}\theta_{i\bar{j}}&=&\nabla_{0\bar{l}}\nabla_{0 k}\theta_{i\bar{j}}-X^{\bar{s}}_{\bar{l}\bar{j}}\nabla_{0 k}\theta_{i\bar{s}}-\nabla_{0\bar{l}}X^{s}_{ki}\theta_{s\bar{j}}\\\nonumber
&\ &-X^{s}_{ki}\nabla_{0 \bar{l}}\theta_{s\bar{j}}+X^{s}_{ki}X^{\bar{t}}_{\bar{l}\bar{j}}\nabla_{0 \bar{l}}\theta_{s\bar{t}},\\
\label{3.20.8'}\nabla_{0\bar{k}}X^{i}_{jl}&=&\partial_{\bar{k}}X^{i}_{jl}=-R_{\varphi\ l\bar{k}j}^{\ i}+R_{0\ l\bar{k}j}^{\ i}.
\end{eqnarray}
Combining the above equalities, we have
\begin{eqnarray}\label{3.20.9}\nonumber(\frac{d}{dt}-\triangle_{\omega_{\varphi}})|Rm_{\varphi}|^{2}_{\omega_{\varphi}}&\leq& C|Rm_{\varphi}|^{3}_{\omega_{\varphi}}+C|Rm_{\varphi}|^{2}_{\omega_{\varphi}}+C|Rm_{\varphi}|_{\omega_{\varphi}}+C S^{\frac{1}{2}}|Rm_{\varphi}|_{\omega_{\varphi}}\\
&\ &+CS|Rm_{\varphi}|_{\omega_{\varphi}}-|\nabla_{\varphi}Rm_{\varphi}|^{2}_{\omega_{\varphi}}-|\overline{\nabla}_{\varphi}Rm_{\varphi}|_{\omega_{\varphi}}^{2}\\ \nonumber
&\leq&C(|Rm_{\varphi}|^{3}_{\omega_{\varphi}}+1+\frac{|Rm_{\varphi}|_{\omega_{\varphi}}}{r^{2}})-|\nabla_{\varphi}Rm_{\varphi}|_{\omega_{\varphi}}^{2}-|\overline{\nabla}_{\varphi}Rm_{\varphi}|_{\omega_{\varphi}}^{2}.
\end{eqnarray}

Next, we show that $|Rm_{\varphi}|^{2}_{\omega_{\varphi}}$ is uniformly bounded. We fix a smaller radius $r_{2}$ satisfying $r_{1}>r_{2}>\frac{r}{2}$. Let $\rho$ be a cut-off function identically equal to $1$ on $\overline{B_{r_{2}}}(p)$ and identically equal to $0$ outside $B_{r_{1}}$. We also let $\rho$ satisfy
$$|\partial\rho|_{\omega_{0}}^{2},\ |\sqrt{-1}\partial\bar{\partial}\rho|_{\omega_{0}}\leq\frac{C}{r^{2}}$$
for some uniform constant $C$. From the former part we know that $S$ is bounded by $\frac{C}{r^{2}}$ on $B_{r_{1}}(p)$. Let
\begin{eqnarray}\label{3.20.10}K=\frac{\hat{C}}{r^{2}},\end{eqnarray}
where $\hat{C}$ is a constant to be determined later and is large enough that $\frac{K}{2}\leq K-S\leq K$. We consider
\begin{eqnarray}\label{3.20.11}F=\rho^{2}\frac{|Rm_{\varphi}|^{2}_{\omega_{\varphi}}}{K-S}+AS.\end{eqnarray}
By computing, we have
\begin{eqnarray}\label{3.20.12}\nonumber(\frac{d}{dt}-\triangle_{\omega_{\varphi}})F&=&(-\triangle_{\omega_{\varphi}}\rho^{2})
\frac{|Rm_{\varphi}|^{2}_{\omega_{\varphi}}}{K-S}+\rho^{2}\frac{|Rm_{\varphi}|^{2}_{\omega_{\varphi}}}{(K-S)^{2}}(\frac{d}{dt}-\triangle_{\omega_{\varphi}})S\\ \nonumber
&\ &+\rho^{2}\frac{1}{K-S}(\frac{d}{dt}-\triangle_{\omega_{\varphi}})|Rm_{\varphi}|^{2}_{\omega_{\varphi}}-4Re\langle\rho\frac{\nabla_{\varphi}\rho}{K-S},\nabla_{\varphi}|Rm_{\varphi}|^{2}_{\omega_{\varphi}}\rangle_{\omega_{\varphi}} \\
&\ &-4Re\langle\rho\frac{|Rm_{\varphi}|^{2}_{\omega_{\varphi}}}{(K-S)^2}\nabla_{\varphi}S,\nabla_{\varphi}\rho\rangle_{\omega_{\varphi}}-2\frac{\rho^{2}|Rm_{\varphi}|^{2}_{\omega_{\varphi}}}{(K-S)^3}|\nabla_{\varphi}S|^{2}_{\omega_{\varphi}}\\  \nonumber
&\ &-2Re\langle\rho^{2}\frac{\nabla_{\varphi}S}{(K-S)^2},\nabla_{\varphi}|Rm_{\varphi}|^{2}_{\omega_{\varphi}}\rangle_{\omega_{\varphi}}+A(\frac{d}{dt}-\triangle_{\omega_{\varphi}})S.
\end{eqnarray}

As in the previous part, we only consider an inner point $(x_{0},t_{0})$ which is a maximum point of $F$ achieved on $\overline{B_{r_{1}}(p)}\times[0,T]$. We use the fact that $\nabla F=0$ at this point, then we get
\begin{eqnarray}\label{3.20.13} 2\rho\nabla_{\varphi}\rho\frac{|Rm_{\varphi}|^{2}_{\omega_{\varphi}}}{K-S}+\rho^2\frac{\nabla_{\varphi}|Rm_{\varphi}|^{2}_{\omega_{\varphi}}}{K-S}+\rho^2\frac{|Rm_{\varphi}|^{2}_{\omega_{\varphi}}\nabla_{\varphi}S}{(K-S)^2}+A\nabla_{\varphi}S=0.\end{eqnarray}

Our goal is to show that at $(x_{0},t_{0})$ we have $|Rm_{\varphi}|^{2}_{\omega_{\varphi}}\leq\frac{C}{r^{4}}$. Without loss of generality, we assume that $|Rm_{\varphi}|^{3}_{\omega_{\varphi}}\geq1+\frac{|Rm_{\varphi}|_{\omega_{\varphi}}}{r^{2}}$. By (\ref{3.20.9}) and $(\ref{3.20.1})$, we have
\begin{eqnarray}\label{3.20.15}(\frac{d}{dt}-\triangle_{\omega_{\varphi}})|Rm_{\varphi}|^{2}_{\omega_{\varphi}}&\leq & C|Rm_{\varphi}|^{3}_{\omega_{\varphi}}-|\nabla_{\varphi}Rm_{\varphi}|_{\omega_{\varphi}}^{2}-|\overline{\nabla}_{\varphi}Rm_{\varphi}|_{\omega_{\varphi}}^{2},\\
\label{3.20.17}(\frac{d}{dt}-\triangle_{\omega_{\varphi}})S&\leq & \frac{C}{r^{2}}-|\nabla_{\varphi}X|_{\omega_{\varphi}}^{2}-|\overline{\nabla}_{\varphi}X|_{\omega_{\varphi}}^{2}
\end{eqnarray}
on $B_{r_{1}}(p)$.
We also note that
\begin{eqnarray}\label{3.20.16}|\nabla_{\varphi}|Rm_{\varphi}|_{\omega_{\varphi}}^{2}|_{\omega_{\varphi}}&\leq&|Rm_{\varphi}|_{\omega_{\varphi}}(|\nabla_{\varphi}Rm_{\varphi}|_{\omega_{\varphi}}+|\overline{\nabla}_{\varphi}Rm_{\varphi}|_{\omega_{\varphi}}),\\
\label{3.20.16'}|\nabla_{\varphi}S|^{2}_{\omega_{\varphi}}&\leq&2S(|\nabla_{\varphi}X|_{\omega_{\varphi}}^{2}+|\overline{\nabla}_{\varphi}X|_{\omega_{\varphi}}^{2}).
\end{eqnarray}
Putting $(\ref{3.20.13})$-$(\ref{3.20.16'})$ into $(\ref{3.20.12})$, then at $(x_{0},t_{0})$, we have
\begin{eqnarray}\label{3.20.18}\nonumber(\frac{d}{dt}-\triangle_{\omega_{\varphi}})F&\leq&-A(|\nabla_{\varphi}X|_{\omega_{\varphi}}^{2}+|\overline{\nabla}_{\varphi}X|_{\omega_{\varphi}}^{2})
+\frac{AC}{r^{2}}+\frac{C|Rm_{\varphi}|_{\omega_{\varphi}}^{2}}{Kr^{2}}+\frac{C\rho^{2}|Rm_{\varphi}|_{\omega_{\varphi}}^{2}}{K^{2}r^{2}}\\
&\ &-\frac{\rho^{2}|Rm_{\varphi}|_{\omega_{\varphi}}^{2}(|\nabla_{\varphi}X|_{\omega_{\varphi}}^{2}+|\overline{\nabla}_{\varphi}X|_{\omega_{\varphi}}^{2})}{K^{2}}
+\frac{C\rho^{2}|Rm_{\varphi}|_{\omega_{\varphi}}^{3}}{K}\\ \nonumber
&\ &-\frac{\rho^{2}(|\nabla_{\varphi}Rm_{\varphi}|_{\omega_{\varphi}}^{2}+|\overline{\nabla}_{\varphi}Rm_{\varphi}|_{\omega_{\varphi}}^{2})}{K}+\frac{C|Rm_{\varphi}|_{\omega_{\varphi}}^{2}}{Kr^{2}}\\ \nonumber
&\ &+\frac{\rho^{2}(|\nabla_{\varphi}Rm_{\varphi}|_{\omega_{\varphi}}^{2}+|\overline{\nabla}_{\varphi}Rm_{\varphi}|_{\omega_{\varphi}}^{2})}{K}+\frac{8AS(|\nabla_{\varphi}X|_{\omega_{\varphi}}^{2}+|\overline{\nabla}_{\varphi}X|_{\omega_{\varphi}}^{2})}{K}.
\end{eqnarray}
Let $\hat{C}$ in (\ref{3.20.10}) be sufficiently large so that $\frac{8ASQ}{K}\leq
\frac{AQ}{2}$, where we denote $Q=|\nabla_{\varphi}X|_{\omega_{\varphi}}^{2}+|\overline{\nabla}_{\varphi}X|_{\omega_{\varphi}}^{2}$. By $(\ref{3.20.8'})$, we have
\begin{eqnarray}\label{3.20.19}\nonumber\frac{C\rho^{2}|Rm_{\varphi}|_{\omega_{\varphi}}^{3}}{K}&\leq&\frac{\rho^{2}|Rm_{\varphi}|_{\omega_{\varphi}}^{4}}{2K^{2}}+C\rho^{2}|Rm_{\varphi}|_{\omega_{\varphi}}^{2}\\
&\leq&\frac{\rho^{2}|Rm_{\varphi}|_{\omega_{\varphi}}^{2}Q}{K^{2}}+C\rho^{2}|Rm_{\varphi}|_{\omega_{\varphi}}^{2}
\end{eqnarray}
So the evolution equation of $F$ can be controlled as follows,
\begin{eqnarray}\label{3.20.20}(\frac{d}{dt}-\triangle_{\omega_{\varphi}})F&\leq&-\frac{AQ}{2}+\frac{AC}{r^{2}}+C|Rm_{\varphi}|_{\omega_{\varphi}}^{2}\\ \nonumber
&\leq&-\frac{AQ}{2}+\frac{AC}{r^{2}}+\tilde{C}Q+C.
\end{eqnarray}
Now we choose a sufficiently large $A$ such that $A\geq2(\tilde{C}+1)$ and obtain
$$Q\leq\frac{C}{r^{2}}$$
at $(x_{0},t_{0})$. This implies that $|Rm_{\varphi}|_{\omega_{\varphi}}^{2}\leq\frac{C}{r^{2}}$ at this point, where $C$ depends only on $N$, $\lambda$, $S$, $\|\theta\|_{C^2(B_{r}(p))}$ and $\omega_0$. Following that we conclude that $F$ is bounded by $\frac{C}{r^{2}}$ at $(x_{0},t_{0})$, where the constant $C$ is independent of $T$. Hence on $\overline{B_{r_{2}}(p)}\times [0,T]$, we obtain
\begin{eqnarray}\label{3.20.21}|Rm_{\varphi}|_{\omega_{\varphi}}^{2}\leq\frac{C}{r^{4}},\end{eqnarray}
where $C$ depends only on $N$, $\lambda$, $\|\varphi(\cdot,0)\|_{C^{4}(B_{r}(p))}$, $\|\theta\|_{C^2(B_{r}(p))}$ and $\omega_0$.

Now, we prove the $C^{\infty}$ estimates of the metric potential $\varphi$ on $B_{\frac{r}{2}}(p)$, combining with the higher order derivative estimates of the Riemann curvature tensors. Here, when we say that $\varphi$ is $C^{k,\alpha}$, we mean that its $C^{k,\alpha}$ norm can be controlled by a constant depending only on $\omega_{0}$, $N$, $\lambda$, $r$, $\|\theta\|_{C^{k-1}(B_{r}(p))}$, $\|f\|_{C^0(B_{r}(p))}$, $\|\varphi(\cdot,0)\|_{C^{k+1}(B_{r}(p))}$ and $\|\varphi\|_{C^{0}(B_{r}(p)\times [0,T])}$. Likewise replacing $\varphi$ by $\dot{\varphi}$, it means the $C^{k,\alpha}$ norm of $\dot{\varphi}$ controlled by a constant that depends only on $\omega_{0}$, $N$, $\lambda$, $r$, $\|\theta\|_{C^{k+1}(B_{r}(p))}$, $\|f\|_{C^0(B_{r}(p))}$, $\|\varphi(\cdot,0)\|_{C^{k+3}(B_{r}(p))}$ and $\|\varphi\|_{C^{0}(B_{r}(p)\times [0,T])}$. Since $|Rm_{\varphi}|_{\omega_{\varphi}}\leq C$ on $\overline{B_{r_{2}}(p)}$ along the flow $(\ref{GKRF2})$, we know that $\dot{\varphi}$ is $C^{1,\alpha}$. Differentiating the equation $(\ref{GKRF3})$ with respect to $z^{k}$, we get
\begin{eqnarray}\label{GKRF4}\frac{d}{dt}\frac{\partial\varphi}{\partial z^{k}}=\triangle_{\omega_{\varphi}}\frac{\partial\varphi}{\partial z^{k}}+g_{\varphi}^{i\bar{j}}\frac{\partial g_{0i\bar{j}}}{\partial z^{k}}-g_{0}^{i\bar{j}}\frac{\partial g_{0i\bar{j}}}{\partial z^{k}}+\frac{\partial f}{\partial z^{k}}+\lambda\frac{\partial\varphi}{\partial z^{k}}.\end{eqnarray}
From the above Calabi's $C^{3}$ estimate, we know that $\varphi $ is $C^{2, \alpha }$ and then the coefficients of $\triangle_{\omega_{\varphi}}$ are $C^{0,\alpha}$. Since $f$ is the twisted Ricci potential, then
\begin{eqnarray}
\triangle_{\omega_{0}}f=-tr_{\omega_{0}}Ric(\omega_{0})+\lambda n +tr_{\omega_{0}}\theta .
\end{eqnarray}
Hence the $C^{1, \alpha }$-norm of $f$ on $B_{r_{2}}(p)$ only depends on $\omega_{0}$, $\|\theta\|_{C^0(B_{r}(p))}$ and $\|f\|_{C^0(B_{r}(p))}$.
By the standard elliptic Schauder estimates, we conclude that $\varphi$ is $C^{3,\alpha}$ on $B_{r_{3}}(p)\times[0,T]$, where $\frac{r}{2}<r_{3}<r_2$.
By computing, we have
\begin{eqnarray}\label{3.20.22}\nonumber(\frac{d}{dt}-\triangle_{\omega_{\varphi}})|\nabla_{\varphi} Rm_{\varphi}|^{2}_{\omega_{\varphi}}&\leq&-|\nabla_{\varphi}\nabla_{\varphi} Rm_{\varphi}|^{2}_{\omega_{\varphi}}-|\bar{\nabla}_{\varphi}\nabla_{\varphi} Rm_{\varphi}|^{2}_{\omega_{\varphi}}\\
&\ &+C|\nabla_{\varphi} Rm_{\varphi}|^{2}_{\omega_{\varphi}}+C|\nabla_{\varphi}\theta|_{\omega_{\varphi}}|\nabla_{\varphi} Rm_{\varphi}|_{\omega_{\varphi}}\\ \nonumber
&\ &+C|\nabla_{\varphi}\bar{\nabla}_{\varphi}\nabla_{\varphi} \theta|_{\omega_{\varphi}}|\nabla_{\varphi} Rm_{\varphi}|_{\omega_{\varphi}},
\end{eqnarray}
where $C$ depends only on $N$, $\lambda$, $\|\theta\|_{C^{0}(B_{r}(p))}$ and $|Rm_{\varphi}|^{2}_{\omega_{\varphi}}$. By $(\ref{3.20.8})$ and $(\ref{3.20.8'})$, we know
\begin{eqnarray}\label{3.20.23}|\nabla_{\varphi}\theta|_{\omega_{\varphi}}&\leq&C,\\
|\nabla_{\varphi}\bar{\nabla}_{\varphi}\nabla_{\varphi} \theta|_{\omega_{\varphi}}&\leq&C(1+|\nabla_{\varphi} Rm_{\varphi}|_{\omega_{\varphi}}+|\nabla_{\varphi} X|_{\omega_{\varphi}}).
\end{eqnarray}
So we have
\begin{eqnarray}\label{3.20.24}(\frac{d}{dt}-\triangle_{\omega_{\varphi}})|\nabla_{\varphi} Rm_{\varphi}|^{2}_{\omega_{\varphi}}&\leq&-|\nabla_{\varphi}\nabla_{\varphi} Rm_{\varphi}|^{2}_{\omega_{\varphi}}-|\bar{\nabla}_{\varphi}\nabla_{\varphi} Rm_{\varphi}|^{2}_{\omega_{\varphi}}\\ \nonumber
&\ &+C|\nabla_{\varphi} Rm_{\varphi}|^{2}_{\omega_{\varphi}}+|\nabla_{\varphi} X|^{2}_{\omega_{\varphi}}+C.
\end{eqnarray}

Let $\varrho$ be a cut-off function, identically equal to $1$ on $\overline{B_{r_{3}}}(p)$ and identically equal to $0$ outside $B_{r_{2}}$. As before we can assume
$$|\partial\varrho|_{\omega_{0}}^{2},\ |\sqrt{-1}\partial\overline{\partial}\varrho|_{\omega_{0}}\leq C$$
for some uniform constant $C$ depending only on $\omega_0$, $N$ and $r$. From the former part we know that $S$ and $|Rm_{\varphi}|_{\omega_{\varphi}}^{2}$  on $B_{r_{2}}(p)$ are bounded by a uniform constant. Define $H=\varrho^{2}|\nabla_{\varphi} Rm_{\varphi}|^{2}_{\omega_{\varphi}}+S+B|Rm_{\varphi}|^{2}_{\omega_{\varphi}}$, where $B$ will be determined later,
\begin{eqnarray}\label{3.20.25}\nonumber(\frac{d}{dt}-\triangle_{\omega_{\varphi}})H&\leq&-\varrho^{2}|\nabla_{\varphi}\nabla_{\varphi} Rm_{\varphi}|^{2}_{\omega_{\varphi}}-\varrho^{2}|\overline{\nabla}_{\varphi}\nabla_{\varphi} Rm_{\varphi}|^{2}_{\omega_{\varphi}}+C|\nabla_{\varphi} Rm_{\varphi}|^{2}_{\omega_{\varphi}}\\ \nonumber
&\ &+|\nabla_{\varphi} X|^{2}_{\omega_{\varphi}}-2Re\langle\nabla_{\varphi}\varrho^2,\nabla_{\varphi}|\nabla Rm_{\varphi}|^{2}_{\omega_{\varphi}}\rangle_{\omega_{\varphi}}-|\nabla_{\varphi}X|_{\omega_{\varphi}}^{2}\\ \nonumber
&\ &
-B|\nabla_{\varphi}Rm_{\varphi}|_{\omega_{\varphi}}^{2}-B|\overline{\nabla}_{\varphi}Rm_{\varphi}|_{\omega_{\varphi}}^{2}+C\\
&\leq&(C-2B)|\nabla_{\varphi} Rm_{\varphi}|^{2}_{\omega_{\varphi}}+C.
\end{eqnarray}

Let $(x_{0},t_{0})$ be the maximum point of $H$ on $\overline{B_{r_{2}}(p)}\times[0,T]$. We assume that $t_{0}>0$ and that $x_{0}$ doesn't lie in the boundary of $B_{r_2}(p)$. We choose $2B=C+1$, by maximum principle, at this point, we have
\begin{eqnarray}\label{3.20.26}|\nabla_{\varphi}Rm_{\varphi}|_{\omega_{\varphi}}^{2}\leq C,\end{eqnarray}
where $C$ depends only on $N$, $\lambda$, $r$, $\|\theta\|_{C^3(B_{r}(p))}$, $|Rm_{\varphi}|^{2}_{\omega_{\varphi}}$ and $\omega_0$. Thus, at $(x_{0},t_{0})$, $H$ is bounded by $C$ independent of $T$. Following the above argument, on $\overline{B_{r_{3}}(p)}\times [0,T]$, we obtain
\begin{eqnarray}\label{3.20.27}|\nabla_{\varphi} Rm_{\varphi}|_{\omega_{\varphi}}^{2}\leq C,\end{eqnarray}
where $C$ depends only on $N$, $\lambda$, $r$, $\|\varphi(\cdot,0)\|_{C^{5}(B_{r}(p))}$, $\|\theta\|_{C^3(B_{r}(p))}$ and $\omega_0$.

Differentiating equation $(\ref{GKRF2})$, we have
$$D\sqrt{-1}\partial\bar{\partial}\dot{\varphi}=D Ric(\omega_{\varphi})+D\theta,$$
where $D$ denotes the covariant derivative with respect to the metric $\omega_{\varphi}$. Taking trace on both side with the metric $\omega_{\varphi}$, we have
\begin{eqnarray}\label{3.21.1}|\triangle_{\omega_{\varphi}} D\dot{\varphi}|\leq |Rm_{\varphi}|_{\omega_{\varphi}}|\overline{\nabla}\dot{\varphi}|+|D Rm_{\varphi}|_{\omega_{\varphi}}+C|X|_{\omega_{\varphi}}+C.\end{eqnarray}
Since $\dot{\varphi}$ is $C^{1,\alpha}$, $|Rm_{\varphi}|_{\omega_{\varphi}}$, $|D Rm_{\varphi}|_{\omega_{\varphi}}$ and $|X|_{\omega_{\varphi}}$ is uniformly bounded, we conclude that $D\dot{\varphi}$ is $C^{1,\alpha}$, and $\dot{\varphi}$ is $C^{2,\alpha}$. Differentiating equation $(\ref{GKRF3})$ two times and using the elliptic Schauder estimates, we know that $\varphi$ is $C^{4,\alpha}$ on $B_{r_{4}}(p)\times[0,T]$, where $\frac{r}{2}<r_{4}<r_3$.

Now we claim that $|D^{k} Rm_{\varphi}|_{\omega_{\varphi}}^{2}\leq C$, $\dot{\varphi}$ is $C^{k+1,\alpha}$ and $\varphi$ is $C^{k+3,\alpha}$ are established for the same $k$ on $\overline{B_{r_{k+3}}(p)}\times [0,T]$, where $C$ depends only on $N$, $\lambda$, $r$, $\|\varphi(\cdot,0)\|_{C^{k+4}(B_{r}(p))}$, $\|\varphi\|_{C^{0}(B_{r}(p)\times [0,T])}$, $\|\theta\|_{C^{k+2}(B_{r}(p))}$, $\|f\|_{C^0(B_{r}(p))}$ and $\omega_0$, $r_{k}>r_{k+1}>\frac{r}{2}$ for any $k\geq0$. We argue it by induction. First, when $k=0,1$, this claim is established. Assume that
\begin{eqnarray}\label{2.6}|D^{j} Rm_{\varphi}|_{\omega_{\varphi}}^{2}\leq C,\ \ \ \|\dot{\varphi}\|\ is\ C^{j+1,\alpha},\ \ \ \|\varphi\|\ is\ C^{j+3,\alpha}\end{eqnarray}
hold on $\overline{B_{r_{j+3}}(p)}\times [0,T]$ for all $j\leq k$.

Now we estimate $|D^{k+1} Rm_{\varphi}|_{\omega_{\varphi}}^{2}$, since any covariant derivative of $Rm_{\varphi}$ of order $k+1$ differs from covariant
derivatives of the form $\nabla_{\varphi}^{r}\overline{\nabla}_{\varphi}^{s}Rm_{\varphi}$ by $D^{i}Rm_{\varphi}\ast D^{r+s-2-i}Rm_{\varphi}$ with
$i\geq0$ and $r+s=k+1$, we should only estimate $|\nabla_{\varphi}^{r}\overline{\nabla}_{\varphi}^{s} Rm_{\varphi}|_{\omega_{\varphi}}^{2}$.
\begin{eqnarray}\label{3.21.2}\nonumber&\ &(\frac{d}{dt}-\triangle_{\omega_{\varphi}})\mid\nabla_{\varphi}^{r}\overline{\nabla}_{\varphi}^{s}Rm_{\varphi}\mid^{2}\\ \nonumber
&=&-\mid\nabla_{\varphi}^{r+1}\overline{\nabla}_{\varphi}^{s}Rm_{\varphi}\mid^{2}
-\mid\overline{\nabla}_{\varphi}\nabla_{\varphi}^{r}\overline{\nabla}_{\varphi}^{s}Rm_{\varphi}\mid^{2}-(r+s+2)\mid\nabla_{\varphi}^{r}\overline{\nabla}_{\varphi}^{s}Rm_{\varphi}\mid^{2}\\
&\ &+\sum_{\substack{i+j=s \\ p+l=r}}\nabla_{\varphi}^{p}\overline{\nabla}_{\varphi}^{i}(Rm_{\varphi}+\theta)\ast\nabla_{\varphi}^{l}\overline{\nabla}_{\varphi}^{j}Rm_{\varphi}\ast\nabla_{\varphi}^{r}\overline{\nabla}_{\varphi}^{s}Rm_{\varphi}\\ \nonumber
&\ &+\sum_{\substack{i+j=s \\ p+l=r}}\overline{\nabla}_{\varphi}^{p}\nabla_{\varphi}^{i}(Rm_{\varphi}+\theta)\ast\overline{\nabla}_{\varphi}^{l}\nabla_{\varphi}^{j}Rm_{\varphi}\ast\overline{\nabla}_{\varphi}^{r}\nabla_{\varphi}^{s}Rm_{\varphi}\\ \nonumber
&\ &+\langle\nabla_{\varphi}^{r}\overline{\nabla}_{\varphi}^{s+1}\nabla_{\varphi}\theta,\overline{\nabla}_{\varphi}^{r}\nabla_{\varphi}^{s}Rm_{\varphi}\rangle+\langle\overline{\nabla}_{\varphi}^{r}\nabla_{\varphi}^{s+1}\overline{\nabla}_{\varphi}\theta,\nabla_{\varphi}^{r}\overline{\nabla}_{\varphi}^{s}Rm_{\varphi}\rangle,
\end{eqnarray}
where the $\ast$ symbol indicates general pairings of these tensors. Since $\varphi$ is $C^{k+3,\alpha}$ on $\overline{B_{r_{k+3}}(p)}\times [0,T]$,
\begin{eqnarray}\label{3.21.3}|D^{k+1}\theta|_{\omega_{\varphi}}\leq C\sum_{i=1}^{k}|D^{i} X|_{\omega_{\varphi}}+C\leq C.
\end{eqnarray}
In the case of $r,s\neq0$, combining with $(\ref{3.20.8})$ and $(\ref{3.20.8'})$, we have
\begin{eqnarray}\label{3.21.4}|\nabla_{\varphi}^{r}\overline{\nabla}_{\varphi}^{s+1}\nabla_{\varphi}\theta|_{\omega_{\varphi}}\leq C|\nabla_{\varphi}^{r}\overline{\nabla}_{\varphi}^{s}Rm_{\varphi}|_{\omega_{\varphi}}+C.
\end{eqnarray}
When $r=0$ or $s=0$, without loss of generality, we assume $s=0$,
\begin{eqnarray}\label{3.21.5}|\nabla_{\varphi}^{k+1}\overline{\nabla}_{\varphi}\nabla_{\varphi}\theta|_{\omega_{\varphi}}\leq C|\nabla_{\varphi}^{k+1}X|_{\omega_{\varphi}}+C|\nabla_{\varphi}^{k+1}Rm_{\varphi}|_{\omega_{\varphi}}+C.
\end{eqnarray}
The corresponding evolution equation are as follows.
\begin{eqnarray}\label{3.21.6}\nonumber(\frac{d}{dt}-\triangle_{\omega_{\varphi}})|\nabla_{\varphi}^{r}\overline{\nabla}_{\varphi}^{s}Rm_{\varphi}|^{2}&\leq&-|\nabla_{\varphi}^{r+1}\overline{\nabla}_{\varphi}^{s}Rm_{\varphi}|^{2}
-|\overline{\nabla}_{\varphi}\nabla_{\varphi}^{r}\overline{\nabla}_{\varphi}^{s}Rm_{\varphi}|^{2}\\
&\ &+C|\nabla_{\varphi}^{r}\overline{\nabla}_{\varphi}^{s}Rm_{\varphi}|^{2}+C,
\end{eqnarray}
\begin{eqnarray}\label{3.21.7}\nonumber(\frac{d}{dt}-\triangle_{\omega_{\varphi}})|\nabla_{\varphi}^{k+1}Rm_{\varphi}|^{2}&\leq&-|\nabla_{\varphi}^{k+2}Rm_{\varphi}|^{2}
-|\overline{\nabla}_{\varphi}\nabla_{\varphi}^{k+1}Rm_{\varphi}|^{2}\\
&\ &+C|\nabla_{\varphi}^{k+1}Rm_{\varphi}|^{2}+|\nabla_{\varphi}^{k+1}X|_{\omega_{\varphi}}^{2}+C.
\end{eqnarray}
Let $\vartheta$ be a cut-off function, identically equal to $1$ on $\overline{B_{r'_{k+3}}}(p)$ and identically equal to $0$ outside $B_{r_{k+3}}$, where $\frac{r}{2}<r'_{k+3}<r_{k+3}$. As before we can assume
$$|\partial\vartheta|_{\omega_{0}}^{2},\ |\sqrt{-1}\partial\partial\vartheta|_{\omega_{0}}\leq C$$
for some constant $C$ depending only on $\omega_{0}$, $N$ and $r$. From the former part we know that $|\nabla_{\varphi}^{k}X|_{\omega_{\varphi}}^{2}$ and $|D^{k} Rm_{\varphi}|_{\omega_{\varphi}}^{2}$ are bounded by a uniform constant on $B_{r_{k+3}}(p)$. Then we talk about it in the following two case:

$(1)$ When $r,s\neq0$, we define $G_{1}=\vartheta^{2}|\nabla_{\varphi}^{r}\overline{\nabla}_{\varphi}^{s}Rm_{\varphi}|_{\omega_{\varphi}}^{2}+A_1|\nabla_{\varphi}^{r-1}\overline{\nabla}_{\varphi}^{s}Rm_{\varphi}|_{\omega_{\varphi}}^{2}$;

$(2)$ When $s=0$, we define $G_2=\vartheta^{2}|\nabla_{\varphi}^{k+1}Rm_{\varphi}|_{\omega_{\varphi}}^{2}+A_2|\nabla_{\varphi}^{k}Rm_{\varphi}|_{\omega_{\varphi}}^{2}+|\nabla_{\varphi}^{k}X|_{\omega_{\varphi}}^{2}$.
We first analysis the evolution of $|\nabla_{\varphi}^{k}X|_{\omega_{\varphi}}^{2}$. By direct computation, we have
\begin{eqnarray}\label{3.21.8}(\frac{d}{dt}-\triangle_{\omega_{\varphi}})X^{\beta}_{ml}&=&\nabla_{\varphi m}\theta^{\beta}_{\ l}-\overline{\nabla}^{q}_{\varphi}R^{\ \beta}_{0\ l\bar{q}m},\\
\frac{d}{dt}\Gamma^{\beta}_{\varphi ml}&=&-g_{\varphi}^{\beta\bar{t}}\nabla_{\varphi m}(R_{g_{\varphi}\bar{t}l}-\theta_{\bar{t}l}).
\end{eqnarray}
Since there exists no $Rm_{\varphi}$ in the evolution equation of $X$ and there only exists derivative of $Rm_{\varphi}$ of order $1$ in the evolution equation of Christoffel $\Gamma_\varphi$, we know that there exists derivative of $Rm_{\varphi}$ no more than of order $k$ in the evolution equation of $\nabla_{\varphi}^{k}X$. Combining $\varphi$ is $C^{k+3,\alpha}$, we obtain

\begin{eqnarray}\label{3.21.9}(\frac{d}{dt}-\triangle_{\omega_{\varphi}})|\nabla_{\varphi}^{k}X|_{\omega_{\varphi}}^{2}\leq-|\nabla_{\varphi}^{k+1}X|_{\omega_{\varphi}}^{2}
-|\overline{\nabla}\nabla_{\varphi}^{k}X|_{\omega_{\varphi}}^{2}+C.
\end{eqnarray}
Then by choosing suitable $A_1$ and $A_2$, we have
\begin{eqnarray}\label{3.21.10}\nonumber(\frac{d}{dt}-\triangle_{\omega_{\varphi}})G_1&\leq&-\vartheta^{2}|\nabla_{\varphi}^{r+1}\overline{\nabla}_{\varphi}^{s}Rm_{\varphi}|_{\omega_{\varphi}}^{2}-\vartheta^{2}|\overline{\nabla}_{\varphi}\nabla_{\varphi}^{r}\overline{\nabla}_{\varphi}^{s}Rm_{\varphi}|_{\omega_{\varphi}}^{2}\\ \nonumber
&\ &+C|\nabla_{\varphi}^{r}\overline{\nabla}_{\varphi}^{s}Rm_{\varphi}|_{\omega_{\varphi}}^{2}-2Re\langle\nabla\vartheta^{2},\nabla|\nabla_{\varphi}^{r}\overline{\nabla}_{\varphi}^{s}Rm_{\varphi}|_{\omega_{\varphi}}^{2}\rangle_{\omega_{\varphi}}\\ \nonumber
&\ &-A_1|\nabla_{\varphi}^{r}\overline{\nabla}_{\varphi}^{s}Rm_{\varphi}|_{\omega_{\varphi}}^{2}+C\\
&\leq&-\vartheta^{2}|\nabla_{\varphi}^{r+1}\overline{\nabla}_{\varphi}^{s}Rm_{\varphi}|_{\omega_{\varphi}}^{2}-\vartheta^{2}|\overline{\nabla}_{\varphi}\nabla_{\varphi}^{r}\overline{\nabla}_{\varphi}^{s}Rm_{\varphi}|_{\omega_{\varphi}}^{2}\\ \nonumber
&\ &+C|\nabla_{\varphi}^{r}\overline{\nabla}_{\varphi}^{s}Rm_{\varphi}|_{\omega_{\varphi}}^{2}+\vartheta^{2}|\nabla_{\varphi}^{r+1}\overline{\nabla}_{\varphi}^{s}Rm_{\varphi}|_{\omega_{\varphi}}^{2}\\ \nonumber
&\ &+\vartheta^{2}|\overline{\nabla}_{\varphi}\nabla_{\varphi}^{r}\overline{\nabla}_{\varphi}^{s}Rm_{\varphi}|_{\omega_{\varphi}}^{2}-A_1|\nabla_{\varphi}^{r}\overline{\nabla}_{\varphi}^{s}Rm_{\varphi}|_{\omega_{\varphi}}^{2}+C\\ \nonumber
&\leq&-|\nabla_{\varphi}^{r}\overline{\nabla}_{\varphi}^{s}Rm_{\varphi}|_{\omega_{\varphi}}^{2}+C,\\
\label{3.21.11}\ \nonumber(\frac{d}{dt}-\triangle_{\omega_{\varphi}})G_2&\leq&-\vartheta^{2}|\nabla_{\varphi}^{k+2}Rm_{\varphi}|^{2}
-\vartheta^{2}|\overline{\nabla}_{\varphi}\nabla_{\varphi}^{k+1}Rm_{\varphi}|^{2}\\ \nonumber
&\ &+C|\nabla_{\varphi}^{k+1}Rm_{\varphi}|^{2}+|\nabla_{\varphi}^{k+1}X|_{\omega_{\varphi}}^{2}-|\nabla_{\varphi}^{k+1}X|_{\omega_{\varphi}}^{2}+C\\ \nonumber
&\ &-2Re\langle\nabla\vartheta^{2},\nabla|\nabla_{\varphi}^{k+1}Rm_{\varphi}|_{\omega_{\varphi}}^{2}\rangle_{\omega_{\varphi}}-A_2|\nabla_{\varphi}^{k+2}Rm_{\varphi}|^{2}\\
&\leq&-\vartheta^{2}|\nabla_{\varphi}^{k+2}Rm_{\varphi}|^{2}
-\vartheta^{2}|\overline{\nabla}_{\varphi}\nabla_{\varphi}^{k+1}Rm_{\varphi}|^{2}\\ \nonumber
&\ &+C|\nabla_{\varphi}^{k+1}Rm_{\varphi}|^{2}+|\nabla_{\varphi}^{k+1}X|_{\omega_{\varphi}}^{2}-|\nabla_{\varphi}^{k+1}X|_{\omega_{\varphi}}^{2}+C\\ \nonumber
&\ &+\vartheta^{2}|\nabla_{\varphi}^{k+2}Rm_{\varphi}|^{2}
+\vartheta^{2}|\overline{\nabla}_{\varphi}\nabla_{\varphi}^{k+1}Rm_{\varphi}|^{2}-A_2|\nabla_{\varphi}^{k+2}Rm_{\varphi}|^{2}\\ \nonumber
&\leq&-|\nabla_{\varphi}^{k+2}Rm_{\varphi}|^{2}+C.
\end{eqnarray}
Let $(x_{1},t_{1})$ and $(x_{2},t_{2})$ be the maximum point of $G_1$ and $G_2$ on $\overline{B_{r_{k+3}}(p)}\times[0,T]$ respectively. We assume that $t_{i}>0$ and that $x_{i}$ doesn't lie in the boundary of $B_{r_{k+3}}(p)$ for $i=1,\ 2$. By maximum principle, we have
\begin{eqnarray}\label{3.21.12}|\nabla_{\varphi}^{r}\overline{\nabla}_{\varphi}^{s}Rm_{\varphi}|_{\omega_{\varphi}}^{2}(x_1,t_1)\leq C,\ \ \ |\nabla_{\varphi}^{k+1}Rm_{\varphi}|_{\omega_{\varphi}}^{2}(x_2,t_2)\leq C,\end{eqnarray}
where $C$ depends only on $N$, $\lambda$, $r$, $\|\varphi(\cdot,0)\|_{C^{k+4}(B_{r}(p))}$, $\|\varphi\|_{C^{0}(B_{r}(p)\times [0,T])}$, $\|f\|_{C^0(B_{r}(p))}$, $\|\theta\|_{C^{k+3}(B_{r}(p))}$ and $\omega_0$. Thus $G_i$ is bounded from above uniformly at $(x_{i},t_{i})$. Following this argument, on $\overline{B_{r'_{k+3}}(p)}\times [0,T]$, we have
\begin{eqnarray}\label{3.21.13}|\nabla_{\varphi}^{r}\overline{\nabla}_{\varphi}^{s}Rm_{\varphi}|_{\omega_{\varphi}}^{2}\leq C,\ \ \ |\nabla_{\varphi}^{k+1}Rm_{\varphi}|_{\omega_{\varphi}}^{2}\leq C,\end{eqnarray}
where $C$ depends only on $N$, $\lambda$, $r$, $\|\varphi(\cdot,0)\|_{C^{k+5}(B_{r}(p))}$, $\|\varphi\|_{C^{0}(B_{r}(p)\times [0,T])}$, $\|f\|_{C^0(B_{r}(p))}$, $\|\theta\|_{C^{k+3}(B_{r}(p))}$ and $\omega_0$. Then we prove that $|D^{j} Rm_{\varphi}|_{\omega_{\varphi}}^{2}\leq C$ established for $k+1$ on $\overline{B_{r'_{k+3}}(p)}\times [0,T]$.

From equation $(\ref{GKRF2})$, we have
\begin{eqnarray}\label{3.21.14}|\Delta_{\omega_{\varphi}}D^{k+1}\dot{\varphi}|\leq C(\sum_{i=1}^{k+2}|D^{i-1} Rm_{\varphi}|_{\omega_{\varphi}}|D^{k-i+2}\dot{\varphi}|_{\omega_{\varphi}}+\sum_{i=1}^{k}|D^{i} X|_{\omega_{\varphi}}+1).\end{eqnarray}
By $(\ref{2.6})$, we know that $|\Delta_{\omega_{\varphi}}D^{k+1}\dot{\varphi}|\leq C$, so $D^{k+1}\dot{\varphi}$ is $C^{1,\alpha}$. Then by the assumption, it is easy to see that $\dot{\varphi}$ is $C^{k+2,\alpha}$. By differentiating the parabolic Monge-Amp\'ere equation $(\ref{GKRF3})$ $k+2$ times and using the elliptic Schauder estimates again, we know that $\varphi$ is $C^{k+4,\alpha}$ on $\overline{B_{r_{k+4}}(p)}\times [0,T]$, where $r'_{k+3}>r_{k+4}>\frac{r}{2}$. Hence we get $C^{\infty}$ estimates of $\varphi$ on $B_{\frac{r}{2}}(p)\times[0, T]$.\QEDB

\begin{rem}\label{2015.1.25.1} Considering only the regularity estimates for a single flow $(\ref{GKRF3})$, we can get the local uniform $C^{\infty}$ estimates of $\varphi$ by the standard Schauder estimate of the parabolic equation (see \cite{GLIE}) after getting the Calabi's $C^3$ estimate and the curvature estimate. Since we want to get the conical K\"ahler-Ricci flow by limiting a sequence of the twisted K\"ahler-Ricci flows $(\ref{CMAF1})$ as $\varepsilon\rightarrow0$, we need to get the uniform $C^{\infty}$ estimates of $\varphi_{\varepsilon}(\cdot,t)$ on $\overline{B_{r}}\times [0,T]$, where $\overline{B_{r}}\subset\subset M\setminus D$. But by applying the parabolic Schauder estimates, we can only get the uniform $C^{\infty}$ estimates of $\varphi_{\varepsilon}(\cdot,t)$ on $\overline{B_{r}}\times [\delta,T]$. Here $\delta>0$ and the uniform estimates depends on $\delta$. This is the reason why we apply the elliptic estimates in the proof of Proposition $2.2$. We can also note that the estimates are independent of time $T$, so the results hold also for time intervals $[0,+\infty)$.
\end{rem}

\section{The long-time solution to the conical K\"ahler-Ricci flow}
\setcounter{equation}{0}

In this section, we use the estimates obtained in the preceding section to give a long-time solution to the conical K\"ahler-Ricci flow. We prove the following theorem:
\begin{thm}\label{3.1} Assume $\beta\in(0,1).$ Then there esists a sequence $\{\varepsilon_{i}\}$ satisfying $\varepsilon_{i}\rightarrow0$ as $i\rightarrow+\infty$, such that the flow $(\ref{CMAF1})$ converges to the following equation
\begin{equation}\label{CMAF4}
\begin{cases}
  \frac{\partial \varphi}{\partial t}=\log\frac{\omega_{\varphi}^{n}}{\omega_{0}^{n}}+F_{0}+\beta(k|s|_{h}^{2}+\varphi)+\log|s|_{h}^{2(1-\beta)}\\
  \\
  \varphi|_{t=0}=c_{0}, \\
  \end{cases}
\end{equation}
in the $C_{loc}^{\infty}$ topology outside divisor D. Furthermore, $\omega_{\varphi }=\omega^\ast+\sqrt{-1}\partial\bar{\partial}\varphi$ is a long-time solution to the conical K\"ahler-Ricci flow (\ref{CKRF2}) with initial metric $\omega^*$.
\end{thm}

{\bf Proof:}\ \ Differentiating equation $(\ref{CMAF1})$ with respect to $t$, we have
\begin{equation}\label{1.7.11}\frac{d}{dt}\dot{\varphi}_{\varepsilon}(t)=\triangle_{\omega_{\varphi_{\varepsilon}(t)}}\dot{\varphi}_{\varepsilon}(t)+\beta\dot{\varphi}_{\varepsilon}(t).\end{equation}
According to the maximum principle, we have
\begin{equation}\label{1.7.12}\sup\limits_{M}|\dot{\varphi}_{\varepsilon}(t)|\leq\sup\limits_{M}|e^{\beta t}\dot{\varphi}_{\varepsilon}(0)|,
\end{equation}
where $\dot{\varphi}_{\varepsilon}(0)=\log\frac{\omega_{\varepsilon}^{n}(\varepsilon^{2}+|s|_{h}^{2})^{1-\beta}}{\omega_{0}^{n}}+F_{0}+\beta(k\chi
+c_{\varepsilon0})$, so $\sup\limits_{M}|\dot{\varphi}_{\varepsilon}(t)|\leq Ce^{\beta t}$ for some uniform constant $C$. Then on $M\times[0,T]$, we have $\|\varphi_{\varepsilon}(t)\|_{C^{0}}\leq Ce^{\beta T}$. By proposition \ref{2.1}, there exists constant $C(T)$ satisfying
\begin{eqnarray}\label{3.21.15}C^{-1}(T)\omega_{\varepsilon}\leq\omega_{\varphi_{\varepsilon}}\leq C(T)\omega_{\varepsilon}
\end{eqnarray}
on $M\times[0,T]$. For any $K\subset\subset M\setminus D$, we have
\begin{eqnarray}\label{3.21.16}\frac{1}{N}\omega_{0}\leq\omega_{\varphi_{\varepsilon}}\leq N\omega_{0},\end{eqnarray}
where the uniform constant $N$ depends only on $K$ and $C(T)$. Since the initial data $k\chi+c_{\varepsilon}(0)$, the twisted Ricci potential $F_{0}+\log (\varepsilon^{2} +|s|_{h}^{2})^{1-\beta}$ of $\omega_0$ and the twist form $\theta_\varepsilon$ are $C^{\infty}_{loc}$ uniformly bounded away from divisor $D$, then by Proposition \ref{2.2}, $\varphi_{\varepsilon}+k\chi$ is $C^{\infty}$ bounded uniformly (independent of $\varepsilon$) on $K\times[0,T]$. Let $K$ approximate to $M\setminus D$ and $T$ approximate to $\infty$, by diagonal rule, we get a sequence which we denote $\{\varepsilon_{i}\}$, such that $\varphi_{\varepsilon_{i}}(t)$ converges in $C_{loc}^{\infty}$ topology outside divisor $D$ to a function $\varphi(t)$ that is smooth on $M\setminus D$. From (\ref{3.21.15}), we know that every $\omega_{\varphi (t)}$ is conical K\"ahler metric with cone angle $2\pi \beta $ along the divisor $D$.

Next, we prove that the limit $\varphi(t)$ satisfies the conical K\"ahler-Ricci flow $(\ref{CMAF4})$ globally on $M\times[0,+\infty)$ in the sense of currents. Since
$\log\frac{\omega_{\varphi_{\varepsilon_{i}}}^{n}(\varepsilon_{i}^{2}+|s|_{h}^{2})^{1-\beta}}{\omega_{0}^{n}}$, $ k\chi(\varepsilon_{i}^{2}+|s|_{h}^{2})$ and $\varphi_{\varepsilon_{i}}$ are bounded by some constant which is independent of $\varepsilon$, then for any $(n-1,n-1)$-form $\eta$,  by dominated convergence theorem
\begin{eqnarray*}&\ &\int_M\sqrt{-1}\partial\bar{\partial}\frac{\partial\varphi_{\varepsilon_{i}}}{\partial t}\wedge\eta\\
&=&\int_M\sqrt{-1}\partial\bar{\partial}(\log\frac{\omega_{\varphi_{\varepsilon_{i}}}^{n}(\varepsilon_{i}^{2}+|s|_{h}^{2})^{1-\beta}}{\omega_{0}^{n}}+F_{0}+\beta(k\chi(\varepsilon_{i}^{2}+|s|_{h}^{2})+\varphi_{\varepsilon_{i}})
  )\wedge\eta\\
&=&\int_M\log\frac{\omega_{\varphi_{\varepsilon_{i}}}^{n}(\varepsilon_{i}^{2}+|s|_{h}^{2})^{1-\beta}}{\omega_{0}^{n}}+F_{0}+\beta(k\chi(\varepsilon_{i}^{2}+|s|_{h}^{2})+\varphi_{\varepsilon_{i}})
  \sqrt{-1}\partial\bar{\partial}\eta\\
&\xrightarrow{\varepsilon_{i}\rightarrow 0}&\int_{M}(\log\frac{\omega_{\varphi}^{n}}{\omega_{0}^{n}}+F_{0}+\beta(k|s|_{h}^{2\beta}+\varphi)
  +\log|s|_{h}^{2(1-\beta)})\sqrt{-1}\partial\bar{\partial}\eta \\
&=&\int_M\sqrt{-1}\partial\bar{\partial}(\log\frac{\omega_{\varphi}^{n}}{\omega_{0}^{n}}+F_{0}+\beta(k|s|_{h}^{2\beta}+\varphi)
  +\log|s|_{h}^{2(1-\beta)})\wedge\eta.
\end{eqnarray*}

On the other hand, let $K\subset\subset M\setminus D$ be a compact subset, $\int_{M\setminus K}\sqrt{-1}\partial\bar{\partial}\eta=\delta$, and $\delta\rightarrow0$ as $K\rightarrow M\setminus D$. By the fact that both $\frac{\partial\varphi_{\varepsilon_{i}}}{\partial t}$ and $\frac{\partial\varphi}{\partial t}$ are uniformly bounded,
\begin{eqnarray*}
&\ &|\int_{M}(\frac{\partial\varphi_{\varepsilon_{i}}}{\partial t}-\frac{\partial\varphi}{\partial t})\sqrt{-1}\partial\bar{\partial}\eta|\\
&=&|\int_{K}(\frac{\partial\varphi_{\varepsilon_{i}}}{\partial t}-\frac{\partial\varphi}{\partial t})\sqrt{-1}\partial\bar{\partial}\eta+\int_{M\setminus K}(\frac{\partial\varphi_{\varepsilon_{i}}}{\partial t}-\frac{\partial\varphi}{\partial t})\sqrt{-1}\partial\bar{\partial}\eta|\\
&\leq&|\int_{K}(\frac{\partial\varphi_{\varepsilon_{i}}}{\partial t}-\frac{\partial\varphi}{\partial t})\sqrt{-1}\partial\bar{\partial}\eta|+C(T)\delta\\
&\xrightarrow{\varepsilon_{i}\rightarrow 0}&C(T)\delta.
\end{eqnarray*}
When $K\rightarrow M\setminus D$, we have
\begin{eqnarray*}
\int_M\sqrt{-1}\partial\bar{\partial}\frac{\partial\varphi_{\varepsilon_{i}}}{\partial t}\wedge\eta=\int_{M}\frac{\partial\varphi_{\varepsilon_{i}}}{\partial t}\sqrt{-1}\partial\bar{\partial}\eta
&\xrightarrow{\varepsilon_{i}\rightarrow 0}&\int_{M}\frac{\partial\varphi}{\partial t}\sqrt{-1}\partial\bar{\partial}\eta.
\end{eqnarray*}
Hence the limit $\omega_{\varphi} (\cdot , t)$ satisfies flow (\ref{CKRF2}) on $M\times[0,+\infty)$ in the current sense.\QEDB

\begin{pro}\label{3.2} For any $t\in[0,+\infty)$, the potential $\varphi(t)$ is H$\ddot{o}$lder continuous with respect to the metric $\omega_{0}$ on $M$.
\end{pro}

{\bf Proof:}\ \ Let $\phi=\varphi+k|s|_{h}^{2\beta}$. For any $t$, we fix $T>t$. From Theorem \ref{3.1}, we have $\|\dot{\phi}(t)\|_{C^{0}}\leq C(T)$ and $\|\phi(t)\|_{C^{0}}\leq C(T)$ on $M\setminus D\times[0,T]$. Flow $(\ref{CMAF4})$ can be written as
\begin{eqnarray}(\omega_{0}+\sqrt{-1}\partial\bar{\partial}\phi)^{n}=e^{\dot{\phi}-F_{0}-\beta\phi}\frac{\omega_{0}^{n}}{|s|_{h}^{2(1-\beta)}}
\end{eqnarray}
on $M\setminus D$. Since $\beta\in(0,1)$, there exists $\delta$ such that $2(1-\beta)(1+\delta)<2$.
\begin{eqnarray*}\int_{M}e^{(\dot{\phi}-F_{0}-\beta\phi-\log|s|_{h}^{2(1-\beta)})(1+\delta)}dV_{0}\leq C(T)\int_{M}\frac{1}{|s|_{h}^{2(1-\beta)(1+\delta)}}dV_{0}\leq C(T).
\end{eqnarray*}
Then by the $L^p$ estimate of S. Kolodziej \cite{K}, we conclude that the potential $\varphi(t)$ is H$\ddot{o}$lder continuous with respect to the metric $\omega_{0}$ on $M$.\QEDB

\medskip

\begin{rem}\label{3.3} From Theorem \ref{3.1} and Proposition \ref{3.2}, we have
\begin{eqnarray}\ \ \ \ \|\dot{\varphi}\|_{C^{0}}\leq C(T),\ \ \ \ \|\varphi\|_{C^{\alpha,\frac{\alpha}{2}}}\leq C(T),\ \ \ \ \ C^{-1}(T)\omega\leq\omega_{\varphi}\leq C(T)\omega
\end{eqnarray}
on $M\setminus D\times[0,T]$. By the uniqueness theorem of the weak conical K\"ahler-Ricci flow (see Lemma $3.2$ in \cite{YQW}) and the existence of  long-time solution to the strong conical K\"ahler-Ricci flow proved in \cite{CW1}, we conclude that the conical K\"ahler-Ricci flow constructed in Theorem \ref{3.1} must be the strong conical K\"ahler-Ricci flow.
\end{rem}

\section{Uniform Perelman's estimates along the twisted K\"ahler-Ricci flows}
\setcounter{equation}{0}

In this section, we first obtain a uniform lower bound for the twisted scalar curvature $R(g_{\varepsilon}(t))-tr_{g_{\varepsilon}(t)}\theta_{\varepsilon}$ in some time interval $[\delta, +\infty)$, where $\delta>0$. This conclusion is very important to get the uniform Perelman's estimates along the twisted K\"ahler-Ricci flow (\ref{GKRF1}), because we have no uniform lower bound of the initial twisted scalar curvature $R(g_{\varepsilon}(0))-tr_{g_{\varepsilon}(0)}\theta_{\varepsilon}$ when $\beta\in(\frac{1}{2}, 1)$.

\begin{pro}\label{1.8.1} $t^2(R(g_{\varepsilon}(t))-tr_{g_{\varepsilon}(t)}\theta_{\varepsilon})$ is uniformly bounded from below along the flow $(\ref{GKRF1})$, i.e. there exists a uniform constant $C$, such that
\begin{equation}\label{3.22.9}t^2(R(g_{\varepsilon}(t))-tr_{g_{\varepsilon}(t)}\theta_{\varepsilon})\geq-C \end{equation}
for any $t$ and $\varepsilon$, while the constant $C$ only depends on $\beta$ and $n$. In particular,
\begin{equation}\label{3.22.911}R(g_{\varepsilon}(t))-tr_{g_{\varepsilon}(t)}\theta_{\varepsilon}\geq-C \end{equation}
when $t\geq1$.
\end{pro}

{\bf Proof:}\ \ First, we derive the evolution equation of $t^2(R(g_{\varepsilon}(t))-tr_{g_{\varepsilon}(t)}\theta_{\varepsilon})$ as follows.
\begin{eqnarray*}&\ &(\frac{d}{dt}-\triangle_{g_{\varepsilon}(t)})(t^2(R(g_{\varepsilon}(t))-tr_{g_{\varepsilon}(t)}\theta_{\varepsilon}))\\
&=&t^2|R_{\varepsilon i\bar{j}}-\theta_{\varepsilon i\bar{j}}|^{2}-\beta t^2(R(g_{\varepsilon}(t))-tr_{g_{\varepsilon}(t)}\theta_{\varepsilon})+2t(R(g_{\varepsilon}(t))-tr_{g_{\varepsilon}(t)}\theta_{\varepsilon})
\end{eqnarray*}
Assume that $(t_0, x_0)$ is the minimum point of $t^2(R(g_{\varepsilon}(t))-tr_{g_{\varepsilon}(t)}\theta_{\varepsilon})$ on $[0,T]\times M$.

$Case\ 1$, $t_0=0$, then we have $t^2(R(g_{\varepsilon}(t))-tr_{g_{\varepsilon}(t)}\theta_{\varepsilon})\geq0$.

$Case\ 2$, $t_0\geq\frac{2}{\beta}$, then at $(t_0, x_0)$
\begin{eqnarray*}0&\geq&(2t_0-\beta t_{0}^2)(R(g_{\varepsilon}(t_0))-tr_{g_{\varepsilon}(t_0)}\theta_{\varepsilon}( x_0))).
\end{eqnarray*}
Hence $R(g_{\varepsilon}(t_0))-tr_{g_{\varepsilon}(t_0)}\theta_{\varepsilon}(x_0)\geq 0$, and then $t^2(R(g_{\varepsilon}(t))-tr_{g_{\varepsilon}(t)}\theta_{\varepsilon})\geq 0$.

$Case\ 3$, $0<t_0\leq\frac{2}{\beta}$, without loss of generality, we can assume $R(g_{\varepsilon}(t_0))-tr_{g_{\varepsilon}(t_0)}\theta_{\varepsilon}( x_0)\leq 0$. By inequality
$$|R_{\varepsilon i\bar{j}}-\theta_{\varepsilon i\bar{j}}|^{2}\geq\frac{(R(g_{\varepsilon}(t))-tr_{g_{\varepsilon}(t)}\theta_{\varepsilon})^2}{n},$$
at $(t_0, x_0)$, we have
\begin{eqnarray*}0&\geq& t_0^2\frac{(R(g_{\varepsilon}(t_0))-tr_{g_{\varepsilon}(t_0)}\theta_{\varepsilon}(x_0))^2}{n}+2t_0(R(g_{\varepsilon}(t_0))-tr_{g_{\varepsilon}(t_0)}\theta_{\varepsilon}(x_0))\\
&=&(t_0\frac{R(g_{\varepsilon}(t_0))-tr_{g_{\varepsilon}(t_0)}\theta_{\varepsilon}(x_0)}{\sqrt{n}}+\sqrt{n})^2-n.
\end{eqnarray*}
So $t_0^{2}R(g_{\varepsilon}(t_0))-tr_{g_{\varepsilon}(t_0)}\theta_{\varepsilon}( x_0)\geq -2t_0n>-\frac{4n}{\beta}$. Hence $t^2(R(g_{\varepsilon}(t))-tr_{g_{\varepsilon}(t)}\theta_{\varepsilon})\geq-\frac{4n}{\beta}.$

By the above argument, we conclude that there exists a uniform constant only depending on $n$ and $\beta$, such that $t^2(R(g_{\varepsilon}(t))-tr_{g_{\varepsilon}(t)}\theta_{\varepsilon})\geq-C$ for any $t\geq 0$ and $\varepsilon$. When $t\geq1$, we have $R(g_{\varepsilon}(t))-tr_{g_{\varepsilon}(t)}\theta_{\varepsilon}\geq-C$. \QEDB

\medskip

Now, using the above uniform estimate, we prove the uniform Perelman's estimates along the flows (\ref{GKRF1}) for $t\geq1$ by following the argument of N. Sesum and G. Tian in \cite{NSGT} (see also the twisted case in \cite{JWL}). Because we have no uniform bound on the  initial data (i.e. $t=0$), we will make some small changes in the  argument. In the following sections, we let $\nabla$ be the $(1,0)$-type covariant derivative with respect to the metric $g_\varepsilon(t)$.

\begin{thm}\label{1.8.2} Let $g_{\varepsilon}(t)$ be a solution of the twisted K$\ddot{a}$hler Ricci flow, i.e. the corresponding form $\omega_{\varepsilon}(t)$ satisfies the equation $(\ref{GKRF1})$ with initial metric $\omega_{\varepsilon}$, $u_{\varepsilon}(t)\in C^{\infty}(M)$ is the twisted Ricci potential satisfying
\begin{equation}\label{4.1}-Ric(\omega_{\varepsilon}(t))+\beta\omega_{\varepsilon}(t)+\theta_{\varepsilon}=\sqrt{-1}\partial\bar{\partial}u_{\varepsilon}(t)\end{equation}
and $\frac{1}{V}\int_{M}e^{-u_{\varepsilon}(t)} dV_{\varepsilon t}=1$, where $\theta_{\varepsilon}=(1-\beta)(\omega_{0}+\sqrt{-1}\partial\overline{\partial}\log(\varepsilon^{2}+|s|_{h}^{2}))$. Then for any $\beta\in(0,1)$, there exists a uniform constant $C$, such that
\begin{eqnarray*}|R(g_{\varepsilon}(t))-tr_{g_{\varepsilon}(t)}\theta_{\varepsilon}|&\leq& C,\\
\|u_{\varepsilon}(t)\|_{C^{1}(g_{\varepsilon}(t))}&\leq& C,\\
diam(M,g_{\varepsilon}(t))&\leq&C\end{eqnarray*}
hold for any $t\geq1$ and $\varepsilon$, where $R(g_{\varepsilon}(t))$ and $diam(M,g_{\varepsilon}(t))$ are the scalar curvature and diameter of the manifold respectively with respect to the metric $g_{\varepsilon}(t)$.
\end{thm}

Now we start to prove Theorem \ref{1.8.2}. Firstly, through differentiating equation $(\ref{4.1})$ and $\frac{1}{V}\int_{M}e^{-u_{\varepsilon}(t)}dV_{\varepsilon t}=1$, we conclude
\begin{equation}\label{3.22.6}\frac{d}{dt}u_{\varepsilon}(t)=\triangle_{\omega_\varepsilon(t)} u_{\varepsilon}(t)+\beta u_{\varepsilon}(t)-a_{\varepsilon}(t),\end{equation}
where
\begin{equation}\label{3.22.7}a_{\varepsilon}(t)=\frac{\beta}{V}\int_{M}u_{\varepsilon}(t)e^{-u_{\varepsilon}(t)}dV_{\varepsilon t}.\end{equation}

It is obvious that $a_{\varepsilon}(t)\leq0$ by Jensen's inequality. When $\beta\in(0, \frac{1}{2}]$, by the analogous argument in \cite{NSGT} or \cite{JWL}, the lower bound of $a_\varepsilon(t)$ can be derived by using the functional $\mu_{\theta_\varepsilon}(g_\varepsilon, 1)$, because the term $\max\limits_{M}(R(g_{\varepsilon})-tr_{g_{\varepsilon}}\theta_{\varepsilon})^{-}$ in lower bound of $\mu_{\theta_\varepsilon}(g_\varepsilon, 1)$ can be uniformly bounded when $\beta\in(0, \frac{1}{2}]$. However, this method does not work when $\beta\in(\frac{1}{2}, 1)$. Here, for any $\beta\in(0,1)$, we use the uniform Poincar\'e inequality to get a uniform lower bound of $a_\varepsilon(t)$. This lower bound is independent of the lower bound of $\mu_{\theta_\varepsilon}(g_\varepsilon, 1)$.

\begin{lem}\label{1.8.3} Let $u_{\varepsilon}(t)$ satisfy $(\ref{4.1})$. Then for every $f\in C^{\infty}(M)$, we have inequality
\begin{eqnarray}\label{3.22.37'}\\ \nonumber
\frac{1}{V}\int_{M}f^{2}e^{-u_{\varepsilon}(t)}dV_{\varepsilon t}\leq\frac{1}{\beta V}\int_{M}|\overline{\nabla}f|_{g_\varepsilon(t)}^{2}e^{-u_{\varepsilon}(t)}dV_{\varepsilon t}+(\frac{1}{V}\int_{M}fe^{-u_{\varepsilon}(t)}dV_{\varepsilon t})^{2}.
\end{eqnarray}
\end{lem}

\medskip

{\bf Proof:}\ \ It suffices to show the lowest strictly positive eigenvalue $\mu$ of operator $L$ satisfying $\mu\geq\beta$, where
\begin{eqnarray}\label{3.22.38}Lf=-g^{i\bar{j}}_{\varepsilon}(t)\nabla_{i}\nabla_{\bar{j}}f+g^{i\bar{j}}_{\varepsilon}(t)\nabla_{i}u_{\varepsilon}(t)\nabla_{\bar{j}}f.
\end{eqnarray}
Note that $L$ is self-adjoint with respect to the inner product
\begin{eqnarray}\label{3.22.39}(f,g)=\frac{1}{V}\int_{M}f\bar{g}e^{-u_{\varepsilon}(t)}dV_{\varepsilon t},
\end{eqnarray}
and $Ker\ L=\mathbb{C}$.
Suppose that $f$ is the eigenfunction of eigenvalue $\mu$, $f \not \equiv Constant$.
$$-g^{i\bar{j}}_{\varepsilon}(t)\nabla_{i}\nabla_{\bar{j}}f+g^{i\bar{j}}_{\varepsilon}(t)\nabla_{i}u_{\varepsilon}(t)\nabla_{\bar{j}}f=\mu f.
$$
By applying $\nabla_{\bar{k}}$ on both sides and combining Ricci identity, we have
$$\mu\nabla_{\bar{k}} f=-g^{i\bar{j}}_{\varepsilon}(t)\nabla_{i}\nabla_{\bar{k}}\nabla_{\bar{j}}f-g^{i\bar{j}}_{\varepsilon}(t)R^{\bar{s}}_{g_\varepsilon(t)\bar{j}i\bar{k}}\nabla_{\bar{s}}f+
g^{i\bar{j}}_{\varepsilon}(t)\nabla_{i}u_{\varepsilon}(t)\nabla_{\bar{k}}\nabla_{\bar{j}}f+g^{i\bar{j}}_{\varepsilon}(t)\nabla_{\bar{j}}f\nabla_{\bar{k}}\nabla_{i}u_{\varepsilon}(t).
$$
Integrating after multiplying $g^{l\bar{k}}_{\varepsilon}(t)\nabla_{l}fe^{-u_{\varepsilon}}dV_{\varepsilon t}$ on both sides, and then using the facts that $-Ric(\omega_\varepsilon(t))+\beta\omega_\varepsilon(t)+\theta_\varepsilon=\sqrt{-1}\partial\bar{\partial}u_\varepsilon(t)$ and $\theta_{\varepsilon}$ is semi-positive, we get
\begin{eqnarray*}\beta\int_{M}|\overline{\nabla}f|_{g_\varepsilon(t)}^{2}e^{-u_{\varepsilon}(t)}d V_{\varepsilon t} &\leq&\int_{M}|\overline{\nabla}\overline{\nabla}f|_{g_\varepsilon(t)}^{2}e^{-u_{\varepsilon}(t)}d V_{\varepsilon t}+\beta\int_{M}|\overline{\nabla}f|_{g_\varepsilon(t)}^{2}e^{-u_{\varepsilon}(t)}dV_{\varepsilon t}\\
&\ &+\frac{1}{2}\int_{M}\theta_{\varepsilon}(grad\ f,\mathcal{J}(grad\ f))e^{-u_{\varepsilon}(t)}dV_{\varepsilon t}\\
&=&\mu\int_{M}|\overline{\nabla}f|_{g_\varepsilon(t)}^{2}e^{-u_{\varepsilon}(t)}d V_{\varepsilon t}.
\end{eqnarray*}
Hence $\mu\geq\beta$.\QEDB

\begin{lem}\label{1.8.4} There exists a uniform constant $C$, such that
\begin{equation}\label{3.22.8}|a_{\varepsilon}(t)|\leq C\end{equation}
for any $t$ and $\varepsilon$.
\end{lem}

{\bf Proof:}\ \ We only need to prove that $a_{\varepsilon}(t)$ can be uniformly bounded from below. By Lemma \ref{1.8.3}, (\ref{3.22.6}) and (\ref{3.22.7}), we have
\begin{eqnarray*}\frac{d}{dt}a_\varepsilon(t)&=&\frac{\beta}{V}\int_M |\nabla u_{\varepsilon}(t)|_{\omega_\varepsilon(t)}^{2}e^{-u_{\varepsilon}(t)}dV_{\varepsilon t}-\frac{\beta^2}{V}\int_M u_{\varepsilon}^{2}(t)e^{-u_{\varepsilon}(t)}dV_{\varepsilon t}+a_\varepsilon^{2}(t)\\
&\geq&0.
\end{eqnarray*}
So $a_\varepsilon(t)$ is nondecreasing in time $t$ and $a_\varepsilon(t)\geq a_\varepsilon(0)=\frac{\beta}{V}\int_M  u_{\varepsilon}e^{-u_{\varepsilon}}dV_{\varepsilon0}$. For any $\beta\in(0,1)$, $\log\frac{\omega_{\varepsilon}^{n}(\varepsilon^{2}+|s|_{h}^{2})^{1-\beta}}{\omega_{0}^{n}}+k\beta\chi+F_{0}$ can be uniformly bounded (see section $4.5$ in \cite{CGP}). Hence by adjusting $\chi$ with a constant (whose variation with respect to $\varepsilon$ is bounded), we can assume that
\begin{eqnarray}\label{1.88.44}u_{\varepsilon}=\log\frac{\omega_{\varepsilon}^{n}(\varepsilon^{2}+|s|_{h}^{2})^{1-\beta}}{\omega_{0}^{n}}+k\beta\chi+F_{0}.
\end{eqnarray}
Then there exists a uniform constant $C$ such that (\ref{3.22.8}) holds.
\QEDB

\begin{pro}\label{1.8.5}  The twisted Ricci potential $u_{\varepsilon}(t)$ is uniformly bounded from below along the flow $(\ref{GKRF1})$.
\end{pro}

{\bf Proof:}\ \ By equation $(\ref{4.1})$, we have $\triangle_{g_{\varepsilon}(t)} u_{\varepsilon}(t)=-R(g_{\varepsilon}(t))+\beta n+tr_{g_{\varepsilon}(t)}\theta_{\varepsilon}$. From Proposition \ref{1.8.1} and Lemma \ref{1.8.4}, when $t\geq1$, there exists a uniform constant $C_{1}$ satisfying
\begin{eqnarray}\label{3.22.10}\triangle_{g_{\varepsilon}(t)} u_{\varepsilon}(t)-a_{\varepsilon}(t)\leq C_{1}.\end{eqnarray}

We conjecture $u_{\varepsilon}(t)\geq-\frac{2C_{1}}{\beta}$ for any $t\geq1$ and $\varepsilon$. If not, then there exists $(\varepsilon_{0},\ y_{0},\ t_{0})$, where $t_0\geq1$, such that $u_{\varepsilon_{0}}(y_{0},t_{0})<-\frac{2C_{1}}{\beta}$. By $(\ref{3.22.6})$
$$\frac{du_{\varepsilon_{0}}(t)}{dt}|_{(y_{0},t_{0})}\leq \beta u_{\varepsilon_{0}}(y_{0},t_{0})+C_{1}<-C_{1}.$$
So there exists $U(y_{0})\times[t_{0},t_{0}+\delta)$ such that $u_{\varepsilon_{0}}(t)$ satisfies
\[
\left\{
  \begin{array}{lll}
  \frac{du_{\varepsilon_{0}}(t)}{dt}|_{(y,t)}<0\\
  \ \ \ \ \ \ \ \ \ \ \ \ \ \ \ \ \ \ \ \ \ \ \ \ \ \ \ \ \ \ \ \ \ \ \ \ \ \ \ \ (y,t)\in U(y_{0})\times[t_{0},t_{0}+\delta).\\
  u_{\varepsilon_{0}}(y,t)<-\frac{2C_{1}}{\beta}\\
  \end{array}
  \right.
\]
By the continuity of $u_{\varepsilon_{0}}(t)$ with respect to time $t$, $u_{\varepsilon_{0}}(y)\ll 0$ on $U(y_{0})$ when $t\geq t_{0}$. Now we denote $U(y_{0})$ as $U$ for simplicity.

For any $z\in U,t\geq t_{0}$, $\frac{du_{\varepsilon_{0}}(t)}{dt}|_{(z,t)}\leq \beta u_{\varepsilon_{0}}(z,t)+C_{1}$, so
\begin{equation}\label{3.22.11}u_{\varepsilon_{0}}(z,t)\leq e^{\beta t}(u_{\varepsilon_{0}}(z,t_{0})e^{-\beta t_{0}}-\frac{C_{1}}{\beta}e^{-\beta t}+\frac{C_{1}}{\beta}e^{-\beta t_{0}})\leq-C_{2}e^{\beta t},\end{equation}
where $C_{2}$ depends only on $C_{1}$, $\beta$ and $t_{0}$. From equation $(\ref{4.1})$ and the flow $(\ref{GKRF1})$, we have
\begin{eqnarray}\label{3.22.111}u_{\varepsilon}(t)=\frac{d}{dt}(\varphi_{\varepsilon}(t)-e^{\beta t}\varphi_{\varepsilon}(0))+\tilde{c}_{\varepsilon}(t),
\end{eqnarray}
where $\tilde{c}_{\varepsilon}(t)$ depends on $\varepsilon$ and $t$.
$$u_{\varepsilon_{0}}(t)=\frac{d}{dt}(\varphi_{\varepsilon_0}(t)-e^{\beta t}\varphi_{\varepsilon_0}(0))+\tilde{c}_{\varepsilon_{0}}(t)=\frac{d}{dt}(\varphi_{\varepsilon_{0}}(t)
-e^{\beta t}\varphi_{\varepsilon_0}(0)+\int_{0}^{t}\tilde{c}_{\varepsilon_{0}}(s)ds)=\frac{d}{dt}\phi_{\varepsilon_{0}}(t).$$
When $z\in U$ and $t$ is sufficiently large,
\begin{equation}\label{3.22.12}\phi_{\varepsilon_{0}}(z,t)\leq\phi_{\varepsilon_{0}}(z,t_{0})-\frac{C_{2}}{\beta}e^{\beta t}+\frac{C_{2}}{\beta}e^{\beta t_{0}}\leq-C_{3}e^{\beta t},
\end{equation}
where $C_{3}$ depends only on $C_{1}$, $\beta$, $t_{0}$ and $\varepsilon_{0}$.

On the other hand,
$$1=\frac{1}{V}\int_{M}e^{-u_{\varepsilon_{0}}(t)}dV_{_{\varepsilon_{0}t}}\geq e^{-\sup\limits_{M}u_{\varepsilon_{0}}(t)}$$
for any $t$, hence we have
\begin{equation}\label{3.22.13}\sup_{M}u_{\varepsilon_{0}}(t)\geq 0.\end{equation}
Let $u_{\varepsilon_{0}}(x_{t},t)=\sup\limits_{M}u_{\varepsilon_{0}}(t)$. Combining $(\ref{3.22.13})$ with

$$\frac{d}{dt}(u_{\varepsilon_{0}}(t)-\beta\phi_{\varepsilon_{0}}(t))=\triangle_{g_{\varepsilon_{0}}(t)} u_{\varepsilon_{0}}(t)-a_{\varepsilon_{0}}(t)\leq C_{1},$$
we have
$$u_{\varepsilon_{0}}(x_{t},t)-\beta\phi_{\varepsilon_{0}}(x_{t},t)-(u_{\varepsilon_{0}}(x_{t},0)-\beta\phi_{\varepsilon_{0}}(x_{t},0))\leq C_{1}t,$$
\begin{equation}\label{3.22.14}\sup_{M}\phi_{\varepsilon_{0}}(\cdot,t)\geq-C_{4}-C_{1}t,
\end{equation}
where $C_{4}$ depends only on $\beta$ and $\omega_0$. Applying the Green formula with respect to metric $g_{0}$, for sufficiently large $t$, we have
\begin{eqnarray*}\phi_{\varepsilon_{0}}(x,t)+k\chi(\varepsilon_{0}^{2}+|s|^{2}_{h})&=&\frac{1}{Vol_{0}(M)}\int_{M}\phi_{\varepsilon_{0}}(y,t)+k\chi(\varepsilon_{0}^{2}+|s|^{2}_{h})dV_{0}\\
&\ &-\frac{1}{Vol_{0}(M)}\int_{M}\triangle_{g(0)}(\phi_{\varepsilon_{0}}(y,t)+k\chi(\varepsilon_{0}^{2}+|s|^{2}_{h}))G_{0}(x,y)dV_{0}\\
&\leq&\frac{Vol_{0}(M\setminus U)}{Vol_{0}(M)}\sup_{M}\phi_{\varepsilon_{0}}(\cdot,t)+\frac{1}{Vol_{0}(M)}\int_{U}\phi_{\varepsilon_{0}}(y,t)dV_{0}+C_{5}\\
&\leq&\frac{Vol_{0}(M\setminus U)}{Vol_{0}(M)}\sup_{M}\phi_{\varepsilon_{0}}(\cdot,t)-C_{6}e^{\beta t}+C_{5}.
\end{eqnarray*}
Then we obtain
\begin{equation}\label{3.22.15}\sup_{M}\phi_{\varepsilon_{0}}(\cdot,t)\leq-C_{7}e^{\beta t}+C_{8},
\end{equation}
where $C_{7}$ depends only on $C_{1}$,  $\beta$, $\varepsilon_{0}$, $t_{0}$ and $\omega_0$ while $C_{8}$ depends only on $\omega_0$.
From $(\ref{3.22.14})$ and $(\ref{3.22.15})$,
\begin{equation}\label{4.3.3'}-C_{7}e^{\beta t}+C_{8}\geq-C_{4}-C_{1}t.\end{equation}
When $t$ is sufficiently large, the inequality $(\ref{4.3.3'})$ is not correct, so $u_{\varepsilon}(t)$ is bounded uniformly from below along the flow $(\ref{GKRF1})$ when $t\geq 1$. By equation (\ref{3.22.6}) and $a_\varepsilon(t)\leq0$, we have
$$(\frac{d}{dt}-\triangle_{g_{\varepsilon}(t)})u_{\varepsilon}(t)\geq\beta u_{\varepsilon}(t).$$
Applying the maximum principle and the uniform bound of $u_\varepsilon$ (see (\ref{1.88.44})), we deduce that $u_\varepsilon(t)$ is uniformly bounded from below on $[0,2]\times M$. \QEDB

\medskip

Denote $\Box=\frac{d}{dt}-\triangle_{g_{\varepsilon}(t)}$, as the computations in \cite{JWL}, we have
\begin{equation}\label{3.22.16}\Box(\triangle_{g_{\varepsilon}(t)} u_{\varepsilon}(t))=-|\nabla\overline{\nabla}u_{\varepsilon}(t)|_{g_{\varepsilon}(t)}^{2}+\beta\triangle_{g_{\varepsilon}(t)} u_{\varepsilon}(t),\end{equation}
\begin{eqnarray}\label{3.22.17}\nonumber\Box(|\nabla u_{\varepsilon}(t)|_{g_{\varepsilon}(t)}^{2})&=&-|\nabla\overline{\nabla}u_{\varepsilon}(t)|_{g_{\varepsilon}(t)}^{2}-|\nabla\nabla u_{\varepsilon}(t)|_{g_{\varepsilon}(t)}^{2}\\
&\ &+\beta|\nabla u_{\varepsilon}(t)|_{g_{\varepsilon}(t)}^{2}-\frac{1}{2}\theta_{\varepsilon}(grad\ u_{\varepsilon}(t),\mathcal{J}(grad\ u_{\varepsilon}(t))),\end{eqnarray}
where $\mathcal{J}$ is the complex structure on $M$.

\begin{lem}\label{1.8.6}  For any $\beta\in(0,1)$, there exists a uniform constant $C$ independent of $t$ and $\varepsilon$, such that
\begin{eqnarray}\label{1.5.2}\|u_\varepsilon(t)\|_{C^1(g_{\varepsilon}(t))}&\leq C,\\
\label{1.5.3}|R(g_{\varepsilon}(t))-tr_{g_{\varepsilon}(t)}\theta_{\varepsilon}|&\leq C.\end{eqnarray}
on $[1,2]\times M$.
\end{lem}

{\bf Proof:}\ \ In Proposition \ref{1.8.1} and Proposition \ref{1.8.5}, we have got the lower bound of $u_\varepsilon(t)$ and $R(g_{\varepsilon}(t))-tr_{g_{\varepsilon}(t)}\theta_{\varepsilon}$ on $[1,2]\times M$. Since $a_\varepsilon(t)\geq-C$ for some uniform $C$, then by equation (\ref{3.22.6}), we have
$$(\frac{d}{dt}-\triangle_{g_{\varepsilon}(t)})(e^{-\beta t}u_{\varepsilon}(t)+\frac{C}{\beta}e^{-\beta t})\leq0.$$
By maximum principle and the uniform bound of $u_\varepsilon$, it follows that
$$u_\varepsilon(t)\leq e^{\beta t}(u_\varepsilon+\frac{C}{\beta})\leq C$$
on $[0,2]\times M$.

Let $H_\varepsilon(t,x)=t|\nabla u_\varepsilon(t)|^{2}_{g_{\varepsilon}(t)}+Au^{2}_\varepsilon(t)$ and $(t_0, x_0)$ be the maximum point of $H_\varepsilon(t,x)$ on $[0,2]\times M$. By (\ref{3.22.6}) and (\ref{3.22.17}), we obtain
\begin{eqnarray}(\frac{d}{dt}-\triangle_{g_{\varepsilon}(t)})H_\varepsilon(t,x)\leq(2\beta+1-2A)|\nabla u_\varepsilon(t)|^{2}_{g_{\varepsilon}(t)}+C,\end{eqnarray}
where constant $A$ will be determined later, and $C$ is a uniform constant depending only on $\|u_\varepsilon(t)\|_{C^0([0,2]\times M)}$ and $\sup\limits_{[0,2]}|a_\varepsilon(t)|$,

$Case\ 1$, $t_0=0$. Then $t|\nabla u_\varepsilon(t)|^{2}_{g_{\varepsilon}(t)}\leq Au^{2}_\varepsilon(0, x_0)\leq C,$ where $C$ is a uniform constant.

$Case\ 2$, $t_0>0$. Let $A=\beta+1$. By maximum principle, we have $|\nabla u_\varepsilon(t_0, x_0)|^{2}_{g_\varepsilon(t_0)}\leq C.$ Hence $t|\nabla u_\varepsilon(t)|^{2}_{g_\varepsilon(t)}\leq C.$

By the above two cases, we conclude that $t|\nabla u_\varepsilon(t)|^{2}_{g_\varepsilon(t)}\leq C$ on $[0,2]\times M$. Obviously $|\nabla u_\varepsilon(t)|^{2}_{g_\varepsilon(t)}\leq C$ on $[1,2]\times M$.

Next, since $\triangle_{g_{\varepsilon}(t)} u_{\varepsilon}(t)=-R(g_{\varepsilon}(t))+\beta n+tr_{g_{\varepsilon}(t)}\theta_{\varepsilon}$, we only need to prove the uniform upper bound of $-\triangle_{g_{\varepsilon}(t)} u_{\varepsilon}(t)$. We take $G_\varepsilon(t,x)=t^2(-\triangle_{g_{\varepsilon}(t)} u_{\varepsilon}(t))+2t^2|\nabla u_\varepsilon(t)|^{2}_{g_\varepsilon(t)}$. According to (\ref{3.22.16}), (\ref{3.22.17}) and
$$|\nabla \bar{\nabla}u_\varepsilon(t)|^{2}_{g_\varepsilon(t)}\geq\frac{(\triangle_{g_{\varepsilon}(t)} u_{\varepsilon}(t))^2}{n},$$
the evolution equation of $G_\varepsilon(t,x)$ can be written as
\begin{eqnarray*}(\frac{d}{dt}-\triangle_{g_\varepsilon(t)})G_\varepsilon(t,x)\leq(\beta t^2+2t)(-\triangle_{g_{\varepsilon}(t)} u_{\varepsilon}(t))-\frac{t^2}{n}(\triangle_{g_{\varepsilon}(t)} u_{\varepsilon}(t))^2+C\end{eqnarray*}
for some uniform constant $C$ depending on $\sup\limits_{[0,2]\times M}(t|\nabla u_\varepsilon(t)|^{2}_{g_\varepsilon(t)})$.
Assuming that $(t_0, x_0)$ is the maximum point of $G_\varepsilon(t,x)$ on $[0,2]\times M$:

$Case\ 1$, $t_0=0$, then $-t^2\triangle_{g_{\varepsilon}(t)} u_{\varepsilon}(t)\leq 0$.

$Case\ 2$, $t_0>0$. We assume $-\triangle_{g_{\varepsilon}(t)} u_{\varepsilon}(t)> 0$ at $(t_0,x_0)$. Then we claim that $t_{0}^{2}(-\triangle_{g_{\varepsilon}(t_0)} u_{\varepsilon}(t_0,x_0))\leq Bn$, where $B$ is a uniform constant to be determined later. If not,
$t_{0}^{2}(-\triangle_{g_{\varepsilon}(t_0,)} u_{\varepsilon}(t_0,x_0))> Bn$. By maximum principle, we have
\begin{eqnarray*}0&\leq&(-\triangle_{g_{\varepsilon}(t_0)} u_{\varepsilon}(t_0,x_0))((4+4\beta)-\frac{1}{n}(-t_{0}^2\triangle_{g_{\varepsilon}(t_0)} u_{\varepsilon}(t_0,x_0)))+C\\
&<&\frac{Bn}{4}(4+4\beta-B)+C.
\end{eqnarray*}
We can get a contradiction when we let $B=4(1+\beta+C)$. From these two cases, we conclude that $-t^2\triangle_{g_{\varepsilon}(t)} u_{\varepsilon}(t)\leq C$ for some uniform constant on $[0,2]\times M$. Furthermore, $-\triangle_{g_{\varepsilon}(t)} u_{\varepsilon}(t)\leq C$ on $[1,2]\times M$.\QEDB

\begin{lem}\label{1.8.7} There exists a uniform constant $C$ independent of $t$ and $\varepsilon$, such that
\begin{eqnarray}\label{3.22.18}|\nabla u_{\varepsilon}(t)|_{g_{\varepsilon}(t)}^{2}&\leq& C(u_{\varepsilon}(t)+C),\\\label{3.22.19}R(g_{\varepsilon}(t))-tr_{g_{\varepsilon}(t)}\theta_{\varepsilon}&\leq& C(u_{\varepsilon}(t)+C)\end{eqnarray}
for any $t\geq1$ and $\varepsilon$.
\end{lem}

{\bf Proof:}\ \ It follows from Proposition \ref{1.8.5} that there exists a uniform constant $B>1$ such that $u_{\varepsilon}(t)>-B$. Define
\begin{eqnarray}\label{3.22.20}H_{\varepsilon}(t)=\frac{|\nabla u_{\varepsilon}(t)|_{g_\varepsilon(t)}^{2}}{u_{\varepsilon}(t)+2B}.\end{eqnarray}
As the same argument in \cite{JWL}, we have
\begin{eqnarray*}\Box H_{\varepsilon}(t)&\leq&\frac{-|\nabla\overline{\nabla}u_{\varepsilon}(t)|_{g_\varepsilon(t)}^{2}-|\nabla\nabla u_{\varepsilon}(t)|_{g_\varepsilon(t)}^2}{u_{\varepsilon}(t)+2B}+\frac{|\nabla u_{\varepsilon}(t)|_{g_\varepsilon(t)}^{2}(2B\beta+a_{\varepsilon}(t))}{(u_{\varepsilon}(t)+2B)^{2}}\\
&\ &+\frac{\delta}{2}\frac{|\nabla u_{\varepsilon}(t)|_{g_\varepsilon(t)}^4}{(u_{\varepsilon}(t)+2B)^{3}}+\delta\frac{|\nabla\nabla u_{\varepsilon}(t)|_{g_\varepsilon(t)}^{2}+|\nabla\overline{\nabla}u_{\varepsilon}(t)|_{g_\varepsilon(t)}^{2}}{u_{\varepsilon}(t)+2B}
-\delta\frac{|\nabla u_{\varepsilon}(t)|_{g_\varepsilon(t)}^4}{(u_{\varepsilon}(t)+2B)^{3}}\\
&\ &-\frac{1}{2}\frac{\theta_{\varepsilon}(grad\ u_{\varepsilon}(t),\mathcal{J}(grad\ u_{\varepsilon}(t)))}{u_{\varepsilon}(t)+2B}+(2-\delta)Re\frac{\overline{\nabla}u_{\varepsilon}(t)\cdot\nabla H_{\varepsilon}(t)}{u_{\varepsilon}(t)+2B}.
\end{eqnarray*}
Taking $\delta<1$ and combining Lemma \ref{1.8.4} with $\theta(grad\ u_{\varepsilon}(t),\mathcal{J}(grad\ u_{\varepsilon}(t)))\geq 0$, we obtain
$$\Box H_{\varepsilon}(t)\leq\frac{|\nabla u_{\varepsilon}(t)|_{g_\varepsilon(t)}^{2}(2B\beta+C_{1})}{(u_{\varepsilon}(t)+2B)^{2}}+(2-\delta)Re\frac{\overline{\nabla}u_{\varepsilon}(t)\cdot\nabla H_{\varepsilon}(t)}{u_{\varepsilon}(t)+2B}-\frac{\delta}{2}\frac{|\nabla u_{\varepsilon}(t)|_{g_\varepsilon(t)}^{4}}{(u_{\varepsilon}(t)+2B)^{3}}.$$

From Lemma \ref{1.8.6}, we have
\begin{eqnarray}\label{1.5.5}\sup\limits_{M}\frac{|\nabla u_{\varepsilon}(1)|_{g_{\varepsilon}(1)}^{2}}{u_{\varepsilon}(1)+2B}\leq C_{3}\end{eqnarray}
for some uniform constant $C_{3}$. Then by maximum principle, we have $H_{\varepsilon}(t)\leq \max\{C_{3},\ 2(2B+C_{1})\delta^{-1}\}$ for any $t\geq 1$ and $\varepsilon$.

Now we prove the second inequality. Since $\triangle_{g_{\varepsilon}(t)} u_{\varepsilon}(t)=\beta n-R(g_{\varepsilon}(t))+tr_{g_{\varepsilon}(t)}\theta_{\varepsilon}$, we only need to prove the existence of the uniform constant $C$ such that $-\triangle_{g_{\varepsilon}(t)} u_{\varepsilon}(t)$ can be controlled by $C(u_{\varepsilon}(t)+C)$.\\
Let $G_{\varepsilon}=\frac{-\triangle_{g_{\varepsilon}(t)} u_{\varepsilon}(t)}{u_{\varepsilon}(t)+2B}+2H_{\varepsilon}$.
\begin{eqnarray*}\Box G_{\varepsilon}&=&\frac{-2|\nabla\nabla u_{\varepsilon}(t)|_{g_\varepsilon(t)}^{2}-|\nabla\overline{\nabla}u_{\varepsilon}(t)|_{g_\varepsilon(t)}^{2}}{u_{\varepsilon}(t)+2B}+\frac{(-\triangle_{g_{\varepsilon}(t)} u_{\varepsilon}(t)+2|\nabla u_{\varepsilon}(t)|_{g_\varepsilon(t)}^{2})(2B\beta+a_{\varepsilon}(t))}{(u_{\varepsilon}(t)+2B)^{2}}
\\
&\ &+2Re\frac{\overline{\nabla}u_{\varepsilon}(t)\cdot\nabla G_{\varepsilon}}{u_{\varepsilon}(t)+2B}-\frac{1}{2}\frac{\theta_{\varepsilon}(grad\ u_{\varepsilon}(t),\mathcal{J}(grad\ u_{\varepsilon}(t)))}{u_{\varepsilon}(t)+2B}.\end{eqnarray*}
Since $\theta_{\varepsilon}$ is semi-positive,
$$\Box G_{\varepsilon}\leq\frac{-|\nabla\overline{\nabla}u_{\varepsilon}(t)|_{g_\varepsilon(t)}^{2}}{u_{\varepsilon}(t)+2B}+\frac{(-\triangle_{g_{\varepsilon}(t)} u_{\varepsilon}(t)+2|\nabla u_{\varepsilon}(t)|_{g_\varepsilon(t)}^{2})(2B\beta+a_{\varepsilon}(t))}{(u_{\varepsilon}(t)+2B)^{2}}+2\frac{\overline{\nabla}u_{\varepsilon}(t)\cdot\nabla G_{\varepsilon}}{u_{\varepsilon}(t)+2B}.$$
In local coordinates,
\begin{eqnarray}\label{3.22.25}(\triangle_{g_\varepsilon(t)} u_{\varepsilon}(t))^{2}=(\sum_{i}u_{\varepsilon i\bar{i}})^{2}\leq n\sum_{i}u_{\varepsilon i\bar{i}}^{2}=n|\nabla\overline{\nabla}u_{\varepsilon}(t)|_{g_\varepsilon(t)}^{2}.\end{eqnarray}
So we have
$$\Box G_{\varepsilon}
\leq\frac{-\triangle_{g_\varepsilon(t)} u_{\varepsilon}(t)}{u_{\varepsilon}(t)+2B}(\frac{2B\beta+a_{\varepsilon}(t)}{u_{\varepsilon}(t)+2B}-\frac{-\triangle_{g_{\varepsilon}(t)} u_{\varepsilon}(t)}{n(u_{\varepsilon}(t)+2B)})+2\frac{|\nabla u_{\varepsilon}(t)|_{g_\varepsilon(t)}^{2}(2B\beta+a_{\varepsilon}(t))}{(u_{\varepsilon}(t)+2B)^{2}}+2\frac{\overline{\nabla}u_{\varepsilon}(t)\cdot\nabla G_{\varepsilon}}{u_{\varepsilon}(t)+2B}.$$
Since $\frac{-\triangle_{g_{\varepsilon}(1)} u_{\varepsilon}(1)}{u_{\varepsilon}(1)+2B}$ is bounded uniformly from the argument in Lemma \ref{1.8.6}, there exists a uniform constant $C>0$ by maximum principle, such that $G_{\varepsilon}\leq C$ for any $t\geq 1$ and $\varepsilon$. Hence we get $\frac{-\triangle u_{\varepsilon}(t)}{u_{\varepsilon}(t)+2B}\leq G_{\varepsilon}\leq C$.\QEDB

\medskip

From $(\ref{3.22.18})$ in Lemma \ref{1.8.7} and the same discussion in \cite{NSGT} (see Claim $8$), we have the following lemma.

\begin{lem}\label{1.8.8} There exists a uniform constant $C$, such that
\begin{eqnarray}\label{3.22.26}
u_{\varepsilon}(y,t)&\leq&\ C\ dist_{\varepsilon t}^{2}(x,y)+C,\\
R(g_{\varepsilon}(t))-tr_{g_{\varepsilon}(t)}\theta_{\varepsilon}&\leq&\ C\ dist_{\varepsilon t}^{2}(x,y)+C,\\
|\nabla u_{\varepsilon}(t)|_{g_{\varepsilon}(t)}&\leq&\ C\ dist_{\varepsilon t}(x,y)+C,\end{eqnarray}
for any $t\geq 1$ and $\varepsilon$, where $u_{\varepsilon}(x,t)=\inf\limits_{y\in M}u_{\varepsilon}(y,t)$.
\end{lem}

By Lemma \ref{1.8.8}, the statements in Theorem \ref{1.8.2} will be true if the $diam(M,g_{\varepsilon}(t))$ is uniformly bounded when $t\geq1$. In order to prove this, we will give a proof of a twisted version of uniform Perelman's noncollapsing theorem by the argument in \cite{KL} and \cite{STW}. Before this, we review the twisted $\mathcal{W}_{\theta}$ functional and $\mu_{\theta}$\ functional.
$$\mathcal{W}_{\theta}(g,f,\tau)=\int_{M}e^{-f}\tau^{-n}(\tau(R-tr_{g}\theta+|\nabla f|_{g}^{2})+\beta f)dV_{g},$$
where $g$ is a K$\ddot{a}$hler metric, $f$ is a smooth function on $M$, $\tau$ is a positive scale parameter and $n$ is the complex dimension of the K$\ddot{a}$hler manifold. Let
$$\mu_{\theta}(g,\tau)=\inf\{\mathcal{W}_{\theta}(g,f,\tau)|f\in C^{\infty}(M),\frac{1}{V}\int_{M}e^{-f}\tau^{-n}dV=1\}$$
be the $\mu_{\theta}$ functional with respect to the metric $g$. From \cite{JWL}, we have the monotonicity of the twisted $\mathcal{W}_{\theta}$ and $\mu_{\theta}$ functional along the twisted K\"ahler-Ricci flow.

\begin{lem}\label{1.8.9}(Theorem 2.4 in \cite{JWL})  Along the evolution equation
\begin{equation}
\begin{cases}
  \frac{\partial g_{i\bar{j}}}{\partial t}=-R_{i\bar{j}}+\beta g_{i\bar{j}}+\theta_{i\bar{j}}\\
  \frac{\partial f}{\partial t}\ =\beta n\tau^{-1}-R+tr_{g(t)}\theta-\triangle f+|\nabla f|^{2}\\
  \frac{\partial\tau}{\partial t}\ =\beta(\tau-1)
  \end{cases},
  \end{equation}
$\mathcal{W}_{\theta}(g(t),f(t),\tau(t))$ is nondecreasing.
\end{lem}

\begin{lem}\label{1.8.10} (Theorem 2.5 in \cite{JWL}) If $\tau$ satisfies the following equality
$$\frac{\partial\tau}{\partial t}\ =\beta(\tau-1),$$
then $\mu_{\theta}(g,\tau)$ is nondecreasing along the twisted K$\ddot{a}$hler-Ricci flow.
\end{lem}

In the process of proving the uniform Perelman's noncollapsing theorem and the fact that $diam(M,g_{\varepsilon}(t))$ can be uniformly bounded, we use the lower bound of the functional $\mu_{\theta_{\varepsilon}}(g_{\varepsilon}(1),\tau)$, which depends on the Sobolev constant $\mathcal{C}_{S}(M,g_{\varepsilon}(1))$ and $\max\limits_{M}(R(g_{\varepsilon}(1))-tr_{g_{\varepsilon}(1)}\theta_{\varepsilon})^{-}$ with respect to the metric $g_{\varepsilon}(1)$. In the following proposition, we obtain a uniform control of the Sobolev constant $\mathcal{C}_{S}(M,g_{\varepsilon}(t))$ for any $t\in [0, 2]$ and $\varepsilon>0$. To the reader's convenience, we will present the proof of Proposition \ref{1.8.12} in the appendix.

\begin{pro}\label{1.8.12}
Let $g_{\varepsilon}(t)$ be a solution of the twisted K$\ddot{a}$hler Ricci flow $(\ref{GKRF1})$. Then there exists some uniform constant $C$, such that
\begin{eqnarray}\label{1.5.6}(\int_{M}v^{\frac{2n}{n-1}}dV_{\varepsilon t})^{\frac{n-1}{n}}\leq C(\int_{M}|dv|^{2}_{g_{\varepsilon}(t)}dV_{\varepsilon t}+\int_{M}|v|^{2}dV_{\varepsilon t})\end{eqnarray}
holds for any smooth function $v$ on $M$, $t\in[0,2]$ and $\varepsilon>0$.
\end{pro}

In \cite{JWL} (see Theorem $2.2$), we have
\begin{eqnarray}\label{3.22.4}\nonumber\mu_{\theta_{\varepsilon}}(g_{\varepsilon}(1),\tau)&\geq&-\tau V\max\limits_{M}(R(g_{\varepsilon}(1))-tr_{g_{\varepsilon}(1)}\theta_{\varepsilon})^{-}-n\beta V(\log 2+\alpha V^{-\frac{1}{n}}-1)\\
&\ &+\beta n v\log\alpha-\beta V\log V-\beta n v\log\tau,
\end{eqnarray}
where $V$ is the volume of $(M,g_{\varepsilon}(1))$ and $\alpha$ satisfies $4\tau\geq\beta n\alpha\mathcal{C}_{S}(M,g_{\varepsilon}(1))$.
Since $Vol(M,g_{\varepsilon}(1))$ is fixed while $\max\limits_{M}(R(g_{\varepsilon}(1))-tr_{g_{\varepsilon}(1)}\theta_{\varepsilon})^{-}$ and $\mathcal{C}_{S}(M,g_{\varepsilon}(1))$ are uniformly bounded by Lemma $4.4$ and Proposition \ref{1.8.12}, we know that by choosing a suitable $\alpha$ there exists a uniform constant $C$ independent of $\varepsilon$, such that
\begin{eqnarray}\label{3.22.5}\mu_{\theta_{\varepsilon}}(g_{\varepsilon}(1),\tau)\geq-C.
\end{eqnarray}

Next, let us state the uniform Perelman's noncollapsing theorem and prove it.

\begin{pro}\label{1.8.13} Let $g_{\varepsilon}(t)$ be a solution of the flow $(\ref{GKRF1})$, there exists a uniform constant $C$, such that
$$Vol_{g_{\varepsilon}(t)}(B_{g_{\varepsilon}(t)}(x,r))\geq Cr^{2n}$$
for every $g_{\varepsilon}(t)$ satisfying $R(g_{\varepsilon}(t))-tr_{g_{\varepsilon}(t)}\theta_{\varepsilon}\leq\frac{m}{r^{2}}$ on $B_{g_{\varepsilon}(t)}(x,r)$ when $t\geq1$, where $\partial B_{g_{\varepsilon}(t)}(x,r)\neq\emptyset$ and $0<r<1$.
\end{pro}

{\bf Proof:}\ \  We argue it by contradiction, that is, there exist $\varepsilon_{k},\ p_{k},\ t_{k}\geq1,\ r_{k}$ satisfying $R_{g_{\varepsilon_{k}}(t_{k})}-tr_{g_{\varepsilon_{k}}(t_{k})}\theta_{\varepsilon_{k}}\leq\frac{m}{r_{k}^{2}}$ on $B_{g_{\varepsilon_k}(t_k)}(p_k,r_k)$, $\varepsilon_{k}\rightarrow 0$ and $Vol_{g_{\varepsilon_{k}}(t_{k})}(B_{g_{\varepsilon_{k}}(t_{k})}(p_{k},r_{k}))\cdot$
$r_{k}^{-2n}\rightarrow0$ when $k\rightarrow+\infty$. Define $B_{g_{\varepsilon_{k}}(t_{k})}(p_{k},r_{k})=B_{k}$, $Vol_{g_{\varepsilon_{k}}(t_{k})}(B_{g_{\varepsilon_{k}}(t_{k})}(p_{k},r_{k}))=V(r_{k})$ in the following.

Setting $\tau=r_{k}^{2}$ at $t_{k}$, we define function as
\begin{eqnarray}\label{3.22.26}u_{k}(x)=e^{C_{k}}\phi(r_{k}^{-1}dist_{g_{\varepsilon_{k}}(t_{k})}(x,p_{k})),\end{eqnarray}
where $\phi$ is smooth function on $\mathbb{R}$, equal to $1$ on $[0,\frac{1}{2}]$, decreasing
on $[\frac{1}{2},1]$ and equal to 0 on $[1,+\infty)$, $C_{k}$ is a constant to make $u_{k}$ satisfy the constraint
$$\frac{1}{V}\int_{M}r_{k}^{-2n}u_{k}^{2}dV_{\varepsilon_{k}t_{k}}=1.$$
Hence
$$1=\frac{1}{V}e^{2C_{k}}r_{k}^{-2n}\int_{B_{k}}\phi^{2}dV_{\varepsilon_{k}t_{k}}\leq\frac{1}{V}e^{2C_{k}}r_{k}^{-2n}V(r_{k}).$$
By assumption, $V(r_{k})r_{k}^{-2n}\rightarrow0$ when  $k\rightarrow+\infty$, which shows that $C_{k}\rightarrow+\infty$ when $k\rightarrow+\infty$.
So we claim that

$(a)\ V(r_{k})r_{k}^{-2n}\rightarrow0$ as $k\rightarrow+\infty$;

$(b)\ (R(g_{\varepsilon_{k}}(t_{k}))-tr_{g_{\varepsilon_{k}}(t_{k})}\theta_{\varepsilon_{k}})r_{k}^{2}\leq m$;

$(c)\ \frac{V(r_{k})}{V(\frac{r_{k}}{2})}$ is uniformly bounded.

We only need to prove $(c)$. If $\frac{V(r_{k})}{V(\frac{r_{k}}{2})}<5^{n}$ for any $k$, then the claim $(c)$ is testified. If not, for a given $k$, we have $\frac{V(r_{k})}{V(\frac{r_{k}}{2})}\geq5^{n}$.
Let $r_{k}^{\prime}=\frac{r_{k}}{2}$. We have $(r_{k}^{\prime})^{-2n}V(r_{k}^{\prime})\leq(r_{k}^{\prime})^{-2n}\frac{1}{5^{n}}V(r_{k})=(\frac{4}{5})^{n}
r_{k}^{-2n}V(r_{k}),\ (r_{k}^{\prime})^{2}(R(g(_{\varepsilon_{k}}(t_{k}))-tr_{g_{\varepsilon_{k}}(t_{k})}\theta_{\varepsilon_{k}})=\frac{r_{k}^{2}}{4}
(R(g_{\varepsilon_{k}}(t_{k}))-tr_{g(t_{k})}\theta_{\varepsilon_{k}})$. Combining $(a)$ and $(b)$, we obtain $(r_{k}^{\prime})^{2}(R(g_{\varepsilon_{k}}(t_{k}))-tr_{g_{\varepsilon_{k}}(t_{k})}\theta_{\varepsilon_{k}})\leq m$ and $(r_{k}^{\prime})^{-2n}V(r_{k}^{\prime}) \rightarrow 0$ when $k\rightarrow+\infty$. Replacing $r_{k}$ by $r_{k}^{\prime}$. If $\frac{V(r_{k}^{\prime})}{V(\frac{r_{k}^{\prime}}{2})}<5^{n}$, the demonstration will be terminated. If not, the above process will be repeated. By the identity $\lim_{r\rightarrow0}\frac{V(r)}{V(\frac{r}{2})}=4^{n}$ proved in \cite{WYH} (see ($6.9$)), we should get $\frac{V(r_{k}^{\prime})}{V(\frac{r_{k}^{\prime}}{2})}<5^{n}$ at some step. Then we consider $\{p_{k},r_{k}\}$ obtained from the above.\\
Considering the function $\frac{1}{r_{k}^{2}+1}(u_{k}^{2}+r_{k}^{2n+2})$, we have its integral average
$$\frac{1}{V}\int_{M}r_{k}^{-2n}\frac{1}{r_{k}^{2}+1}(u_{k}^{2}+r_{k}^{2n+2})dV_{\varepsilon_{k}t_{k}}=1.$$
Computing the functional $\mathcal{W}_{\theta_{\varepsilon_{k}}}(g_{\varepsilon_{k}}(t_{k}),-\log \frac{1}{r_{k}^{2}+1}(u_{k}^{2}+r_{k}^{2n+2}),r_{k}^{2})$, we have
\begin{eqnarray*}
&\ &\mathcal{W}_{\theta_{\varepsilon_{k}}}(g_{\varepsilon_{k}}(t_{k}),-\log \frac{1}{r_{k}^{2}+1}(u_{k}^{2}+r_{k}^{2n+2}),r_{k}^{2})\\
&=&\frac{1}{r_{k}^{2}+1}\int_{M}r_{k}^{-2n}(u_{k}^{2}+r_{k}^{2n+2})(R(g_{\varepsilon_{k}}(t_{k}))
-tr_{g_{\varepsilon_{k}}(t_{k})}\theta_{\varepsilon_{k}})r_{k}^{2}dV_{\varepsilon_{k}t_{k}}\ \ \ (1)\\
&\ &+\frac{1}{r_{k}^{2}+1}\int_{M}r_{k}^{-2n}(u_{k}^{2}+r_{k}^{2n+2})\frac{4u_{k}^{2}|\nabla u_{k}|_{g_{\varepsilon_{k}}(t_{k})}^{2}}{(u_{k}^{2}+r_{k}^{2n+2})^{2}}r_{k}^{2}dV_{\varepsilon_{k}t_{k}}\ \ \ \ \ \ \ \ \ \ \ \ \ \ (2)\\
&\ &-\frac{\beta}{r_{k}^{2}+1}\int_{M}r_{k}^{-2n}(u_{k}^{2}+r_{k}^{2n+2})\log\frac{1}{r_{k}^{2}+1}dV_{\varepsilon_{k}t_{k}}\ \ \ \ \ \ \ \ \ \ \ \ \ \ \ \ \ \ \ \ \ \ \ \ (3)\\
&\ &-\frac{\beta}{r_{k}^{2}+1}\int_{M}r_{k}^{-2n}(u_{k}^{2}+r_{k}^{2n+2})\log(u_{k}^{2}+r_{k}^{2n+2})dV_{\varepsilon_{k}t_{k}}.\ \ \ \ \ \ \ \ \ \ \ \ \ \ \ \ (4)\\
(1)&=&\frac{1}{r_{k}^{2}+1}\int_{M}r_{k}^{-2n}(u_{k}^{2}+r_{k}^{2n+2})(R(g_{\varepsilon_{k}}(t_{k}))
-tr_{g_{\varepsilon_{k}}(t_{k})}\theta_{\varepsilon_{k}})r_{k}^{2}dV_{\varepsilon_{k}t_{k}}\\
&\leq&\frac{m}{r_{k}^{2}+1}\int_{M}r_{k}^{-2n}(u_{k}^{2}+r_{k}^{2n+2})dV_{\varepsilon_{k}t_{k}}\leq m V,\\
(2)&=&\frac{1}{r_{k}^{2}+1}\int_{M}r_{k}^{-2n}(u_{k}^{2}+r_{k}^{2n+2})\frac{4u_{k}^{2}|\nabla u_{k}|_{g_{\varepsilon_{k}}(t_{k})}^{2}}{(u_{k}^{2}+r_{k}^{2n+2})^{2}}r_{k}^{2}dV_{\varepsilon_{k}t_{k}}\\
&\leq&\frac{n}{r_{k}^{2}+1}\int_{M}r_{k}^{-2n}e^{2C_{k}}4|\phi'|^{2}dV_{\varepsilon_{k}t_{k}}\\
&\leq&C r_{k}^{-2n}e^{2C_{k}}(V(r_{k})-V(\frac{r_{k}}{2})),\\
(3)&=&-\frac{\beta}{r_{k}^{2}+1}\int_{M}r_{k}^{-2n}(u_{k}^{2}+r_{k}^{2n+2})\log\frac{1}{r_{k}^{2}+1}dV_{\varepsilon_{k}t_{k}}\\
&=&\beta V\log(r_{k}^{2}+1)\leq\beta V\log2,\\
(4)&=&-\frac{\beta}{r_{k}^{2}+1}\int_{M}r_{k}^{-2n}(u_{k}^{2}+r_{k}^{2n+2})\log(u_{k}^{2}+r_{k}^{2n+2})dV_{\varepsilon_{k}t_{k}}\\
&\leq&-\frac{\beta}{r_{k}^{2}+1}\int_{M}r_{k}^{-2n}u_{k}^{2}\log u_{k}^{2}dV_{\varepsilon_{k}t_{k}}
-\frac{\beta(n+1)}{r_{k}^{2}+1}\int_{M}r_{k}^{2}\log r_{k}^{2}dV_{\varepsilon_{k}t_{k}}\\
\end{eqnarray*}
\begin{eqnarray*}
&\leq&-\frac{2\beta C_{k}}{r_{k}^{2}+1}\int_{M}r_{k}^{-2n}u_{k}^{2}dV_{\varepsilon_{k}t_{k}}
-\frac{\beta}{r_{k}^{2}+1}\int_{M}r_{k}^{-2n}e^{2C_{k}}\phi^{2}\log \phi^{2}dV_{\varepsilon_{k}t_{k}}
+\frac{C\beta(n+1)V}{r_{k}^{2}+1}\\
&\leq&-\beta VC_{k}+C r_{k}^{-2n}e^{2C_{k}}(V(r_{k})-V(\frac{r_{k}}{2}))+C,
\end{eqnarray*}
where the all above constants $C$ are all uniform. Combining all these inequalities together and making use of condition $(b)$ $(c)$, we have
\begin{eqnarray*}&\ &\mathcal{W}_{\theta_{\varepsilon_{k}}}(g_{\varepsilon_{k}}(t_{k}),-\log \frac{1}{r_{k}^{2}+1}(u_{k}^{2}+r_{k}^{2n+2}),r_{k}^{2})\\
&\leq&C-\beta V C_{k}+C r_{k}^{-2n}e^{2C_{k}}(V(r_{k})-V(\frac{r_{k}}{2}))\\
&\leq&C-C_{k}\beta V+Cr_{k}^{-2n}e^{2C_{k}}V(\frac{r_{k}}{2})\\
&=&C-C_{k}\beta V+C\int_{B_{g_{\varepsilon_{k}}(t_{k})}(p_{k},\frac{r_{k}}{2})}r_{k}^{-2n}e^{2C_{k}}\phi^{2}dV_{\varepsilon_{k}t_{k}}\\
&\leq&C-C_{k}\beta V,
\end{eqnarray*}
where $C$ is uniform constant independent of $t_{k}$ and $\varepsilon_{k}$. Considering $\tau=1-(1-r_{k}^{2})e^{-\beta t_{k}}e^{\beta t}$, by Lemma \ref{1.8.4}, we conclude
\begin{eqnarray*}\mu_{\theta_{\varepsilon_{k}}}(g_{\varepsilon_{k}}(1),1-(1-r_{k}^{2})e^{-\beta (t_{k}-1)})
&\leq&\mathcal{W}_{\theta_{\varepsilon_{k}}}(g_{\varepsilon_{k}}(t_{k}),-\log \frac{1}{r_{k}^{2}+1}(u_{k}^{2}+r_{k}^{2n+2}),r_{k}^{2})\\
&\leq& C-2C_{k}\beta V.
\end{eqnarray*}
Since $0<1-(1-r_{k}^{2})e^{-\beta (t_{k}-1)}<1$, we conclude by $(\ref{3.22.5})$ that
\begin{eqnarray}\label{3.22.27}\mu_{\theta_{\varepsilon_{k}}}(g_{\varepsilon_{k}}(1),1-(1-r_{k}^{2})e^{-\beta (t_{k}-1)})\geq -C,\end{eqnarray} where $C$ is independent of $\varepsilon_{k}$ and $t_{k}$. Then we get $-C\leq C-C_{k}\beta V$ which does not work when $k\rightarrow+\infty$. So the lemma is proved.\QEDB
\medskip

Denote $d_{\varepsilon t}(z)=dist_{\varepsilon t}(x,z)$, $B_{\varepsilon t}(k_{1},k_{2})=\{z|2^{k_{1}}\leq d_{\varepsilon t}(z)\leq 2^{k_{2}}\}$, where $u_{\varepsilon}(x,t)=\inf\limits_{M}u_{\varepsilon}(y,t)$. Considering an annulus $B_{\varepsilon t}(k,k+1)$, then by Lemma \ref{1.8.8}, we have $R(g_{\varepsilon}(t))-tr_{g_{\varepsilon}(t)}\theta_{\varepsilon}\leq C2^{2k}$ on $B_{\varepsilon t}(k,k+1)$ when $t\geq1$. Interval $[2^{k}, 2^{k+1}]$ fits
$2^{2k}$ balls of radii $\frac{1}{2^{k}}$. By Proposition \ref{1.8.13}, when $t\geq1$, we have
\begin{eqnarray}\label{4.17}Vol_{g_{\varepsilon}(t)}(B_{\varepsilon t}(k,k+1))\geq\sum_{i}Vol_{g_{\varepsilon}(t)}(B_{\varepsilon t}(x_{i},\frac{1}{2^{k}}))\geq C2^{2k-2nk}.\end{eqnarray}

\begin{lem}\label{1.8.14} When $t\geqslant 1$, for every $\delta>0$, there exists $B_{\varepsilon t}(k_{1},k_{2})$, such that if $diam(M,g_{\varepsilon}(t))$ is large enough, then
\begin{eqnarray*}(a)Vol_{g_{\varepsilon}(t)}(B_{\varepsilon t}(k_{1},k_{2}))&<&\delta,\\
(b)Vol_{g_{\varepsilon}(t)}(B_{\varepsilon t}(k_{1},k_{2}))&\leq& 2^{20n}Vol_{g_{\varepsilon}(t)}(B_{\varepsilon t}(k_{1}+2,k_{2}-2)).\end{eqnarray*}
\end{lem}

{\bf Proof:}\ \ First, we fix any $\delta>0$. Since $Vol_{g_{\varepsilon}(t)}(M)$ is a constant $V$ along the twisted K\"ahler-Ricci flow, it can be uniformly bounded. Let $k\gg1$.
$$V=Vol_{g_{\varepsilon}(t)}(B_{\varepsilon t}(0,k))+Vol_{g_{\varepsilon}(t)}(B_{\varepsilon t}(k,3k))+\cdots+Vol_{g_{\varepsilon}(t)}(B_{\varepsilon t}(3^{\alpha-1}k,3^{\alpha}k))+\cdots,$$
where $\alpha>m[\frac{V}{\delta}]+1$, $m$ will be determined later and $diam(M,g_{\varepsilon}(t))>2^{3^{\alpha}k+1}$. We claim that there must exist a $0\leq i\leq\alpha-1$, such that $Vol_{g_{\varepsilon}(t)}(B_{\varepsilon t}(3^{i}k,3^{i+1}k))<\delta$. If not, then we have
$$V>Vol_{g_{\varepsilon}(t)}(B_{\varepsilon t}(k,3k))+\cdots+Vol_{g_{\varepsilon}(t)}(B_{\varepsilon t}(3^{\alpha-1}k,3^{\alpha}k))\geq\alpha\delta>m\delta[\frac{V}{\delta}]+\delta.$$
When we take $m$ satisfying $m\delta[\frac{V}{\delta}]+\delta>V$, the above inequality leads to a contradiction. So the claim is proved.

Then we determine $k_{1}$ and $k_{2}$. If estimate $(b)$ does not hold, then
$$Vol_{g_{\varepsilon}(t)}(B_{\varepsilon t}(3^{i}k,3^{i+1}k))>2^{20n}Vol_{g_{\varepsilon}(t)}(B_{\varepsilon t}(3^{i}k+2,3^{i+1}k-2)).$$
We would consider $Vol_{g_{\varepsilon}(t)}(B_{\varepsilon t}(3^{i}k+2,3^{i+1}k-2))$ instead and discuss whether $(b)$ holds for
that ball. If for any $p$, at the p-th step we are still not able to find suitable radii to satisfy $(a)$ and $(b)$. In that case, at the p-th step we would have
$$Vol_{g_{\varepsilon}(t)}(B_{\varepsilon t}(3^{i}k,3^{i+1}k))>2^{20np}Vol_{g_{\varepsilon}(t)}(B_{\varepsilon t}(3^{i}k+2p,3^{i+1}k-2p)).$$
In particular, if $3^{i}k+2p=\frac{3}{2}3^{i}k$, then we have $3^{i+1}k-2p=\frac{5}{2}3^{i}k$. By $(\ref{4.17})$,
\begin{eqnarray*}\delta>Vol_{g_{\varepsilon}(t)}(B_{\varepsilon t}(3^{i}k,3^{i+1}k))&>&2^{5n\cdot3^{i}k}Vol_{g_{\varepsilon}(t)}(B_{\varepsilon t}(\frac{3}{2}3^{i}k,\frac{5}{2}3^{i}k))\\
&>&2^{5n\cdot3^{i}k}Vol_{g_{\varepsilon}(t)}(B_{\varepsilon t}(\frac{3}{2}3^{i}k,\frac{3}{2}3^{i}k+1))\\
&\geq&C2^{(2n+3)\cdot3^{i}k}.\end{eqnarray*}
This leads to contradiction if we let $k\gg1$. So there exists some $1\leq j\leq p-1$, such that
$$Vol_{g_{\varepsilon}(t)}(B_{\varepsilon t}(3^{i}k+2j,3^{i+1}k-2j))\leq2^{20n}Vol_{g_{\varepsilon}(t)}(B_{\varepsilon t}(3^{i}k+2(j+1),3^{i+1}k-2(j+1))).$$
Let $k_{1}=3^{i}k+2j$, $k_{2}=3^{i+1}k-2j$ and then we have $k_{2}-k_{1}=2\cdot3^{i}k-4j\geq3^{i}k\gg1$. Till now, the proof of the
lemma is finished.\QEDB

\medskip

As the argument in \cite{NSGT} (see Lemma $11$), we have the following lemma.

\begin{lem}\label{1.8.15} There must exist $r_{1}\in[2^{k_{1}},2^{k_{1}+1}]$, $r_{2}\in[2^{k_{2}-1},2^{k_{2}}]$ and a uniform constant $C$, such that
$$\int_{B(r_{1},r_{2})}(R(g_{\varepsilon}(t))-tr_{g_{\varepsilon}(t)}\theta_{\varepsilon})dV_{\varepsilon t}\leq CV<C\delta,
$$
for $t\geq1$, where $\delta>0$ and $V=Vol_{g_{\varepsilon}(t)}(B_{\varepsilon t}(k_{1},k_{2}))$ are obtained in Lemma \ref{1.8.14}.
\end{lem}

Finally, we prove that $diam(M,g_{\varepsilon}(t))$ can be uniformly bounded along the flows $(\ref{GKRF1})$ when $t\geq1$ by Perelman's argument. There exists a few differences between this proof and the original one, so we present a proof of Proposition \ref{1.8.16} in the appendix to reader's convenience.

\begin{pro}\label{1.8.16} $diam(M,g_{\varepsilon}(t))$ is uniformly bounded along the flow $(\ref{GKRF1})$ for ant $t\geq1$ and $\varepsilon$.
\end{pro}

{\bf Proof of Theorem \ref{1.8.2}:}\ \ By Lemma \ref{1.8.8} and Proposition \ref{1.8.16}, we obtain that $R(g_{\varepsilon}(t))-tr_{g_{\varepsilon}(t)}\theta_{\varepsilon}$ and $u_{\varepsilon}$ is uniformly bounded from above, while $|\nabla u_{\varepsilon}|_{g_\varepsilon(t)}$ is uniformly bounded.
Combining Proposition \ref{1.8.1} with Proposition \ref{1.8.5}, we prove the theorem.\QEDB

\section{The $C^{0}$ estimates for metric potential $\varphi_{\varepsilon}(t)$}
\setcounter{equation}{0}

In this section, for any $\beta\in(0,1)$, we obtain a uniform Sobolev inequality along the twisted K\"ahler-Ricci flows $(\ref{GKRF1})$. When it is in some finite interval, we have proved it in Proposition \ref{1.8.12}. When $t\geq 1$, from \cite{LW}(see also \cite{RGY} or \cite{QSZ}), we know that the Sobolev constants along the twisted K\"ahler-Ricci flows depend only on $n$, $\max(R(g_{\varepsilon}(1))-tr_{g_{\varepsilon}(1)}\theta_{\varepsilon})^{-}$ and $\mathcal{C}_{S}(M,g_{\varepsilon}(1))$, while the latter two can be uniformly bounded by Theorem \ref{1.8.2} and Proposition \ref{1.8.12}. So we have the following uniform Sobolev inequality (when $t\geq 1$) by  Q.S. Zhang's argument (\cite{QSZ}). To readers' convenience, we will give its proof in the appendix.

\begin{thm}\label{1.8.5.1} Let  $M$ be a compact K\"ahler manifold with complex dimension $n\geq2$ and $g_{\varepsilon}(t)$ be a solution of the twisted K\"ahler-Ricci flows $(\ref{GKRF1})$. Then there exist uniform constant A and B, such that for all $v\in W^{1,2}(M,g_{\varepsilon}(t))$, $\varepsilon>0$ and $t\geq 1$, we have
\begin{eqnarray}\label{3.22.31}(\int_{M}v^{\frac{2n}{n-1}}dV_{\varepsilon t})^{\frac{n-1}{n}}&\leq& B\int_{M}v^{2}dV_{\varepsilon t}+A\int_{M}\mid\nabla v\mid_{g_{\varepsilon}(t)}^{2}dV_{\varepsilon t}\\ \nonumber
&\ &+\frac{A}{4}\int_{M}R(g_{\varepsilon}(t))-tr_{g_{\varepsilon}(t)}\theta_{\varepsilon})v^{2}dV_{\varepsilon t}.\end{eqnarray}
Then by the uniform Perelman's estimates along the flow $(\ref{GKRF1})$ when $t\geq1$ and Proposition \ref{1.8.12}, we have
\begin{eqnarray}\label{3.22.31'}(\int_{M}v^{\frac{2n}{n-1}}dV_{\varepsilon t})^{\frac{n-1}{n}}\leq C\int_{M}\mid\nabla v\mid_{g_{\varepsilon}(t)}^{2}dV_{\varepsilon t}+C\int_{M}v^{2}dV_{\varepsilon t},\end{eqnarray}
for $t\geq0$, where C is a uniform constant.
\end{thm}

Next, we argue the uniform $C^{0}$ estimate for metric potential $\varphi_{\varepsilon}(t)$. We will denote $\phi_{\varepsilon}(t)=\varphi_{\varepsilon}(t)+k\chi(\varepsilon^{2}+|s|^{2}_{h})$ and discuss the $C^{0}$ estimates for $\phi_{\varepsilon}(t)$. First, we recall Aubin's functionals, Ding's functional and the twisted Mabuchi $\mathcal{K}$-energy functional.
\begin{eqnarray}\label{3.22.35}
I_{\omega_0}(\phi)&=&\frac{n!}{V}\int_M\phi(dV_{0}-dV_{\phi}), \\ \nonumber
J_{\omega_0}(\phi)&=&\frac{n!}{V}\int_0^1\int_M\dot{\phi}_t(dV_{0}-dV_{\phi_{t}})dt\\
&=&\frac{1}{V}\sum_{i=0}^{n-1}\frac{i+1}{n+1}\int_M\partial\phi\wedge\bar{\partial}\phi\wedge \omega_0^i\wedge w_\phi^{n-i-1},
\end{eqnarray}
where $\phi_t$ is a path with $\phi_0=c$, $\phi_1=\phi$.
\begin{eqnarray}\label{3.22.36}
F^0_{\omega_0}(\phi)&=&J_{\omega_0}(\phi)-\frac{n!}{V}\int_M\phi dV_0,\\
F_{\omega_0}(\phi)&=&J_{\omega_0}(\phi)-\frac{n!}{V}\int_M\phi dV_0-\log(\frac{n!}{V}\int_Me^{-u_{\omega_{0}}-\phi}dV_0),\\
\label{3.22.36'}\mathcal{M}_{ \omega_{0},\ \theta}(\phi)&=&-\beta(I_{\omega_{0}}(\phi)-J_{\omega_{0}}(\phi))-\frac{n!}{V}\int_{M}u_{\omega_{0}}(dV_{0}-dV_{\phi})\\ \nonumber
&\ &+\frac{n!}{V}\int_{M}\log\frac{\omega_{\phi}^{n}}{\omega_{0}^{n}}dV_{\phi},
\end{eqnarray}
where $u_{\omega_{0}}$ is the twisted Ricci potential of $\omega_0$, $i.e.$ $-Ric(\omega_{0})+\beta\omega_{0}+\theta=\sqrt{-1}\partial\bar{\partial} u_{\omega_{0}}$ and $\frac{1}{V}\int_{M}e^{-u_{\omega_{0}}}dV_{\omega_{0}}=1$.  The time derivatives of $I_{\omega_0}$, $J_{\omega_0}$ and $\mathcal{M}_{ \omega_{0},\ \theta}$ along any path $\phi_{t}$ can be written as follows:
\begin{eqnarray*}\frac{\partial}{\partial t}I_{\omega_0}(\phi_t)&=&\frac{n!}{V}\int_M\dot{\phi_t}(dV_{0}-dV_{\phi_{t}})-\frac{n!}{V}\int_M\phi_t\Delta\dot{\phi_t}dV_{\phi_{t}},\\
\frac{\partial}{\partial t}J_{\omega_0}(\phi_t)&=&\frac{n!}{V}\int_M \dot{\phi}_t(dV_{0}-dV_{\phi_{t}}),
\end{eqnarray*}
\begin{eqnarray*}
\frac{\partial}{\partial t}\mathcal{M}_{ \omega_{0},\ \theta}(\phi_t)&=&-\frac{n!}{V}\int_{M}\dot{\phi}_t(R(\omega_{\phi_t})-\beta n-tr_{\omega_{\phi_t}}\theta)dV_{\phi_t}.
\end{eqnarray*}

\begin{pro}\label{1.8.5.2} The integral $\int_{0}^{+\infty}e^{-\beta t}\|\nabla u_{\varepsilon}(t)\|^{2}_{L^{2}}dt$ is uniformly bounded.
\end{pro}

{\bf  Proof:}\ \ When $t\geq 1$, by Theorem \ref{1.8.2}, we know that $e^{-\beta t}\|\nabla u_{\varepsilon}(t)\|^{2}_{L^{2}}\leq C e^{-\beta t}$ for some uniform constant $C$, so $\int_{1}^{+\infty}e^{-\beta t}\|\nabla u_{\varepsilon}(t)\|^{2}_{L^{2}}dt$ is uniformly bounded. Then we only need to prove that $\int_{0}^{1}\|\nabla u_{\varepsilon}(t)\|^{2}_{L^{2}}dt$ is uniformly bounded. Through computing, we have
\begin{eqnarray}\label{3.22.37}\frac{d}{dt}(\mathcal{M}_{\omega_{0},\ \theta_{\varepsilon} }(\phi_{\varepsilon}(t))-\beta F^0_{\omega_0}(\phi_{\varepsilon}(t))-\frac{n!}{V}\int_M\dot{\phi}_{\varepsilon}(t) dV_{\varepsilon t})=0.
\end{eqnarray}
Hence
\begin{eqnarray}\label{1.5.8}&\ &\ \nonumber \mathcal{M}_{\omega_{0},\ \theta_{\varepsilon} }(\phi_{\varepsilon}(t))-\mathcal{M}_{\omega_{0},\ \theta_{\varepsilon} }(\phi_{\varepsilon}(0))\\ \nonumber
&=&\beta\int_{0}^{t}\int_M\frac{d}{ds}F_{\omega_0}^0(\phi_\varepsilon(s))ds+\frac{n!}{V}\int_M\dot{\phi}_\varepsilon(t)dV_{\varepsilon t}-\frac{n!}{V}\int_M\dot{\phi}_\varepsilon(0)dV_{\varepsilon 0}\\\nonumber
&=&-\beta\frac{n!}{V}\int_{0}^{t}\int_M\dot{\phi}_\varepsilon(s)dV_{\varepsilon s}ds+\frac{n!}{V}\int_M\dot{\phi}_\varepsilon(t)dV_{\varepsilon t}-\frac{n!}{V}\int_M\dot{\phi}_\varepsilon(0)dV_{\varepsilon 0}\\\nonumber
&=&-\frac{n!}{V}\beta\int_{0}^{t}\int_M(\dot{\phi}_\varepsilon(s)-e^{\beta s}\beta\varphi_\varepsilon(0))dV_{\varepsilon s}ds+\frac{n!}{V}\int_M(\dot{\phi}_\varepsilon(t)-e^{\beta t}\beta\varphi_\varepsilon(0))dV_{\varepsilon t}\\
&\ &-\frac{n!}{V}\int_M(\dot{\phi}_\varepsilon(0)-\beta\varphi_\varepsilon(0))dV_{\varepsilon}.
\end{eqnarray}
The evolution equation of $e^{-\beta t}(\dot{\phi}_\varepsilon(t)-e^{\beta t}\beta\varphi_\varepsilon(0))$ satisfies
\begin{eqnarray}(\frac{d}{dt}-\triangle_{g_\varepsilon(t)})(e^{-\beta t}(\dot{\phi}_\varepsilon(t)-e^{\beta t}\beta\varphi_\varepsilon(0)))=0.\end{eqnarray}
By maximum principle, we conclude that
\begin{eqnarray*}\sup_M|e^{-\beta t}(\dot{\phi}_\varepsilon(t)-e^{\beta t}\beta\varphi_\varepsilon(0))|\leq\sup_M|\log\frac{\omega_{\varepsilon}^{n}(\varepsilon^{2}+|s|_{h}^{2})^{1-\beta}}{\omega_{0}^{n}}+k\beta\chi(\varepsilon^2+
|s|_{h}^2)+F_{0}|
\end{eqnarray*}
for any $t\geq 0$. Hence there exists a uniform constant $C$ such that
\begin{eqnarray}\label{1.5.9}\sup\limits_{[0,T]\times M}|(\dot{\phi}_\varepsilon(t)-e^{\beta t}\beta\varphi_\varepsilon(0))|\leq Ce^{\beta T}.
\end{eqnarray}
On the other hand,
\begin{eqnarray}\frac{d}{dt}\mathcal{M}_{\omega_{0},\ \theta_{\varepsilon}}(\phi_{\varepsilon}(t))=
-\frac{n!}{V}\int_{M}|\nabla u_{\varepsilon}(t)|^{2}_{g_{\varepsilon}(t)}dV_{\varepsilon t}.
\end{eqnarray}
Integrating from $0$ to $1$ on both sides, we obtian
\begin{eqnarray}\int_{0}^{1}\|\nabla u_{\varepsilon}(t)\|^{2}_{L^{2}}dt=\mathcal{M}_{\omega_{0},\ \theta_{\varepsilon} }(\phi_{\varepsilon}(0))-\mathcal{M}_{\omega_{0},\ \theta_{\varepsilon} }(\phi_{\varepsilon}(1)).
\end{eqnarray}
Then the uniform bound of $\int_{0}^{1}\|\nabla u_{\varepsilon}(t)\|^{2}_{L^{2}}dt$ follows from (\ref{1.5.8}) and (\ref{1.5.9}).
\QEDB

\medskip

Now $\phi_{\varepsilon}(t)$ evolves along the following equation:
\begin{equation}\label{CMAF6}
\begin{cases}
  \frac{\partial \phi_{\varepsilon}(t)}{\partial t}=\log\frac{\omega_{\phi_{\varepsilon}(t)}^{n}}{\omega_{\varepsilon}^{n}}+F_{\varepsilon}+\beta\phi_{\varepsilon}(t)\\
  \phi_{\varepsilon}(t)|_{t=0}=c_{\varepsilon0}+k\chi(\varepsilon^{2}+|s|^{2}_{h})\\
  \end{cases}
\end{equation}
where $c_{\varepsilon0}=\frac{1}{\beta}(\int_{0}^{+\infty}e^{-\beta t}\|\nabla u_{\varepsilon}(t)\|^{2}_{L^{2}}dt-\frac{1}{V}\int_{M}F_{\varepsilon}dV_{\varepsilon}-\frac{1}{V}\int_{M}k\beta\chi(\varepsilon^{2}+|s|^{2}_{h})dV_{\varepsilon})$
and $F_{\varepsilon}=\log(\frac{\omega_{\varepsilon}^{n}}{\omega_{0}^{n}}(\varepsilon^{2}+|s|_{h}^{2})^{1-\beta})+F_{0}.$

 \begin{pro}\label{1.8.5.3}There exists a uniform constant $C$ such that
$$\|\dot{\phi}_{\varepsilon}(t)\|_{C^{0}}\leq C$$
for any $\varepsilon$ and $t$.
\end{pro}

{\bf  Proof:}\ \ As in \cite{PSS}, we let
\begin{eqnarray}\label{3.22.32}\alpha_{\varepsilon}(t)=\frac{1}{V}\int_{M}\dot{\phi}_{\varepsilon}(t)dV_{\phi_{\varepsilon}}
=\frac{1}{V}\int_{M}u_{\varepsilon}(t)dV_{\phi_{\varepsilon}}-c_{\varepsilon}(t).
\end{eqnarray}
Through computing, we have
\begin{eqnarray}\nonumber\frac{d}{dt}\alpha_{\varepsilon}(t)&=&\beta\alpha_{\varepsilon}(t)-\|\nabla\dot{\phi}_{\varepsilon}\|^{2}_{L^{2}},\\
\label{3.22.33}\nonumber e^{-\beta t}\alpha_{\varepsilon}(t)&=&\alpha_{\varepsilon}(0)-\int_{0}^{t}e^{-\beta s}\|\nabla\dot{\phi}_{\varepsilon}\|^{2}_{L^{2}}ds\\
&=&\frac{1}{V}\int_{M}u_{\varepsilon}dV_{\varepsilon}-c_{\varepsilon}(0)-\int_{0}^{t}e^{-\beta s}\|\nabla\dot{\phi}_{\varepsilon}\|^{2}_{L^{2}}ds.
\end{eqnarray}
Putting $u_{\varepsilon}=F_{\varepsilon}+k\beta\chi$ and $-c_{\varepsilon}(0)=\beta\varphi_{\varepsilon}(0)$ into $(\ref{3.22.33})$, we have
\begin{eqnarray*}e^{-\beta t}\alpha_{\varepsilon}(t)&=&\frac{1}{V}\int_{M}F_{\varepsilon}dV_{\varepsilon}+\frac{1}{V}\int_{M}k\beta\chi dV_{\varepsilon}-c_{\varepsilon}(0)-\int_{0}^{t}e^{-\beta s}\|\nabla\dot{\phi}_{\varepsilon}\|^{2}_{L^{2}}ds\\
&=&\frac{1}{V}\int_{M}F_{\varepsilon}dV_{\varepsilon}+\frac{1}{V}\int_{M}k\beta\chi dV_{\varepsilon}+\beta\varphi_{\varepsilon}(0)-\int_{0}^{t}e^{-\beta s}\|\nabla\dot{\phi}_{\varepsilon}\|^{2}_{L^{2}}ds\\
&=&\frac{1}{V}\int_{M}F_{\varepsilon}dV_{\varepsilon}+\frac{1}{V}\int_{M}k\beta\chi dV_{\varepsilon}-\int_{0}^{t}e^{-\beta s}\|\nabla\dot{\phi}_{\varepsilon}\|^{2}_{L^{2}}ds\\
&\ &+\int_{0}^{+\infty}e^{-\beta t}\|\nabla \dot{\varphi}_{\varepsilon}\|^{2}_{L^{2}}dt-\frac{1}{V}\int_{M}F_{\varepsilon}dV_{\varepsilon}-\frac{1}{V}\int_{M}k\beta\chi dV_{\varepsilon}\\
&=&\int_{t}^{+\infty}e^{-\beta s}\|\nabla\dot{\phi}_{\varepsilon}\|^{2}_{L^{2}}ds.
\end{eqnarray*}
When $t\geq 1$, by Theorem \ref{1.8.2}, we conclude that
\begin{eqnarray}\label{3.22.34}0\leq\alpha_{\varepsilon}(t)&=&\int_{t}^{+\infty}e^{\beta(t-s)}\|\nabla\dot{\phi}_{\varepsilon}\|^{2}_{L^{2}}ds\leq C.
\end{eqnarray}
Then $\dot{\phi}_{\varepsilon}(t)$ is bounded uniformly when $t\geq1$. Since $\phi_\varepsilon(0)$ is uniformly bounded, by (\ref{1.7.12}), it is easy to see that $$\|\dot{\phi}_\varepsilon(t)\|_{C^{0}([0,1]\times M)}\leq C$$
for some uniform constant $C$.
\QEDB

\medskip

Now we establish the relationship among the above functionals along the flow (\ref{CMAF6}).

\begin{pro}\label{1.8.5.4} There exists a uniform constant $C$, such that $\phi_{\varepsilon}(t)$ which evolves along the flow $(\ref{CMAF6})$ satisfies:
\begin{eqnarray*}&(i)& \ \ \ \mathcal{M}_{ \omega_{0},\ \theta_{\varepsilon} }(\phi_{\varepsilon}(t))-\beta F^0_{\omega_0}(\phi_{\varepsilon}(t))-\frac{n!}{V}\int_M\dot{\phi}_{\varepsilon}(t) dV_{\varepsilon t}=C_{\varepsilon},\\
&(ii)& \ \ |\beta F_{\omega_0}(\phi_{\varepsilon}(t))-\mathcal{M}_{\omega_{0},\ \theta_{\varepsilon} }(\phi_{\varepsilon}(t))|+|\beta F^0_{\omega_0}(\phi_{\varepsilon}(t))-\mathcal{M}_{\omega_{0},\ \theta_{\varepsilon} }(\phi_{\varepsilon}(t))|\leq C,\\
&(iii)& \ \ \frac{(n-1)!}{V}\int_M(-\phi_{\varepsilon}(t)) dV_{\varepsilon t}-C\leq J_{\omega_0}(\phi_{\varepsilon}(t))\leq \frac{n!}{V}\int_M\phi_{\varepsilon}(t) dV_0+C,\\
&(iv)& \ \ \frac{n!}{V}\int_M\phi_{\varepsilon}(t) dV_0\leq \frac{n\cdot n!}{V}\int_M(-\phi_{\varepsilon}(t)) dV_{\varepsilon t}-(n+1)\mathcal{M}_{ \omega_{0},\ \theta_{\varepsilon} } (\phi_{\varepsilon}(t))+C,\end{eqnarray*}
where $C_{\varepsilon}$ in $(i)$ can be bounded by a uniform constant $C$.
\end{pro}

{\bf  Proof:}\ \ Following the argument in \cite{PS1}, we only need to prove the two facts:
\begin{eqnarray*}&(1)&\ the\ constant\ C_{\varepsilon}\ in\ (i)\ can\ be\ bounded\ by\ a\ uniform\ constant\ C;\\
&(2)&\ \mathcal{M}_{ \omega_{0},\ \theta_{\varepsilon}}(\phi_{\varepsilon}(0))\ is\ uniformly\ bounded.
\end{eqnarray*}
We note that $\dot{\phi}_\varepsilon$ is uniformly bounded. $(i)$, $(ii)$, $(iii)$ and $(iv)$ can be easily deduced from the above two facts. Since
\begin{eqnarray}\label{3.22.37}\frac{d}{dt}(\mathcal{M}_{\omega_{0},\ \theta_{\varepsilon} }(\phi_{\varepsilon}(t))-\beta F^0_{\omega_0}(\phi_{\varepsilon}(t))-\frac{n!}{V}\int_M\dot{\phi}_{\varepsilon}(t) dV_{\varepsilon t})=0,
\end{eqnarray}
we obtain
\begin{eqnarray*}&\ &\mathcal{M}_{\omega_{0},\ \theta_{\varepsilon} }(\phi_{\varepsilon}(t))-\beta F^0_{\omega_0}(\phi_{\varepsilon}(t))-\frac{n!}{V}\int_M\dot{\phi}_{\varepsilon}(t) dV_{\varepsilon t}\\
&=&\mathcal{M}_{\omega_{0},\ \theta_{\varepsilon} }(\phi_{\varepsilon}(0))-\beta F^0_{\omega_0}(\phi_{\varepsilon}(0))-\frac{n!}{V}\int_M\dot{\phi}_{\varepsilon}(0) dV_{\varepsilon}\\
&=&\frac{n!}{V}\int_{M}\log\frac{\omega_{\varepsilon}^{n}(|s|^{2}_{h}+\varepsilon^{2})^{1-\beta}}{e^{-F_{0}}\omega_{0}^{n}}dV_{\varepsilon}+\frac{\beta n!}{V}\int_{M}\phi_{\varepsilon}(0)dV_{\varepsilon}\\
&\ &-\frac{ n!}{V}\int_{M}F_{0}+\log(|s|^{2}_{h}+\varepsilon^{2})^{1-\beta}dV_{0}-\frac{n!}{V}\int_M\dot{\phi}_{\varepsilon}(0) dV_{\varepsilon}.
\end{eqnarray*}
where the last equality can be bounded by a uniform constant.
Then we prove fact $(2)$. By the definition of $\mathcal{M}_{\omega_{0},\ \theta_{\varepsilon} }$, we have
\begin{eqnarray*}\mathcal{M}_{\omega_{0},\ \theta_{\varepsilon} }(\phi_{\varepsilon}(0))
&=&\frac{n!}{V}\int_{M}\log\frac{\omega_{\varepsilon}^{n}(|s|^{2}_{h}+\varepsilon^{2})^{1-\beta}}{e^{-F_{0}}\omega_{0}^{n}}dV_{\varepsilon}-\beta I_{\omega_{0}}(\phi_{\varepsilon}(0))+\beta J_{\omega_{0}}(\phi_{\varepsilon}(0))\\
&\ &-\frac{ n!}{V}\int_{M}F_{0}+\log(|s|^{2}_{h}+\varepsilon^{2})^{1-\beta}dV_{0}.
\end{eqnarray*}
Since $I_{\omega_{0}}(\phi_{\varepsilon}(0))$ is uniformly bounded and $\frac{1}{n}J_{\omega_0}\leq\frac{1}{n+1}I_{\omega_0}\leq J_{\omega_0}$, we prove the second fact.\QEDB

\medskip

Since we have proved that the uniform Sobolev inequality $(\ref{3.22.31'})$, Poincar\'e inequality $(\ref{3.22.37'})$ and $\|u_\varepsilon(t)\|_{C^0}$ can be uniformly bounded along the twisted K\"ahler-Ricci flows $(\ref{CMAF6})$, we obtain the following lemma by the argument in  \cite{PS1} (see Lemma $10$). The proof is completely similar, so we omit it.

\begin{pro}\label{1.8.5.5} We have the following estimate along the twisted K\"ahler-Ricci flow $(\ref{CMAF6})$
\begin{eqnarray}\label{3.22.40}osc(\phi_{\varepsilon}(t))\leq \frac{A}{V}\int_M\phi_{\varepsilon}(t) dV_0+B,
\end{eqnarray}
where the constants $A$ and $B$ are independent of $\varepsilon$ and $t$.
\end{pro}

We define the space of smooth  K\"ahler potentials as
\begin{eqnarray}\label{3.22.41}\mathcal{H}(\omega_0)=\{\phi\in C^{\infty}(M)|\ \omega_0+\sqrt{-1}\partial\bar{\partial}\phi>0\}.
\end{eqnarray}

\begin{thm}\label{1.8.5.6} Let $\phi_{\varepsilon}(t)$ be a solution of the flow $(\ref{CMAF6})$, and $\theta_\varepsilon=(1-\beta)(\omega_{0}+\sqrt{-1}\partial\overline{\partial}\log(\varepsilon^{2}+|s|_{h}^{2}))$ be a smooth closed semi-positive (1,1)-form, where $s$ is the defining section of divisor $D$ and $h$ is a smooth Hermitian metric on the line bundle associated to $D$. If the twisted Mabuchi $\mathcal{K}$-energy functional $\mathcal{M}_{ \omega_{0},\ \theta_{\varepsilon} }$ is uniformly proper on $\mathcal{H}(\omega_0)$, i.e. there exists a uniform function $f$ such that
\begin{eqnarray}\label{3.22.42}\mathcal{M}_{ \omega_{0},\ \theta_{\varepsilon}}(\phi)\geq f(J_{\omega_0}(\phi))
\end{eqnarray}
for any $\varepsilon$ and $\phi\in\mathcal{H}(\omega_0)$, where $f(t):\mathbb{R}^+\rightarrow\mathbb{R}$ is some monotone increasing function satisfying $\lim\limits_{t\rightarrow+\infty}f(t)=+\infty$, then there exists a uniform constant $C$ such that
\begin{eqnarray}\label{3.22.43}\|\phi_{\varepsilon}(t)\|_{C^{0}}\leq C.
\end{eqnarray}
\end{thm}

{\bf  Proof:}\ \ Since $\mathcal{M}_{\omega_{0},\ \theta_{\varepsilon} }(\phi_{\varepsilon}(t))$ decreases along the flow $(\ref{CMAF6})$ and $\mathcal{M}_{\omega_{0},\ \theta_{\varepsilon} }(\phi_{\varepsilon}(0))$ is uniformly bounded proved in Proposition \ref{1.8.5.4}. It follows that $J_{\omega_0}(\phi_{\varepsilon}(t))$ is uniformly bounded from above. Thus by Proposition \ref{1.8.5.4} $(iii)$, we have
\begin{eqnarray}\label{3.22.44}\int_{M}(-\phi_{\varepsilon}(t))dV_{\varepsilon t}\leq C.
\end{eqnarray}
Since $J_{\omega_0}\geq0$, applying $(\ref{3.22.42})$, we know that the twisted Mabuchi $\mathcal K$-energy $\mathcal{M}_{\omega_{0},\theta_{\varepsilon}  }(\phi_{\varepsilon}(t))$ is uniformly bounded from below. By Proposition \ref{1.8.5.4} $(iv)$, we have
\begin{eqnarray}\label{3.22.45}\int_{M}\phi_{\varepsilon}(t) dV_{0}\leq C,
\end{eqnarray}
where $C$ is a uniform constant. By this inequality and Green's formula with respect to the metric $g_{0}$, we get a uniform upper bound of $\sup\limits_{M}\phi_{\varepsilon}(t).$

By the normalization
 $$1=\frac{1}{V}\int_{M}dV_{\phi_{\varepsilon}(t)}=
\frac{1}{V}\int_{M}e^{\dot{\phi}_{\varepsilon}(t)-\beta\phi_{\varepsilon}(t)-F_{\varepsilon}}dV_{\varepsilon}$$
and the fact that $\|\dot{\phi}_{\varepsilon}(t)\|_{C^{0}}$ is uniformly bounded along the flow $(\ref{CMAF6})$, we have
$$0<C_{1}\leq\int_{M}e^{-\beta\phi_{\varepsilon}(t)}dV_{\varepsilon}\leq C_{2},$$
where $C_{1}$ and $C_{2}$ are uniform constants. This inequality easily implies a uniform lower bound for $\sup\limits_{M}\phi_{\varepsilon}(t)$. Combining with $(\ref{3.22.40})$ and $(\ref{3.22.45})$, we obtain a uniform bound for $\|\phi_{\varepsilon}(t)\|_{C^{0}}$. We also conclude that $$\|\varphi_{\varepsilon}(t)\|_{C^{0}}\leq C$$ for a uniform constant because $\chi$ is uniformly bounded.\QEDB

\section{The convergence of  the conical K\"ahlre-Ricci flow}
\setcounter{equation}{0}
In this section, we  consider the convergence of the conical K\"ahler-Ricci flow. The most important step in the convergence is to obtain a uniform $C^0$ estimate for $\phi_\varepsilon(t)$, so we  only need to get a uniform properness of the twisted Mabuchi $\mathcal{K}$-energy functional $\mathcal{M}_{ \omega_{0},\ \theta_{\varepsilon}}$ by Theorem \ref{1.8.5.6}. On the other hand, we notice that $\mathcal{M}_{ \omega_{0},\ \theta_{\varepsilon}}$ is associated with the Log Mabuchi $\mathcal{K}$-energy functional $\mathcal{M}_{\omega_{0},\ (1-\beta)D}$. So we first recall some contents of the Log Mabuchi $\mathcal{K}$-energy functional which are first introduced by C. Li and S. Sun in \cite{LS}.

For any $\phi\in\mathcal{H}(\omega_0)$,
\begin{eqnarray}\label{3.22.46}\nonumber\mathcal{M}_{\omega_{0},\ (1-\beta)D}(\phi)
&=&-\frac{n!}{V}\int_{M}H_{\omega_{0},(1-\beta)D}(dV_{0}-dV_{\phi})\\
&\ &+\frac{n!}{V}\int_{M}\log\frac{\omega_{\phi}^{n}}{\omega_{0}^{n}}dV_{\phi}-\beta(I_{\omega_{0}}(\phi)-J_{\omega_{0}}(\phi)),
\end{eqnarray}
where $H_{\omega_{0},(1-\beta)D}$ satisfies $-Ric(\omega_0)+\beta\omega_0+(1-\beta)\{D\}=\sqrt{-1}\partial\bar{\partial}H_{\omega_{0},(1-\beta)D}$ and  $\frac{1}{V}\int_M e^{-H_{\omega_{0},(1-\beta)D}}dV_0=1$. It is easy to see that up to a constant $H_{\omega_{0},(1-\beta)D}=F_0+(1-\beta)\log|s|_{h}^{2}$.

The Log Mabuchi $\mathcal{K}$-energy functional $\mathcal{M}_{\omega_{0},\ (1-\beta)D}:\mathcal{H}(\omega_0)\rightarrow\mathbb{R}$ is called proper if there is an inequality of the type
\begin{eqnarray}\label{3.22.47}\mathcal{M}_{\omega_{0},\ (1-\beta)D}(\phi)\geq f(J_{\omega_0}(\phi))
\end{eqnarray}
for any $\phi\in\mathcal{H}(\omega_0)$, where $f(t):\mathbb{R}^+\rightarrow\mathbb{R}$ is some monotone increasing function satisfying $\lim\limits_{t\rightarrow+\infty}f(t)=+\infty$. By using the linear property of Log Mabuchi $\mathcal{K}$-energy functional $\mathcal{M}_{\omega_{0},\ (1-\beta)D}$ \cite{LS} and the Donaldson's openness theorem \cite{SD2},  C. Li and S. Sun proved the following lemma.

\begin{lem}\label{1.8.6.1}(Corollary $1.4$ in \cite{LS}) If there is a conical K\"ahler-Einstein metric for $\beta\in(0,1)$, then the Log Mabuchi $\mathcal{K}$-energy functional $\mathcal{M}_{\omega_{0},\ (1-\beta)D}$ is proper.
\end{lem}

J. Song and X.W. Wang proved a similar result in \cite{SW}.  In both L-S and S-W's arguments,  the  Donaldson's openness theorem plays a key role. Recently, C.J. Yao provided an alternative proof of the Donaldson's openness theorem in \cite{CJY}. Here, we give a remark to Yao's paper.

Let's recall Yao's idea. Suppose that $\omega_{\varphi_\beta}$ is a weak conical K\"ahler-Einstein metric (see Definition $2.1$ in \cite{CJY}). Yao considered the following two parameter continuity path $\star^{\beta}_{\varepsilon, t}$ with $\varepsilon\in(0,1]$ and $t\in[0,\beta]$ to deform the K\"ahler metric $\omega_{\psi_{\varepsilon,\beta}}$:
$$ \star^{\beta}_{\varepsilon, t}:\ \ \ \ \left\{
\begin{aligned}
Ric( \omega_{\phi^\beta_{\varepsilon,t}})&=t \omega_{\phi^\beta_{\varepsilon,t}}+(\beta-t)\omega_{\varphi_\varepsilon}+(1-\beta)\chi_\varepsilon,\\
 \phi^\beta_{\varepsilon,0}&= \psi_{\varepsilon, \beta}\\
\end{aligned}
\right.$$
where $\{\omega_{\varphi_\varepsilon}\}$ is a sequence of smooth K\"ahler forms such that $\omega^n_{\varphi_\varepsilon}$ approximates $\omega_{\varphi_\beta}^n$ in $L^p$ ( $p\in(1,\frac{1}{1-\beta})$), $\psi_{\varepsilon,\beta}$ is a solution to $\star^{\beta}_{\varepsilon, 0}$ obtained by using S.T. Yau's result in \cite{STY} , and $\chi_\varepsilon=\omega_0+\sqrt{-1}\partial\bar{\partial}\log(\varepsilon^2+|s|_h^2)$.
By using B. Berndtsson's uniqueness theorem in \cite{BBERN}, Yao proved that there exists $\varepsilon_0>0$ such that the continuity path $\star^{\beta}_{\varepsilon, t}$ is solvable up to $t = \beta$ for all $\varepsilon\in(0,\varepsilon_0]$, and he also got a uniform bound of $\|\phi^\beta_{\varepsilon,\beta}\|_{L^\infty(M)}$ for $\varepsilon\in(0,\varepsilon_0]$ (Proposition $3.11$ in \cite{CJY}). Then he considered the new two parameter family continuity path
$$ \star_{\varepsilon, t}:\ \ \ \ \left\{
\begin{aligned}
Ric( \omega_{u_{\varepsilon,t}})&=t \omega_{u_{\varepsilon,t}}+(1-t)\chi_\varepsilon\\
 u_{\varepsilon,\beta}&= \phi^\beta_{\varepsilon, \beta}\\
\end{aligned}
\right.$$
and proved that there exist $\delta >0$ and $\tilde{\varepsilon}>0$ such that $u_{\varepsilon, t}$ is uniformly $L^{\infty }$ bounded for $(\varepsilon , t) \in (0, \tilde{\varepsilon}]\times (\beta -\delta, \beta +\delta)$ (Proposition 4.1 in \cite{CJY}). Using this uniform $C^0$ estimate of $u_{\varepsilon,t}$ for any $\beta ' \in (\beta -\delta, \beta +\delta)$ and $\varepsilon\in(0,\tilde{\varepsilon}]$, he proved that $\omega_{u_{\varepsilon,\beta'}}$ must converge to a weak conical K\"ahler-Einstein metric $\omega_{\varphi_{\beta'}}$ with angle $2\pi\beta'$ along $D$ as $\varepsilon\rightarrow0$. This gives another proof of the Donaldson's openness theorem.

\begin{rem}In C.J. Yao's argument of the uniform $C^0$ estimate for $u_{\varepsilon,t}$, we should find a fixed $\tilde{\delta }$ and prove that $\star_{\varepsilon,t}$ can be solved for any $t\in(\beta ,\beta+\tilde{\delta})$ and $\varepsilon\in(0,\varepsilon_0]$ beforehand. Since $\chi_\varepsilon$ is a strictly positive $(1,1)$-form, the linearized operator at $t=\beta$, which equals to $\triangle_{\phi_{\varepsilon_1,\beta}^\beta}+\beta$, is invertible for some standard Banach space. Yao used the standard implicit function theorem to perturb $t$ a little bit in both directions on $\star_{\varepsilon,t}$ for $\varepsilon\in(0,\varepsilon_0]$. But in general, the perturbation $(\beta -\delta(\varepsilon), \beta +\delta(\varepsilon))$ of $t$ depends on $\varepsilon$, i.e. we are not sure whether $\delta(\varepsilon)$ converges to $0$ as $\varepsilon\rightarrow0$. To deal with this problem, we can use G. Sz$\acute{e}$kelyhidi's result in \cite{SZ}, where he proved that for any $\omega_0\in C_1(M)$, if there exists a metric $\tilde{\omega}\in C_1(M)$ such that $Ric(\tilde{\omega})>k\tilde{\omega}$, then the equation
\begin{eqnarray}\label{1.29.1}Ric(\omega)=k \omega+(1-k)\omega_0
\end{eqnarray}
is solvable. Since $\star_{\varepsilon,t}$ can be solved at $t=\beta$ for $\varepsilon\in(0,\varepsilon_0]$ while $\chi_\varepsilon\in C_1(M)$ is a K\"ahler form, we have $Ric(\omega_{u_{\varepsilon_0,\beta}})>\beta \omega_{u_{\varepsilon_0,\beta}}$. Obviously, there exists a small number $\tilde{\delta}$ such that
\begin{eqnarray}\label{1.29.2}Ric(\omega_{u_{\varepsilon_0,\beta}})>(\beta+\tilde{\delta}) \omega_{u_{\varepsilon_0,\beta}}.
\end{eqnarray}
Replacing $\omega_0$ and $k$ in (\ref{1.29.1}) with $\chi_\varepsilon$ (any $\varepsilon\in(0,\varepsilon_0]$) and $t$ respectively, we know that $\star_{\varepsilon,t}$ can be solved for any $t\in[0,\beta+\tilde{\delta}]$ and $\varepsilon\in (0,\varepsilon_0]$ by G. Sz$\acute{e}$kelyhidi's result.
\end{rem}

Now, We will connect the uniform properness of the twisted Mabuchi $\mathcal{K}$-energy functional $\mathcal{M}_{\omega_{0},\ \theta_{\varepsilon}}$ with the properness of the Log Mabuchi $\mathcal{K}$-energy functional $\mathcal{M}_{\omega_{0},\ (1-\beta)\{D\}}$ by the following lemma.

\begin{lem}\label{1.8.6.2} If the Log Mabuchi $\mathcal{K}$-energy functional $\mathcal{M}_{\omega_{0},\ (1-\beta)\{D\}}$ is proper on $\mathcal{H}(\omega_0)$, then the twisted Mabuchi $\mathcal{K}$-energy functional $\mathcal{M}_{\omega_{0},\ \theta_{\varepsilon}}$ is uniformly proper on $\mathcal{H}(\omega_0)$.
\end{lem}

{\bf Proof:}\ \ By assumption, we have
$$\mathcal{M}_{\omega_{0},\ (1-\beta)\{D\}}(\phi)\geq C \tilde{f}(J_{\omega_{0}}(\phi))-C.$$
From the definition of Log Mabuchi $\mathcal{K}$-energy functional and twisted Mabuchi $\mathcal{K}$-energy functional, we have
\begin{eqnarray*}&\ &\mathcal{M}_{\omega_{0},\ \theta_{\varepsilon}}(\phi)-\mathcal{M}_{\omega_{0},\ (1-\beta)\{D\}}(\phi)\\
&=&\int_{M}(1-\beta)\log\frac{|s|_{h}^{2}}{\varepsilon^{2}+|s|_{h}^{2}}dV_{0}-
\int_{M}(1-\beta)\log\frac{|s|_{h}^{2}}{\varepsilon^{2}+|s|_{h}^{2}}dV_{\phi}\\
&\geq&\int_{M}(1-\beta)\log\frac{|s|_{h}^{2}}{\varepsilon^{2}+|s|_{h}^{2}}dV_{0}\\
&\geq& -C,
\end{eqnarray*}
where $C$ is independent of $\varepsilon$. Hence we obtain that
\begin{eqnarray*}\mathcal{M}_{\omega_{0},\ \theta_{\varepsilon}}(\phi)&\geq&\mathcal{M}_{\omega_{0},\ (1-\beta)\{D\}}(\phi)-C\\
&\geq&C \tilde{f}(J_{\omega_{0}}(\phi))-C.
\end{eqnarray*}
By setting $f=C\tilde{f}-C$, we get the uniform properness of $\mathcal{M}_{\omega_{0},\ \theta_{\varepsilon}}$ on $\mathcal{H}(\omega_0)$.\QEDB

\medskip

Next, we prove the convergence of the conical K\"ahler-Ricci flow.

\begin{thm}\label{1.8.6.3} Assume that there exists a conical K\"ahler-Einstein meric $\omega_{\beta, D}$, then the flow $(\ref{CKRF2})$ converges to the conical K\"ahler-Einstein meric $\omega_{\beta, D}$ in $C^{\infty}_{loc}$ topology outside $D$ and globally in the sense of currents.
\end{thm}

{\bf Proof:}\ \ First, by computing, we have
\begin{eqnarray}\label{3.22.49}\frac{d}{dt}\mathcal{M}_{\omega_{0},\ \theta_{\varepsilon}}(\phi_{\varepsilon})=
-\frac{n!}{V}\int_{M}|\partial\dot{\phi}_{\varepsilon}|^{2}_{g_\varepsilon(t)}dV_{\varepsilon t}.
\end{eqnarray}
Let $Y_{\varepsilon}(t)=\frac{n!}{V}\int_{M}|\partial\dot{\phi}_{\varepsilon}|^{2}_{g_\varepsilon(t)}dV_{\varepsilon t}$. By Lemma \ref{1.8.6.1} and Lemma \ref{1.8.6.2}, the twisted Mabuchi $\mathcal{K}$-energy functional $\mathcal{M}_{\omega_{0},\ \theta_{\varepsilon}}$ is uniformly proper, hence it can be bounded from below uniformly. For any $T$, we have
\begin{eqnarray}\label{3.22.50}\int_{1}^{T}Y_{\varepsilon}(t)dt\leq\int_{0}^{T}Y_{\varepsilon}(t)dt=\mathcal{M}_{\omega_{0},\ \theta_{\varepsilon}}(\phi_{\varepsilon}(0))-\mathcal{M}_{\omega_{0},\ \theta_{\varepsilon}}(\phi_{\varepsilon}(T))\leq C,\end{eqnarray}
where $C$ is a uniform constant. Define
\begin{eqnarray}\label{3.22.51}Y(t)=\frac{n!}{V}\int_{M}|\partial\dot{\phi}|^{2}_{g(t)}dV_{ t}.\end{eqnarray}
From Theorem \ref{1.8.2}, we know that $|\partial\dot{\phi}_{\varepsilon}|^{2}_{g_\varepsilon(t)}\leq C$ for a uniform constant $C$ when $t\geq 1$, so \begin{eqnarray}\label{3.22.52}\int_{1}^{T}Y_{\varepsilon_{i}}(t)dt\xrightarrow{\varepsilon_{i}\rightarrow0} \int_{1}^{T}Y(t)dt,\end{eqnarray}
where $\{\varepsilon_{i}\}$ is obtained in Theorem $3.1$. Hence we obtain
$$\int_{1}^{T}Y(t)dt\leq C.$$
When we let $T\rightarrow+\infty$, we get
\begin{eqnarray}\label{3.22.53}\int_{1}^{+\infty}Y(t)dt<\infty.\end{eqnarray}
Hence there exists a time sequence $\{t_{m}\}$, where $t_{m}\in[m,m+1)$ such that $Y(t_{m})\rightarrow0$ as $m\rightarrow+\infty$.

Next, $Y_{\varepsilon}(t)$ satisfies the following differential identity,
\begin{eqnarray*}\dot{Y}_{\varepsilon}(t)&=&\beta(n+1)Y_{\varepsilon}(t)-\int_{M}|\nabla \dot{\phi}_{\varepsilon}|_{g_\varepsilon(t)}^{2}(R(g_{\varepsilon}(t))-tr_{g_{\varepsilon}(t)}\theta_{\varepsilon})dV_{\varepsilon t}-\int_{M}|\overline{\nabla}\nabla \dot{\phi}_{\varepsilon}|_{g_\varepsilon(t)}^{2}dV_{\varepsilon t}\\
&\ &-\int_{M}|\nabla\nabla \dot{\phi}_{\varepsilon}|_{g_\varepsilon(t)}^{2}dV_{\varepsilon t}-\frac{1}{2}\int_{M}\theta(\nabla \dot{\phi}_{\varepsilon},\mathcal{J}\nabla \dot{\phi}_{\varepsilon})dV_{\varepsilon t}.\end{eqnarray*}
By Theorem \ref{1.8.2}, we have $|R(g_{\varepsilon}(t))-tr_{g_{\varepsilon}(t)}\theta_{\varepsilon}|\leq C$
for a uniform constant when $t\geqslant 1$. Hence we conclude that
\begin{eqnarray}\label{3.22.54}\dot{Y}_{\varepsilon}(t)\leq CY_{\varepsilon}(t).\end{eqnarray}
So $Y_{\varepsilon_{i}}(t)\leq e^{C(t-s)}Y_{\varepsilon_{i}}(s)$ for any $t>s$. Let $\varepsilon_{i}\rightarrow0$, we have
\begin{eqnarray}\label{3.22.55}Y(t)\leq e^{C(t-s)} Y(s)\end{eqnarray}
when $s,t\geqslant 1$. In particular,
$$Y(t)\leq e^{2C} Y(t_{m})$$
for all $t\in[m+1,m+2)$, and hence $Y(t)\rightarrow0$ as $t\rightarrow+\infty$.

Since the twisted Mabuchi $\mathcal{K}$-energy functional $\mathcal{M}_{\omega_{0},\ \theta_{\varepsilon}}$ is uniformly proper, we conclude that $\|\varphi_\varepsilon\|_{C^0}$ is uniformly bounded. From Proposition \ref{1.8.5.3}, we obtain that $\|\dot{\varphi}_\varepsilon\|_{C^0}$ is also uniformly bounded. By Theorem \ref{3.1}, we have
\begin{eqnarray}\label{3.22.56}\|\varphi\|_{C^0}\leq C,\ \ \ \|\dot{\varphi}\|_{C^0}\leq C, \ \ \ C^{-1}\omega^*\leqslant \omega_\varphi\leqslant C\omega^*\end{eqnarray}
for some uniform constant $C$ on $M\setminus D\times[0,+\infty)$. Then for any $K\subset\subset M\setminus D$, by Proposition \ref{2.2}, there exists a time sequence $\{t_{i}\}$ such that $\varphi(t_{i})$ converges in $C^{\infty}$ topology to a smooth function $\varphi_{\infty}$ on $K$.
\begin{eqnarray*}&\ &\int_{K}|\partial(log\frac{\omega_{\varphi(t_{i})}^{n}}{\omega_{0}^{n}}+F_{0}+\beta(k|s|_{h}^{2}+\varphi(t_{i}))+\log|s|_{h}^{2(1-\beta)})|^{2}_{g_{0}}dV_{0}\\
&\leq&C\int_{K}|\partial(log\frac{\omega_{\varphi(t_{i})}^{n}}{\omega_{0}^{n}}+F_{0}+\beta(k|s|_{h}^{2}+\varphi(t_{i}))+\log|s|_{h}^{2(1-\beta)})|^{2}_{g_(t_i)}dV_{t_i}\\
&\leq&C\int_{M}|\partial(log\frac{\omega_{\varphi(t_{i})}^{n}}{\omega_{0}^{n}}+F_{0}+\beta(k|s|_{h}^{2}+\varphi(t_{i}))+\log|s|_{h}^{2(1-\beta)})|^{2}_{g_(t_i)}dV_{t_i}\\
&=&C\int_{M}|\partial\dot{\phi}(t_{i})|^{2}_{g_(t_i)}dV_{t_i}\rightarrow0.\\
\end{eqnarray*}
On the other hand, we have
\begin{eqnarray*}&\ &\int_{K}|\partial(log\frac{\omega_{\varphi(t_{i})}^{n}}{\omega_{0}^{n}}+F_{0}+\beta(k|s|_{h}^{2}+\varphi(t_{i}))+\log|s|_{h}^{2(1-\beta)})|^{2}_{g_{0}}dV_{0}\\
&\rightarrow&\int_{K}|\partial(log\frac{\omega_{\varphi_{\infty}}^{n}}{\omega_{0}^{n}}+F_{0}+\beta(k|s|_{h}^{2}+\varphi_{\infty})+\log|s|_{h}^{2(1-\beta)})|^{2}_{g_{0}}dV_{0}.\\
\end{eqnarray*}
By the uniqueness of the limit,
$$\int_{K}|\partial(\log\frac{\omega_{\varphi_{\infty}}^{n}}{\omega_{0}^{n}}+F_{0}+\beta(k|s|_{h}^{2}+\varphi_{\infty})+\log|s|_{h}^{2(1-\beta)})|^{2}_{g_{0}}dV_{0}=0.$$
Hence
\begin{eqnarray}\label{3.22.57} Ric(\omega_{\varphi_{\infty}})=\beta\omega_{\varphi_{\infty}},\ \ \ \ \ on\ \ K.\end{eqnarray}

At the same time, there exists a time subsequence denoted also by $\{t_{i}\}$ such that $\varphi(t_{i})$ converges in $C_{loc}^{\infty}$ topology outside $D$ to a function $\varphi_{\infty}$ which is smooth on $M\setminus D$. We also have $\dot{\varphi}(t_{i})$ converge to some constant $C$ in $C^{\infty}_{loc}$ topology outside D. For any $(n-1,n-1)$-form $\eta$, since $\log\frac{\omega_{\varphi}^{n}|s|_{h}^{2(1-\beta)}}{\omega_{0}^{n}}$ and $\|\varphi\|_{C^0}$ are uniformly bounded, in the sense of currents,  we have
\begin{eqnarray*}&\ &\int_M\sqrt{-1}\partial\bar{\partial}\frac{\partial\varphi(t_{i})}{\partial t}\wedge\eta\\
&=&\int_M\sqrt{-1}\partial\bar{\partial}(\log\frac{\omega_{\varphi(t_{i})}^{n}}{\omega_{0}^{n}}+F_{0}+\beta(k|s|_{h}^{2}+\varphi(t_{i}))
+\log|s|_{h}^{2(1-\beta)}\wedge\eta
\end{eqnarray*}
\begin{eqnarray*}
&=&\int_{M}(\log\frac{\omega_{\varphi(t_{i})}^{n}}{\omega_{0}^{n}}+F_{0}+\beta(k|s|_{h}^{2}+\varphi(t_{i}))
+\log|s|_{h}^{2(1-\beta)})\sqrt{-1}\partial\bar{\partial}\eta\\
&\xrightarrow{t_{i}\rightarrow +\infty}&\int_{M}(\log\frac{\omega_{\varphi_{\infty}}^{n}}{\omega_{0}^{n}}+F_{0}+\beta(k|s|_{h}^{2}+\varphi_{\infty})
+\log|s|_{h}^{2(1-\beta)})\sqrt{-1}\partial\bar{\partial}\eta\\
&=&\int_M\sqrt{-1}\partial\bar{\partial}(\log\frac{\omega_{\varphi_{\infty}}^{n}}{\omega_{0}^{n}}+F_{0}+\beta(k|s|_{h}^{2}+\varphi_{\infty})
+\log|s|_{h}^{2(1-\beta)}\wedge\eta\\
&=&\int_M(-Ric(\omega_{\varphi_{\infty}})+\beta\omega_{\varphi_{\infty}}+(1-\beta)[D])\wedge\eta.
\end{eqnarray*}
Let $K\subset\subset M\setminus D$ be a compact subset, $\int_{M\setminus K}\sqrt{-1}\partial\bar{\partial}\eta=\delta$, and $\delta\rightarrow0$ when $K\rightarrow M\setminus D$. Then
\begin{eqnarray*}&\ &|\int_{M}(\frac{\partial\varphi(t_{i})}{\partial t}-C)\sqrt{-1}\partial\bar{\partial}\eta|\\
&=&|\int_{K}(\frac{\partial\varphi(t_{i})}{\partial t}-C)\sqrt{-1}\partial\bar{\partial}\eta+\int_{M\setminus K}(\frac{\partial\varphi(t_{i})}{\partial t}-C)\sqrt{-1}\partial\bar{\partial}\eta|\\
&\leq&|\int_{K}(\frac{\partial\varphi(t_{i})}{\partial t}-C)\sqrt{-1}\partial\bar{\partial}\eta|+\tilde{C}\delta\\
&\xrightarrow{t_{i}\rightarrow +\infty}&\tilde{C}\delta.
\end{eqnarray*}
When letting $K\rightarrow M\setminus D$, we have
\begin{eqnarray*}
\int_M\sqrt{-1}\partial\bar{\partial}\frac{\partial\varphi(t_{i})}{\partial t}\wedge\eta=\int_{M}\frac{\partial\varphi(t_i)}{\partial t}\sqrt{-1}\partial\bar{\partial}\eta
&\xrightarrow{t_{i}\rightarrow +\infty}&0.
\end{eqnarray*}
Hence, we obtain that
\begin{eqnarray}\label{3.22.58}Ric(\omega_{\varphi_{\infty}})=\beta\omega_{\varphi_{\infty}}+(1-\beta)[D]\end{eqnarray}
in the current sense. Since $C^{-1}\omega\leq\omega_{\varphi}\leq C\omega$, we also have
\begin{eqnarray}\label{3.22.59}C^{-1}\omega\leq\omega_{\varphi_{\infty}}\leq C\omega.\end{eqnarray}

By estimates $(\ref{3.22.56})$, from the proof of Proposition $\ref{3.2}$, we know that $\|\varphi\|_{C^{\alpha}}$ is uniformly bounded for some $\alpha\in(0,1)$, so the limit $\varphi_{\infty}$ is also H$\ddot{o}$lder continuous on $M$. On the basis of the properness of the Log Mabuchi $\mathcal{K}$-energy functional, there must exist $\omega_\infty=\omega_{\beta,D}$ by the uniqueness of the conical K\"ahler-Einstein metric with bound potential proved by R. Berman in \cite{RB}.

At last, we use the uniqueness to prove that flow $(\ref{CKRF2})$ certainly converges to $\omega_{\varphi_{\infty}}$ in $C_{loc}^{\infty}$ topology outside D and in current sense as $t\rightarrow+\infty$. If not, there exist $K\subset\subset M\setminus D$, an integer $k>0$ , $\epsilon_0>0$, and a time subsequence $\{t_{i}^{'}\}$ such that
\begin{eqnarray}\label{3.22.60}\|\sqrt{-1}\partial\bar{\partial}(\varphi(t_{i}^{'})-\varphi_{\infty})\|_{C^{k}(K)}\geq\epsilon_0.\end{eqnarray}
Since $\varphi(t_{i}^{'})$ is $C_{loc}^{\infty}$ bounded, there exists a subsequence which we also denote it by $\{t_{i}^{'}\}$, such that $\varphi(t_{i}^{'})$ converges in $C_{loc}^{\infty}$ topology to a function $\hat{\varphi}_{\infty}$ and
\begin{eqnarray}\label{3.22.61}\|\sqrt{-1}\partial\bar{\partial}(\hat{\varphi}_{\infty}-\varphi_{\infty})\|_{C^{k}(K)}\geq\epsilon.\end{eqnarray}
By the same argument above, we know that $\omega_{\hat{\varphi}_{\infty}}$ is also a conical K\"ahler-Einstein metric with H$\ddot{o}$lder continuous potential $\hat{\varphi}_{\infty}$. But $\omega_{\hat{\varphi}_{\infty}}\not \equiv\omega_{\varphi_{\infty}}$ by $(\ref{3.22.61})$, which is impossible by R. Berman's uniqueness results. Hence we get the convergence of the conical K\"ahler-Ricci flow.\QEDB

\section*{Appendix}
\setcounter{equation}{0}
In the appendix, we first give the proof of the uniform Sobolev inequality along the twisted K\"ahler-Ricci flows $(\ref{GKRF1})$, i.e. we first prove Proposition \ref{1.8.12} and Theorem \ref{1.8.5.1}.

\medskip

{\bf Proof of Proposition \ref{1.8.12}:}\ \ Under the appropriate coordinate system (see Lemma \ref{1.8.11}), metric $\omega_{\varepsilon}$ can be written as follows.
\begin{align*}\label{3.22.2}\omega_\varepsilon&=\omega_0+k e^{-\varphi}(\varepsilon^2+|z^n|^2e^{-\varphi})^{\beta-1}\sqrt{-1}dz^n\wedge d\bar{z}^n\\
&\ -k e^{-\varphi}z^n\frac{\partial\varphi}{\partial z^\alpha}(\varepsilon^2+|z^n|^2e^{-\varphi})^{\beta-1}\sqrt{-1}dz^\alpha\wedge d\bar{z}^n\\
&\ -k e^{-\varphi}\bar{z}^n\frac{\partial\varphi}{\partial \bar{z}^\beta}(\varepsilon^2+|z^n|^2e^{-\varphi})^{\beta-1}\sqrt{-1}dz^n\wedge d\bar{z}^\beta \tag{A.1}\\
&\ +k e^{-\varphi}|z^n|^2\frac{\partial\varphi}{\partial z^\alpha}\frac{\partial\varphi}{\partial \bar{z}^\beta}(\varepsilon^2+|z^n|^2e^{-\varphi})^{\beta-1}\sqrt{-1}dz^\alpha\wedge d\bar{z}^\beta\\
&\ -\frac{k}{\beta}((\varepsilon^2+|z^n|^2e^{-\varphi})^{\beta}-\varepsilon^{2\beta})\frac{\partial^2\varphi}{\partial z^\alpha \partial \bar{z}^\beta}\sqrt{-1}dz^\alpha\wedge d\bar{z}^\beta.
\end{align*}

We consider the map
\begin{align*}\label{4.1.0}\Psi_{\varepsilon}:\ (z^{1},\ z^{2},\ \cdots,\ z^{n-1},\ \xi)\longmapsto(z^{1},\ z^{2},\ \cdots,\ z^{n-1},\ z^{n}),\tag{A.2}\end{align*}
where $z^{n}=(\varepsilon^{2\beta}+|\xi|^{2})^{\frac{1}{2\beta}-\frac{1}{2}}\xi$. Now, we want to show that $\Psi_{\varepsilon}^{\ast} (g_{\varepsilon}) $ is uniformly equivalent to the Euclidean metric in a small neighborhood of the divisor $D$.

By a direct calculation, we only need to deal with the following term
\begin{align*}\label{4.1.1}\Psi_{\varepsilon}^*(k e^{-\varphi}(\varepsilon^2+|z^n|^2e^{-\varphi})^{\beta-1}dz^n\cdot d\bar{z}^n).\tag{A.3}\end{align*}
We will show that $(\ref{4.1.1})$ is uniformly equivalent to the Euclidean metric on $\mathbb{C}$.

Now we estimate it by the polar coordinates transformation. Let $z^n=x+\sqrt{-1}y$, $x=r\cos\theta$ and $y=r\sin\theta$, we have
\begin{align*}\label{4.1.2} \begin{array}{lll}dz^n\cdot d\bar{z}^n&=&dz^n\otimes d\bar{z}^n+ d\bar{z}^n\otimes dz^{n}\\&=&2(dr^2+r^2d\theta^2).\tag{A.4}
\end{array}
\end{align*}
We let $\xi=u+\sqrt{-1}v$, $u=\rho\cos\theta_1$ and $v=\rho\sin\theta_1$. By the definition of $\Psi_\varepsilon$, we know that $\theta_1=\theta$ and $r=(\varepsilon^{2\beta}+\rho^{2})^{\frac{1}{2\beta}-\frac{1}{2}}\rho$. Hence we have
\begin{align*}\label{4.1.3}\begin{array}{lll} &\ &\Psi_{\varepsilon}^*(k e^{-\varphi}(\varepsilon^2+|z^n|^2e^{-\varphi})^{\beta-1}dz^n\cdot d\bar{z}^n)\\
 &=&2k e^{-\varphi\circ\Psi_{\varepsilon}}(\varepsilon^2+(\varepsilon^{2\beta}+\rho^{2})^{\frac{1}{\beta}-1}\rho^2 e^{-\varphi\circ\Psi_{\varepsilon}})^{\beta-1}(\varepsilon^{2\beta}+\rho^{2})^{\frac{1}{\beta}-1}\cdot\tag{A.5}\\
&\ &((1+(\frac{1}{\beta}-1)(\varepsilon^{2\beta}+\rho^{2})^{-1}\rho^2)^2d\rho^2+\rho^2d\theta_1^2).
\end{array}
\end{align*}
Because $1\leq(1+(\frac{1}{\beta}-1)(\varepsilon^{2\beta}+\rho^{2})^{-1}\rho^2)^2\leq\frac{1}{\beta^2}$, we only need to prove that the term \begin{align*}\label{2015.1.31}(\varepsilon^2+(\varepsilon^{2\beta}+\rho^{2})^{\frac{1}{\beta}-1}\rho^2 e^{-\varphi\circ\Psi_{\varepsilon}})^{\beta-1}(\varepsilon^{2\beta}+\rho^{2})^{\frac{1}{\beta}-1}\tag{A.6}\end{align*} can be uniformly bounded, and the uniform lower bound is away from $0$. Firstly, we bound it from below,
\begin{eqnarray*}&\ &(\varepsilon^2+(\varepsilon^{2\beta}+\rho^{2})^{\frac{1}{\beta}-1}\rho^2 e^{-\varphi\circ\Psi_{\varepsilon}})^{\beta-1}(\varepsilon^{2\beta}+\rho^{2})^{\frac{1}{\beta}-1}\\
&\geq&((\varepsilon^{2\beta}+\rho^{2})^{\frac{1}{\beta}}+(\varepsilon^{2\beta}+\rho^{2})^{\frac{1}{\beta}-1}(\varepsilon^{2\beta}+\rho^{2}) e^{-\varphi\circ\Psi_{\varepsilon}})^{\beta-1}(\varepsilon^{2\beta}+\rho^{2})^{\frac{1}{\beta}-1}\\
&\geq&(\varepsilon^{2\beta}+\rho^{2})^{\frac{1}{\beta}(\beta-1)}(1+e^{-\varphi\circ\Psi_{\varepsilon}})^{\beta-1}(\varepsilon^{2\beta}+\rho^{2})^{\frac{1}{\beta}-1}\\
&=&(1+e^{-\varphi\circ\Psi_{\varepsilon}})^{\beta-1}\geq c>0,
\end{eqnarray*}
where $c$ is independent of $\varepsilon$. Secondly, we prove that the term $(\ref{2015.1.31})$ can be bounded from above. Let $\varepsilon^\beta=l\cos\vartheta$ and $\rho=l\sin\vartheta$, where $\vartheta\in[0,\frac{\pi}{2}]$, then we have
\begin{eqnarray*}&\ &(\varepsilon^2+(\varepsilon^{2\beta}+\rho^{2})^{\frac{1}{\beta}-1}\rho^2 e^{-\varphi\circ\Psi_{\varepsilon}})^{\beta-1}(\varepsilon^{2\beta}+\rho^{2})^{\frac{1}{\beta}-1}\\
&=&(l^{\frac{2}{\beta}}\cos^{\frac{2}{\beta}}\vartheta+l^{2(\frac{1}{\beta}-1)}l^2\sin^2\vartheta e^{-\varphi\circ\Psi_{\varepsilon}})^{\beta-1}l^{2(\frac{1}{\beta}-1)}\\
&=&(\frac{1}{\cos^{\frac{2}{\beta}}\vartheta+\sin^2\vartheta e^{-\varphi\circ\Psi_{\varepsilon}}})^{1-\beta}\\
&\leq&(\frac{1}{\cos^{\frac{2}{\beta}}\vartheta +\sin^{\frac{2}{\beta}}\vartheta })^{1-\beta}e^{c(1-\beta)}\\
&\leq&2^{(\frac{1}{\beta}-1)(1-\beta)}e^{c(1-\beta)}.
\end{eqnarray*}
In conclusion, it shows that
\begin{eqnarray*}C_1(d\rho^2+\rho^2d\theta_1^2)\leq\Psi_{\varepsilon}^*(k e^{-\varphi}(\varepsilon^2+|z^n|^2e^{-\varphi})^{\beta-1}dz^n\cdot d\bar{z}^n)\leq C_2(d\rho^2+\rho^2d\theta_1^2)
\end{eqnarray*}
for some uniform constants $C_1$ and $C_2$ independent of $\varepsilon$. It is easy to see that the pull-back metric $\Psi_{\varepsilon}^{\ast}(g_{\varepsilon}$) is uniformly equivalent to the Euclidean metric in a small neighborhood of the divisor D. Therefore, the Sobolev inequality holds if the function $v$ is supported in the above coordinate charts. The global case follows in the standard way by using a partition of unity. Following this argument, we prove the following Sobolev inequality
\begin{align*}(\int_{M}v^{\frac{2n}{n-1}}dV_{\varepsilon})^{\frac{n-1}{n}}\leq C(\int_{M}|dv|^{2}_{g_{\varepsilon}}dV_{\varepsilon}+\int_{M}|v|^{2}dV_{\varepsilon}).\tag{A.7}\end{align*}
To prove (\ref{1.5.6}), we only need to prove that $\omega_\varepsilon$ and $\omega_{\varepsilon}(t)$ are uniformly equivalent when $t\in[0,2]$. Noting that the metric $\omega_\varepsilon(t)$ is independent of the choice of the initial constant $\varphi_\varepsilon(0)$, without loss of generality, we assume $\varphi_\varepsilon(0)=0$. By (\ref{1.7.11}), we have $\|\dot{\varphi}_\varepsilon(t)\|_{C^{0}([0,2]\times M)}\leq C$ and $\|\varphi_\varepsilon(t)\|_{C^{0}([0,2]\times M)}\leq C$ for some uniform constant $C$ only depending on $\log\frac{\omega_{\varepsilon}^{n}(\varepsilon^{2}+|s|_{h}^{2})^{1-\beta}}{\omega_{0}^{n}}+k\beta\chi(\varepsilon^2+|s|_{h}^2)+F_{0}$. Then the uniform equivalence between the metrics follows from Proposition \ref{2.1}.\QEDB

\medskip

{\bf  Proof of Theorem \ref{1.8.5.1}:}\ \ In the proof, we only to consider the Sobolev inequality along the twisted K\"ahler-Ricci flow $(\ref{GKRF1})$ for $t\geq 1$. This proof is almost the same as that in \cite{QSZ} with the only difference that we require the constants independent of $\varepsilon$ in addition, so we give the proof briefly here.

Step $1$. By using the monotonicity of the functional $\mu_{\theta_\varepsilon}(g_\varepsilon(s),\tau(s))$ (see Lemma \ref{1.8.10}) and taking $\tau(s)=1-e^{-\beta t}(1-\delta^2)e^{\beta s}$, where $\delta\in(0,1)$, we conclude that
\begin{align*}\label{3.28.1}\nonumber\int_{M}v^{2}\log v^{2}dV_{\varepsilon t}&\leq\int_{M}\delta^{2}((R(g_{\varepsilon}(t))-tr_{g_{\varepsilon}(t)}\theta_{\varepsilon})v^{2}+4|\nabla v|_{g_{\varepsilon}(t)}^{2})dV_{\varepsilon t}\\
&\ -2n\log\delta+L_{1}+\max\limits_{M}(R(g_{\varepsilon}(1))-tr_{g_{\varepsilon}(1)}\theta_{\varepsilon})^{-},\tag{A.8}\end{align*}
where $L_{1}=n\log (n\mathcal{C}_{S}(M,g_\varepsilon(1))^{2})-n\log2-n+\frac{4}{\mathcal{C}_{S}(M,g_\varepsilon(1))^{2}}V^{-\frac{1}{n}}$.

Step $2$. Fixing a time $t_0\geq 1$ during the twisted K\"ahler-Ricci flow, we show that the upper bound of short time heat kernel for the fundamental  solution of equation
\begin{equation*}\label{3.28.2}\triangle_{g_\varepsilon(t_0)} u(x,t)-\frac{1}{4}(R(g_\varepsilon(t_{0}))-tr_{g_\varepsilon(t_{0})}\theta_\varepsilon)u(x,t)-\frac{\partial}{\partial t}u(x,t)=0\tag{A.9}\end{equation*}
under the fixed metric $g_\varepsilon(t_{0})$.

Let $u$ be a positive solution of equation $(\ref{3.28.2})$. From the given $T\in(0,1]$ and $t\in(0,T]$, we take $p(t)=\frac{T}{T-t}$. Differentating $\|u\|_{p(t)}$ and putting (\ref{3.28.1}) into it, then after integrating from $t=0$ to $t=T$ on both sides, we have
\begin{equation*}\label{3.28.3}\log\frac{\|u(\cdot,T)\|_{\infty}}{\|u(\cdot,0)\|_{1}}\leq-n\log T+L+2\max\limits_{M}(R(g_\varepsilon(1))-tr_{g_\varepsilon(1)}\theta_\varepsilon)^{-},\tag{A.10}\end{equation*}
where $L=L_{1}+2n$.

Since $u(x,T)=\int_{M}P_\varepsilon(x,y,T)u(y,0)dV_{\varepsilon t_{0}}$, where $P_\varepsilon$ is the heat kernel of equation $(\ref{3.28.2})$,
\begin{equation*}\label{3.28.4}P_\varepsilon(x,y,T)\leq\frac{\exp(L+2\max\limits_{M}(R(g_\varepsilon(1))-tr_{g_\varepsilon(1)}
\theta_\varepsilon)^{-})}{T^{n}}:=\frac{\Lambda}{T^{n}}.\tag{A.11}\end{equation*}

Step $3$. Let $F_\varepsilon=\max\limits_{M}(R(g_\varepsilon(1))-tr_{g_\varepsilon(1)}\theta_\varepsilon)^{-}$, $\Psi_{\varepsilon,t_0}=\frac{1}{4}(R(g_\varepsilon(t_0))+tr_{g_\varepsilon(t_0)}\theta_\varepsilon)+F_\varepsilon+1\geqslant 1$ and $P_{F_\varepsilon}$ be the heat kernel of operator $\triangle_{g_{\varepsilon}(t_0)}-\Psi_{\varepsilon,t_0}$. For $t>0$ and $y\in M$,
\begin{align*}\label{2015.1.24.10}(\frac{d}{dt}-\triangle_{g_{\varepsilon}(t_0)})P_{F_\varepsilon}(x,y,t+1)=
-\Psi_{\varepsilon,t_0}P_{F_\varepsilon}(x,y,t+1).\tag{A.12}
\end{align*}

By maximum principle and (\ref{3.28.4}), $P_{F_\varepsilon}$ obeys the global upper bound
\begin{align*}\label{3.28.5}P_{F_\varepsilon}(x,y,t)\leq \tilde{C}t^{-n},\ \ \ \ t>0,\tag{A.13}
\end{align*}
where $\tilde{C}$ depends only on $\Lambda$ and $n$. Moreover, by H$\ddot{o}$lder inequality, for any $f\in L^2(M)$, we have
\begin{eqnarray*}|\int_{M}P_{F_\varepsilon}(x,y,t)f(y)dV_{\varepsilon t_{0}}|\leq(\int_{M}P_{F_\varepsilon}^{2}(x,y,t)dV_{\varepsilon t_{0}})^{\frac{1}{2}}\|f\|_{L^2}\leq \tilde{C}^{\frac{1}{2}}t^{-\frac{n}{2}}\|f\|_{L^2}.
\end{eqnarray*}
Then the Sobolev inequality follows Theorem $2.4.2$ in \cite{EBD} and the constants in inequality depend only on $\tilde{C}$, $2n$ and $\frac{1}{2n-2}$. By the expression of $\Lambda$, Lemma \ref{1.8.6} and Proposition \ref{1.8.12}, we know that the constants are independent of $\varepsilon$ and $t$. \QEDB

\medskip

At last, we prove Proposition \ref{1.8.16} by contradiction.

\medskip

{\bf Proof of Proposition \ref{1.8.16}:}\ \ If $diam(M,g_{\varepsilon}(t))$ is not uniformly bounded, there exist $\{t_{i}\}\subset[1,+\infty)$ and $\varepsilon_{i}\rightarrow0$ such that $diam(M,g_{\varepsilon_{i}}(t_{i}))\rightarrow+\infty$. Let $\delta_{i}\rightarrow0$ be a sequence consisting of positive numbers, which corresponds to $\{t_{i}\}$ and $\{\varepsilon_{i}\}$. By Lemma \ref{1.8.14}, we can find sequences $\{k_{1}^{i}\}$ and $\{k_{2}^{i}\}$, such that
\begin{align*}\label{3.22.28}Vol_{g_{\varepsilon_{i}}(t_{i})}(B_{\varepsilon_{i} t_{i}}(k_{1}^{i},k_{2}^{i}))&<\delta_{i},\tag{A.14}\\
Vol_{g_{\varepsilon_{i}}(t_{i})}(B_{\varepsilon_{i} t_{i}}(k_{1}^{i},k_{2}^{i}))&\leq 2^{20n}Vol_{g_{\varepsilon_{i}}(t_{i})}(B_{\varepsilon_{i} t_{i}}(k_{1}^{i}+2,k_{2}^{i}-2)).\tag{A.15}\end{align*}
Let $r_{1}^{i}\in[2^{k_{1}^{i}},2^{k_{1}^{i}+1}]$ and $r_{2}^{i}\in[2^{k_{2}^{i}-1},\ 2^{k_{2}^{i}}]$ given in Lemma \ref{1.8.15} for each $i$, $\phi_{i}$ be cut off functions such that $\phi_{i}=1$ on $[2^{k_{1}^{i}+2},2^{k_{2}^{i}-2}]$ and $\phi_{i}=0$ on $(-\infty,r_{1}^{i}]\bigcup[r_{2}^{i},+\infty).$ Define
\begin{align*}\label{3.22.29}u_{i}(x)=e^{C_{i}}\phi_{i}(dist_{\varepsilon_{i}t_{i}}(x,p_{i})),\tag{A.16}\end{align*}
where $u_{\varepsilon_{i}}(p_{i},t_{i})=\inf_{M}u_{\varepsilon_{i}}(y,t_{i})$, $C_{i}$ is a constant such that $u_{i}(x)$ satisfies $\frac{1}{V}\int_{M}u_{i}^{2}dV_{\varepsilon_{i}t_{i}}=1$.
$$1=\frac{1}{V}\int_{M}e^{2C_{i}}\phi_{i}^{2}dV_{\varepsilon_{i}t_{i}}\leq\frac{1}{V}e^{2C_{i}}Vol(B_{\varepsilon_{i} t_{i}}(k_{1}^{i},k_{2}^{i}))\leq\frac{1}{V}e^{2C_{i}}\delta_{i}.$$
Let $i\rightarrow+\infty$, since $\delta_{i}\rightarrow0$, we conclude that $C_{i}\rightarrow+\infty$. We Consider the function $\frac{1}{2}(u_{i}^{2}+1)$ whose integral average is

$$\frac{1}{V}\int_{M}\frac{1}{2}(u_{i}^{2}+1)dV_{\varepsilon_{i}t_{i}}=1.$$ Computing the $\mathcal{W}_{\theta_{\varepsilon_{i}}}(g_{\varepsilon_{i}}(t_{i}),-\log \frac{1}{2}(u_{i}^{2}+1),1)$ functional, we have
\begin{eqnarray*}
&\ &\mathcal{W}_{\theta_{\varepsilon_{i}}}(g_{\varepsilon_{i}}(t_{i}),-\log \frac{1}{2}(u_{i}^{2}+1),1)\\
&=&\frac{1}{2}\int_{M}(u_{i}^{2}+1)(R(g_{\varepsilon_{i}}(t_{i}))
-tr_{g_{\varepsilon_{i}}(t_{i})}\theta_{\varepsilon_{i}}+L)dV_{\varepsilon_{i}t_{i}}\\
&\ &+\frac{1}{2}\int_{M}(u_{i}^{2}+1)\frac{4u_{i}^{2}|\nabla u_{i}|_{g_{\varepsilon_{i}}(t_{i})}^{2}}{(u_{i}^{2}+1)^{2}}dV_{\varepsilon_{i}t_{i}}-LV\\
&\ &+\frac{\beta}{2}\log2\int_{M}(u_{i}^{2}+1)dV_{\varepsilon_{i}t_{i}}-\frac{\beta}{2}\int_{M}(u_{i}^{2}+1)
\log(u_{i}^{2}+1)dV_{\varepsilon_{i}t_{i}}.
\end{eqnarray*}
where $L$ satisfies $R(g_{\varepsilon_{i}}(t_{i}))
-tr_{g_{\varepsilon_{i}}(t_{i})}\theta_{\varepsilon_{i}}+L\geq0$ for every $i$ uniformly.
\begin{eqnarray*}
&\ &\frac{1}{2}\int_{M}(u_{i}^{2}+1)(R(g_{\varepsilon_{i}}(t_{i}))
-tr_{g_{\varepsilon_{i}}(t_{i})}\theta_{\varepsilon_{i}}+L)dV_{\varepsilon_{i}t_{i}}\\
&\leq&\frac{1}{2}\int_{B_{\varepsilon_{i}t_{i}}(r_{1}^{i},r_{2}^{i})}e^{2C_{i}}(R(g_{\varepsilon_{i}}(t_{i}))
-tr_{g_{\varepsilon_{i}}(t_{i})}\theta_{\varepsilon_{i}}+L)dV_{\varepsilon_{i}t_{i}}\\
&+&\frac{1}{2}\int_{M}(\beta n-\Delta _{g_{\varepsilon_{i}}(t_{i})}u_{\varepsilon_{i}}(t_{i})+L)dV_{\varepsilon_{i}t_{i}}\\
&\leq&C e^{2C_{i}}Vol_{g_{\varepsilon_{i}}(t_{i})}(B_{\varepsilon_{i}t_{i}}(k_{1}^{i},k_{2}^{i}))+\frac{1}{2}(\beta n+L)V,\\
&\ &\frac{1}{2}\int_{M}(u_{i}^{2}+1)\frac{4u_{i}^{2}|\nabla u_{i}|_{g_{\varepsilon_{i}}(t_{i})}^{2}}{(u_{i}^{2}+1)^{2}}dV_{\varepsilon_{i}t_{i}}
\leq2n\int_{M}e^{2C_{i}}|\phi_{i}'|^{2}dV_{\varepsilon_{i}t_{i}}\\
&\leq&C e^{2C_{i}}Vol_{g_{\varepsilon_{i}}(t_{i})}(B_{\varepsilon_{i}t_{i}}(k_{1}^{i},k_{2}^{i})),\\
&\ &-\frac{\beta}{2}\int_{M}(u_{i}^{2}+1)\log(u_{i}^{2}+1)dV_{\varepsilon_{i}t_{i}}\\
&\leq&-\frac{\beta}{2}\int_{M}u_{i}^{2}\log u_{i}^{2}dV_{\varepsilon_{i}t_{i}}\\
&=&-\beta C_{i}\int_{M}u_{i}^{2}dV_{\varepsilon_{i}t_{i}}-\frac{\beta}{2}\int_{M}e^{2C_{i}}\phi_{i}^{2}\log\phi_{i}^{2}dV_{\varepsilon_{i}t_{i}}\\
&\leq&-\beta VC_{i}+C e^{2C_{i}}Vol_{g_{\varepsilon_{i}}(t_{i})}(B_{\varepsilon_{i}t_{i}}(k_{1}^{i},k_{2}^{i})),
\end{eqnarray*}
where constants $C$ are uniform. Combining all these inequalities together, we have
\begin{eqnarray*}\mathcal{W}_{\theta_{\varepsilon_{k}}}(g_{\varepsilon_{k}}(t_{k}),-\log \frac{1}{2}(u_{i}^{2}+1),1)&\leq&C-\beta V C_{i}+C e^{2C_{i}}Vol_{g_{\varepsilon_{i}}(t_{i})}(B_{\varepsilon_{i}t_{i}}(k_{1}^{i},k_{2}^{i}))\\
&\leq&C-\beta V C_{i}+C2^{20n} e^{2C_{i}}Vol_{g_{\varepsilon_{i}}(t_{i})}(B_{\varepsilon_{i}t_{i}}(k_{1}^{i}+2,k_{2}^{i}-2))\\
&=&C-C_{i}\beta V+C2^{20n}\int_{B_{\varepsilon_{i}t_{i}}(k_{1}^{i}+2,k_{2}^{i}-2)}e^{2C_{i}}\phi_{i}^{2}dV_{\varepsilon_{i}t_{i}}\\
&\leq&C-C_{i}\beta V,
\end{eqnarray*}
where $C$ are constants independent of time $t_{i}$ and $\varepsilon_{i}$. Hence, by $(\ref{3.22.5})$, we have
\begin{align*}\label{3.22.30}-C\leq\mu_{\theta_{\varepsilon_{i}}}(g_{\varepsilon_{i}}(1),1)\leq C-\beta VC_{i},\tag{A.17}\end{align*}
where $C$ are positive constants independent of $t_{i}$ and $\varepsilon_{i}$. Let $i\rightarrow+\infty$, we have
$$-C\leq-\infty,$$
which is impossible. Hence, $diam(M,g_{\varepsilon}(t))$ is uniformly bounded when $t\geq1$.\QEDB

\hspace{1.4cm}

\end{document}